\def\Datum{October 10, 2006}
\magnification=\magstephalf
\hsize=16.2truecm
\vsize24truecm
\parskip4pt plus 1pt
\frenchspacing
\lineskiplimit=-2pt
\vglue6truemm

\def \Eins {1}
\def \Zwei {2}

\def \Drei {3}

\def \Vier {4}

\def \Fuenf {5}

\def \Sechs {6}

\def \UH {6.5}

\def \Sieben {7}

\def \Acht {8}

\def \Neun {9}

\def \Zehn {10}

\def \Elf {11}
\def \Zwoelf {12}

\def \Dreizehn {13}

\def \Vierzehn {14}

\def \Fuenfzehn {15}

\def \Sechzehn {16}

\input cls.def\immediate\openout\aux=\jobname.aux

\arxtrue

\ifdraft
\footline={\hss\sevenrm \Datum\hss}
\else\nopagenumbers\fi

\headline={\ifnum\pageno>1\sevenrm\ifodd\pageno Homogeneous Levi
degenerate CR-manifolds
\hss{\tenbf\folio}\else{\tenbf\folio}\hss{Fels-Kaup}
\fi\else\hss\fi}

\centerline{\Gross Classification of  Levi degenerate homogeneous}
\bigskip
\centerline{\Gross CR-manifolds in dimension 5}
\bigskip\bigskip
\centerline{Gregor Fels\quad and\quad Wilhelm Kaup}
\bigskip
{\parindent0pt\footnote{}{\ninerm 2000 Mathematics
Subject Classification: 32M17, 32M25, 32V25.
}}

\def\Partskip{\vskip9mm}

\bigskip\bigskip

\ifarx
{\narrower\narrower\Klein{\bf\noindent Abstract:} In this paper we present new
examples of homogeneous $2$-nondegenerate CR-manifolds of dimension 5
and give, up to local CR-equivalence, a full classification of all
CR-manifolds of this type.\par}

\bigskip\bigskip

{\def \Einsa {1}
\def \Zweia {6}
\def \Dreia {7}
\def \Viera {11}
\def \Fuenfa {14}
\def \Sechsa {16}
\def \Siebena {22}
\def \Achta {26}
\def \Neuna {29}
\def \Zehna {33}
\def \Elfa {36}
\def \Zwoelfa {38}
\def \Dreizehna {41}
\def \Vierzehna {43}
\def \Fuenfzehna {44}
\def \Sechzehna {49}
\def \Refa {52}
 \SkIp=70pt
\advance\baselineskip8pt\Klein\advance\hsize-60pt
\line{\hbox to \SkIp{\hfill\Eins. }Introduction\dotfill\Einsa}
\vskip8pt
\line{\hbox to 90pt{\hfill}\elfrm PART 1: Levi degenerate CR-manifolds\hfill}
\vskip1pt
\line{\hbox to \SkIp{\hfill\Zwei. }Preliminaries\dotfill\Zweia}
\line{\hbox to \SkIp{\hfill\Drei. }Tube manifolds\dotfill\Dreia}
\line{\hbox to \SkIp{\hfill\Vier. }Tube manifolds over cones\dotfill\Viera}
\line{\hbox to \SkIp{\hfill\Fuenf. }Some examples\dotfill\Fuenfa}
\line{\hbox to \SkIp{\hfill\Sechs. }Levi degenerate CR-manifolds associated with an
endomorphism\dotfill\Sechsa}
\line{\hbox to \SkIp{\hfill\Sieben. }Homogeneous 2-nondegenerate manifolds of CR-dimension 2\dotfill\Siebena}
\line{\hbox to \SkIp{\hfill\Acht. }Homogeneous 2-nondegenerate CR-manifolds in dimension 5\dotfill\Achta}
\vskip8pt
\line{\hbox to 90pt{\hfill}\elfrm PART 2: The classification\hfill}
\vskip1pt
\line{\hbox to \SkIp{\hfill\Neun. }Lie-theoretic characterization of locally homogeneous
CR-manifolds\dotfill\Neuna}
\line{\hbox to \SkIp{\hfill\Zehn. }The Lie algebra $\7g$ has small semisimple part\dotfill\Zehna}
\line{\hbox to \SkIp{\hfill\Elf. }The cases $\7s\cong\so(2,3)$ and $\7s\cong\so(1,3)$\dotfill\Elfa}
\line{\hbox to \SkIp{\hfill\Zwoelf. }The case $\7s\cong\7{sl}(2,\RR)$\dotfill\Zwoelfa}
\line{\hbox to \SkIp{\hfill\Dreizehn. }The case $\7s\cong\su(2)$\dotfill\Dreizehna}
\line{\hbox to \SkIp{\hfill\Vierzehn. }Reduction to the case where $\7g$ is solvable and of dimension
5\dotfill\Vierzehna}
\line{\hbox to \SkIp{\hfill\Fuenfzehn. }The existence of a 3-dimensional abelian ideal in
$\7g$ suffices\dotfill\Fuenfzehna}
\line{\hbox to \SkIp{\hfill\Sechzehn. }The final steps\dotfill\Sechzehna}
\vskip4pt
\line{\hbox to \SkIp{\hfill}References\dotfill\Refa}
}\fi

\newwrite\lst\immediate\openout\lst=\jobname.lst

\KAP{Eins}{Introduction}

The main topic of this paper is the study of real-analytic
CR-manifolds $M$ with everywhere degenerate Levi form. In particular,
for homogeneous manifolds of this type, we develop methods for the
computation of the Lie algebra s$\hol(M,a)$ of infinitesimal
CR-transformations at every $a\in M$. We also classify up to local
CR-equivalence all locally homogeneous degenerate
CR-manifolds in dimension 5.

In this context, a well studied example of a homogeneous Levi degenerate
CR-manifold is the quadratic hypersurface
$$
\5M:=\{z\in\CC^{3}: (\Re z_{1})^{2}+(\Re z_{2})^{2}=(\Re
z_{3})^{2},\;\Re z_{3} >0\} \;,
$$ compare
e.g. \Lit{EBFT}, \Lit{FEKA}, \Lit{KAZT} and \Lit{SEVL}.  This
5-dimensional CR-manifold has several remarkable properties and
serves as motivation for various considerations in this paper.
Notice that $\5M$ can also be written as tube manifold
$$\5M=\5F+i\RR^{3}\subset\CC^{3},\steil{where}\5F:=\{x\in\RR^{3}:
x^{2}_{1}+x^{2}_{2}=x^{2}_{3},\,x_{3}>0\}$$ is the future light cone
in 3-dimensional space-time. A glance at this description shows that
$\5M$ is homogeneous under a group of complex-affine
transformations. It is less obvious that the Lie algebras of {\sl
global} and {\sl local} infinitesimal CR-transformations at $a\in\5M$,
$\hol(\5M)$ and $\hol(\5M,a)$, are both isomorphic to $\so(2,3)$, and
hence have dimension 10, compare \Lit{KAZT}. Also the following
`globalization' is known: The group $\SO(2,3)$ acts on the complex
quadric $\fett\3Q_3\subset\PP_{4}(\CC)$ by biholomorphic
transformations and has a hypersurface orbit that contains $\5M$ as a
dense domain.

\ifarx\medskip\fi

The cone $\5F$ clearly is a disjoint union of affine
half-lines. Therefore, $\5M$ is a disjoint union of complex
half-planes, actually $\5M$ is a fiber bundle with typical fiber
$\HH^{+}:=\{z\in\CC:\Re(z)>0\}$. However, the total space $\5M$ is not
even locally CR-equivalent to a product $\HH^+\Times M'$ with $M'$ a
CR-manifold.  Notice that the Levi form in both cases, i.e., for $\5M$
and for a product of $\HH^{+}$ with a Levi nondegenerate 3-dimensional
CR-manifold $M'$, has exactly one non-zero eigenvalue at every
point. Hence, one needs further invariants to distinguish those
CR-manifolds.  While every product $M'\times\CC$ is holomorphically
degenerate, the crucial fact here is that the light cone tube $\5M$ is
nondegenerate in a higher order sense: To be precise, $\5M$ is
2-nondegenerate at every point, and we refer to \Lit{BAHR} and also to
\Lit{BERO}, 11.1, for the notion of $k$-nondegeneracy.

\ifarx\medskip\fi

In the non-homogeneous setting, for every $k\in \NN$ and fixed
manifold dimension it is not difficult to construct large classes of
CR-manifolds, even hypersurfaces, which are $k$-nondegenerate at some
points, but are Levi nondegenerate in a dense open subset.  It seems
to be much harder to construct CR-manifolds which are
$k$-nondegenerate everywhere for $k\ge 2$.  Note that the CR-dimension
of a {\sl homogeneous} $M$ is an upper bound for the degree $k$ of
nondegeneracy.  Hence, the lowest manifold dimension for which
everywhere $2$-nondegenerate CR-manifolds can exist is 5.  which
raises the intriguing question whether besides the light cone tube
there exist further 2-nondegenerate homogeneous CR-manifolds in
dimension 5. So far, all known examples in dimension finally turned
out to be locally CR-isomorphic to $\5M$, compare e.g. \Lit{EBEN},
\Lit{KAZT}, \Lit{GAME}, \Lit{FEKA}, \Lit{BELO}, and the belief arose
that there are no others.

\ifarx\medskip\fi

The main objective of this paper is to show that there are actually
infinitely many locally mutually non-equivalent examples and to
provide a full classification. Starting point is the following simple
observation: Suppose that $F\subset\RR^{n}$ is an affinely homogeneous
(connected) submanifold of dimension say $d<n$. Then the tube
$M:=F+i\RR^{n}$ is a generic CR-submanifold of $\CC^{n}$ of
CR-dimension $d$ and is homogeneous under a group of complex-affine
transformations. Indeed, every real-affine transformation leaving $F$
invariant extends to a complex-affine transformation leaving $M$
invariant and, furthermore, $M$ is invariant under all translations
$z\mapsto z+iv$ with $v\in\RR^{n}$. Clearly, the crucial question is,
when is $M$ $k$-nondegenerate and when are two tubes $M,M'$ of this
type locally CR-equivalent?

\ifarx\medskip\fi

The classification of all affinely homogeneous surfaces
$F\subset\RR^{3}$ can be found in \Lit{DOKR} and \Lit{EAEZ}. In
particular, a complete list (up to local affine equivalence and in a
slightly different formulation) of all degenerate types that are not a
cylinder, is given by the following examples (1) -- (3).

\item{(1)} $F=\5F$ the future light cone as above.

\item{(2a)} $F=\{\,r\,(\cos t,\sin t,e^{\omega
t})\in\RR^{3}:r\in\RRp,t\in\RR\}$ with $\omega>0$ arbitrary.

\item{(2b)} $F=\{\,r\,(1,t,e^{t})\in\RR^{3}:r\in\RRp,t\in\RR\}$.

\item{(2c)} $F=\{\,r\,(1,e^{t},e^{\theta
t})\in\RR^{3}:r\in\RRp,t\in\RR\}$ with $\theta>2$ arbitrary.

\item{(3)} $F=\{\,c(t)+r\re1c'(t)\in\RR^{3}: r\in\RRp,t\in\RR\}$,
where $c(t):=(t,t^{2},t^{3})$ parameterizes the {\sl twisted cubic}
$\{(t,t^{2},t^{3}):t\in\RR\}$ in $\RR^{3}$ and $c'(t)=(1,2t,3t^{2})$.

\noindent Notice that the limit case $\omega=0$ in (2a) gives the
future light cone $\5F$, while the limit case $\theta=2$ in (2c) gives
the linearly homogeneous surface
$\{x\in\RR^{3}:x_{1}x_{3}=x_{2}^{2},\,x_{1}>0,\,x_{2}>0\}$,
which is locally, but not globally linearly equivalent to $\5F$. In
fact, $\5F$ is linearly equivalent to the cone\nline
$\{x\in\RR^{3}:x_{1}x_{3}=x_{2}^{2},\,x_{1}+x_{3}>0\}$.

\ifarx\vfill\eject\else\medskip\fi

As our first main result we show, compare \ruf{UA} and \ruf{UB} for
details:

\noindent{\parskip0pt\bf Theorem I: \sl For every surface $F$ in {\rm(1) -- (3)}
the corresponding tube manifold $M:=F+i\RR^{3}$ is a homogeneous
$2$-nondegenerate CR-submanifold of $\CC^{3}$ and any two of them are
pairwise locally CR-inequivalent. Furthermore, for every $F$ in
{\rm(2a) -- (3)} and every $a\in M=F+i\RR^{3}$ the following holds: \0
The Lie algebra $\hol(M,a)$ is solvable and has dimension $5$. \1 The
stability group $\Aut(M,a)$ is trivial. \1 Every homogeneous
real-analytic CR-manifold $M'$, that is locally CR-equivalent to $M$,
is already globally CR-equivalent to $M$.\par}

\smallskip\noindent Notice that, a priori, there is no reason why the
$F$ in (1) -- (3), although known to be locally affinely inequivalent,
should have locally CR-inequivalent tubes (for nondegenerate affinely
homogeneous surfaces in $\RR^{3}$, for instance, there are
counterexamples).

\medskip We actually prove an analog of Theorem I in every dimension
$n\ge3$, where the same trichotomy occurs as above: Consider the
following surfaces $F\subset\RR^{n}$

\item{(1')} $F=\5F^{n}:=\{x\in\RR^{n}:x_{1}>0,\,x_{2}>0\steil{and}
x_{j}=x^{j-1}_{2}/x_{1}^{j-2} \steil{for}3\le j\le n\}$.

\item{(2')}
$F=\{r\re1e^{t\phi}(a)\in\RR^{n}:r\in\RRp,t\in\RR\}$, where
$\phi$ is an endomorphism of $\RR^{n}$ having $a\in\RR^{n}$ as cyclic
vector (i.e., the iterates $\phi^{k}(a)$, $k\ge0$ span $\RR^{n}$) and
the $n$ eigenvalues of $\phi$ do not form an arithmetic progression in
$\CC$.

\item{(3')} $F=\{c(t)+r\re1c'(t):r\in\RRp,t\in\RR\}$, where
$c(t):=(t,t^{2},\dots,t^{n})\in\RR^{n}$ parameterizes the twisted
$n$-ic in $\RR^{n}$.\par

\noindent In Sections \ruf{Sechs} and \ruf{Sieben} we show among other
statements: For every $F,F'$ from (1') -- (3') the tube manifolds
$M:=F+i\RR^{n}$, $M':=F'+i\RR^{n}$ are affinely homogeneous generic
$2$-nondegenerate submanifolds of $\CC^{n}$ with CR-dimension $2$.
Furthermore, $M$, $M'$ are locally CR-equivalent if and only if $F$,
$F'$ are globally affinely equivalent, and this holds if and only if
for given $a\in M$, $a'\in M'$ the Lie algebras $\hol(M,a)$,
$\hol(M',a')$ are isomorphic.  In case $F=\5F^{n}$ the Lie algebra
$\hol(M,a)$ contains a copy of $\gl(2,\RR)$ and hence is not
solvable.  In all other cases, i.e. $F$ comes from (2') or (3'), the
Lie algebra $\hol(M,a)$ is solvable of dimension $n+2$ and the
stability group $\Aut(M,a)$ has order at most $2$.

\medskip Let us briefly comment on the proof of Theorem I: Once
defining equations for an $F$ under consideration are explicitly known
(this is quite obvious for the types (1) -- (3), compare Section
\ruf{Acht}, but seems to be hard for the types (2'), (3')\hskip1pt),
one can compute by standard methods the order $k$ of
nondegeneracy. However, the amount of calculation necessary to
determine $k$ in such a way grows very fast with $k$ and with the
dimension or codimension of $M\subset \CC^n$. This is one of the
reasons, especially with an eye on possible generalizations, to choose
a different approach, which does not use explicit equations. For
instance, given an arbitrary submanifold $F\subset \RR^n$ which is
(locally) {\sl affinely homogeneous}, we present a simple criterion
(Proposition \ruf{US}) which allows to determine quickly the order $k$
of nondegeneracy for the corresponding tube manifold.

\noindent The hard part of the proof is to show that the various tubes
$F+i\RR^3$ are actually locally inequivalent as CR-manifolds.  Recall
that for real-analytic, not necessary homogeneous, {\sl hypersurfaces
with nondegenerate Levi form} there exist local invariants which
determine each $M$ up to local CR-equivalence due to the
fundamental work of Cartan, Tanaka and Chern-Moser.  However, an
analogues approach is not available for $M$ of higher codimensions or
when the Levi form is degenerate. To distinguish the various tubes
$F+i\RR^3,$ we develop a method (valid also in greater generality)
which enables us to determine explicitly the Lie algebras $\hol(M,a)$
of infinitesimal CR-transformations of the various CR-germs $(M,a)$
(see Section \ruf{Zwei} for basic definitions and Sections \ruf{Sechs}
for further details).

\medskip The CR-manifolds occurring in the previous theorem are quite
special as they all are tube manifolds.  Moreover, all but one (namely
the twisted cubic case in (3)\hskip1pt) are actually conical.  In the
Levi nondegenerate case, many (homogeneous) examples are known which
are not locally CR-equivalent to any tube manifold. For instance, the
unit sphere subbundle of $TS^3$ with its canonical CR-structure is
such an example. Therefore, our second main result came quite
unexpected to us:

\medskip\noindent {\bf Theorem II.} {\sl Every 5-dimensional locally
homogeneous 2-nondegenerate CR-manifold $M$ is locally
CR-equivalent to $F+i\RR^3$ with $F$ being one of the surfaces in {\rm
(1) -- (3)}}.

\medskip\noindent For the precise definition of local homogeneity we
refer to Section \ruf{Zwei}.  A priori, locally homogeneous
CR-manifolds might exist which are not locally CR-equivalent to any
globally homogeneous one.  As a by-product of the above 2 results we
can see that such a pathology does not happen in the case under
consideration. Theorem II gives a classification of all
(abstract) 2-nondegenerate locally homogeneous CR-manifolds in
dimension 5 up to local CR-equivalence.  In fact, using Cartan's
classification \Lit{CART} of the 3-dimensional Levi nondegenerate
homogeneous CR-manifolds, the following result holds.

\medskip\noindent {\bf Classification.}  {\sl Every 5-dimensional
locally homogeneous CR-manifold $M$ with degenerate Levi form is
locally CR-equivalent to one of the following: \0 $M=F+i\RR^3\subset
\CC^3$ and $F$ is one of the surfaces {\rm (1) -- (3)} from above.  \1
$M=\CC\times M'$, where $M'$ is one of the 3-dimensional Levi
nondegenerate CR-manifolds from Cartan's list in \Lit{CART}. \1
$M=\CC^2\times \RR$ or $M=\CC\times \RR^3$ with $\RR^3$ totally real.

\noindent The manifolds in {\rm (i)} are all $2$-nondegenerate and in
{\rm (ii) -- (iii)} are holomorphically degenerate. Also, the
manifolds in {\rm (iii)} are just the Levi flat ones.}

\smallskip With Theorem II the following question arises naturally for
higher codimension: {\sl Are there, up to local CR-equivalence,
further locally homogeneous $2$-nondegenerate CR-manifolds of
CR-dimension 2 besides those that are tubes over the surfaces $F$ in
{\rm (1') -- (3')} above?} Notice in this context that every locally
homogeneous $2$-nondegenerate CR-manifold of dimension $6$ necessarily
has CR-dimension $2$.

\smallskip For 5-dimensional CR-hypersurfaces with {\sl nondegenerate}
Levi form, that is, when Chern-Moser invariants are available, there
already exists a partial classification: Locally homogeneous
CR-hypersurfaces in $\CC^3$ with stability groups of positive
dimension have been classified by Loboda in terms of local equations
in normal form, see \Lit{LOBO}, \Lit{LOBP}, \Lit{LOBQ}.

\medskip The paper is organized as follows. After recalling some
necessary preliminaries in Section \ruf{Zwei} we discuss in Section
\ruf{Drei} tube manifolds $M=F\oplus i\RR^{n}\subset\CC^{n}$ over
real-analytic submanifolds $F\subset\RR^{n}$. It turns out that the
CR-structure of $M$ is closely related to the real-affine structure of
the base $F$. For instance, the Levi form of $M$ is essentially the
sesquilinear extension of the second fundamental form of the
submanifold $F\subset\RR^{n}$. Generalizing the notion of the second
fundamental form we define higher order invariants for $F$ (see
\ruf{HN}). In the uniform case these invariants precisely detect the
$k$-nondegeneracy of the corresponding CR-manifold $M=F\oplus i\RR^n$.
It is known that the (uniform) $k$-nondegeneracy of a real-analytic
CR-manifold $M$ together with minimality ensures that the Lie algebras
$\hol(M,a)$ are finite-dimensional, and is equivalent to this in the
special class of locally homogeneous CR-manifolds. For submanifolds
$F\subset\RR^{n}$ which in addition are homogeneous under a group of
affine transformations, a simple criterion for $k$-nondegeneracy of
the associated tube manifold $M$ is given in Proposition \ruf{US}.

In Section \ruf{Vier} these results are applied to the case where $F$
is conical in $\RR^n$, that is, locally invariant under dilations
$z\mapsto tz$ for $t$ near $1\in\RR$. In this case $M$ is always {\sl
Levi degenerate.} Assuming that $\hol(M,a)$ is finite dimensional
(which automatically is the case for minimal and finitely
nondegenerate CR-manifolds) we develop some basic techniques for the
explicit computation of $\hol(M,a)$.  The main results in Section
\ruf{Vier} are the following: We prove that under the finiteness
assumption $\hol(M,a)$ consists only of polynomial vector fields and
carries a natural graded structure, see Proposition \ruf{DC}. We prove
(under the same assumptions, see Proposition \ruf{ZB}.ii) that local
CR-equivalences between two such tube manifolds are always rational
maps (even if these manifolds are not real-algebraic). Furthermore,
jet determination estimates are provided (\ruf{ZB}.iv). In the special
situation where $\hol(M,a)$ consists of affine germs only, these
results are further strengthened.

In Section \ruf{Fuenf} we illustrate by examples how our methods can
be applied.  In Example \ruf{JP} we present to every $c\ge 1$ and
$k\in \{2,3,4\}$ a homogeneous submanifold of $\CC^{n}$ which is
$k$-nondegenerate, of codimension $c$ and of CR-dimension $k$.  We
close this section investigating the tubes $M^\alpha_{p,q}$ over the
cones $F^\alpha_{p,q}:=\{x\in \RR^{p{+}q}_+: \sum \epsilon_j
x_j^\alpha=0\}$ with $p$ positive and $q$ negative $\epsilon_j$'s and
$\alpha\in \RR^{*}$. Using our methods from the previous section we
explicitly determine the Lie algebras $\hol(M^\alpha_{p,q},a)$ for
arbitrary integers $\alpha\ge2$ and $p\ge q\ge 1$.

In Section \ruf{Sechs}, we construct homogeneous CR-submanifolds
$M=M^{\phi,d}\subset\CC^{n}$ of tube type, depending on the choice of
an endomorphism $\phi\in \End(\RR^n)$, an integer $1<d<n$ and a cyclic
vector $a\in\RR^{n}$ for $\phi$ in the following way: The powers
$\phi^0,\phi^1,\ldots,\phi^{d-1}$ span an abelian Lie algebra $\7h$
and, in turn, give a cone $F:=\exp(\7h)a\subset \RR^n$. For $\phi$
which admits a cyclic vector $a\in \RR^n$ the corresponding tube
manifold $M^{\phi,d}=F+i\RR^{n}$ is 2-nondegenerate, minimal and of
CR-dimension $d$.  The key result here is the explicit determination
of the CR-invariant $\hol(M^{\phi,d},a)$ for $\phi$ in 'general
position' (Propositions \ruf{UH} and \ruf{HU}). The precise meaning
for $\phi$ of being in 'general position' is stated in Lemma \ruf{UL}.
Further results are, again for $\phi$ in general position, that the
tube manifold $M$ is simply connected and has trivial stability group
at every point. As a consequence, the manifolds $M$ of this type have
the following remarkable property: {\sl Every homogeneous
(real-analytic) CR-manifold locally CR-equivalent to $M$ is globally
CR-equivalent to $M$}.

In Section \ruf{Sieben} the results from Section \ruf{Sechs} are
further refined in the case of homogeneous CR-manifolds
$M^\phi:=M^{\phi,2}:=F^\phi+i\RR^n\subset \CC^n$ of CR-dimension 2 but
without restrictions on the codimension.  In fact, every minimal and
locally homogeneous tube CR-manifold $M:=S+i\RR^n$ with a conical
2-dimensional $S\subset \RR^n$ is locally CR-equivalent to
$M^{\phi,2}$ for some cyclic $\phi$.  In this section also the case
when $\phi$ is not in 'general position', i.e., the characteristic
roots do form an arithmetic progression (see Lemma \ruf{UL})
is treated.  The main results here are: Whether a cyclic $\phi$ is in
general position or not, the Lie algebras $\hol(M^\phi,a)$
(Proposition \ruf{ZR}) and the global automorphism groups
$\Aut(M^\phi)$ (Proposition \ruf{GO}) are determined.  Furthermore,
the problem of local and global CR-equivalence among the $M^\phi$'s is
solved (Propositions \ruf{KO} and \ruf{BP}) and a moduli space is
constructed (Subsection \ruf{MO}).

Part I of the paper is concluded with Section \ruf{Acht} where the
examples (1) -- (3) from Theorem I are presented in more detail.  The
results of the preceeding section are applied to this case of
5-dimensional tube manifolds $M$. In particular Propositions \ruf{UA}
and \ruf{IZ} contain some additional information to that stated in
Theorem I and also complete the proof of Theorem I.

Part II of the paper is mainly devoted to prove Theorem II.  In the
preliminary section \ruf{Neun} we explain how the geometric properties
such as $k$-nondegeneracy, minimality or the CR-dimension of an
arbitrary locally homogeneous CR-germ $(M,o)$ can be encoded in a pure
Lie algebraic object, the associated CR-algebra $(\7g,\7q)$. This is
the key for our classification and is based on results taken from
\Lit{FELS}. Specified to 5-dimensional CR-germs, we formulate the
precise algebraic conditions on a CR-algebra $(\7g,\7q)$ ensuring that
the associated CR-germ $(M,o)$ is 2-nondegenerate.

Once the classification of 2-nondegenerate 5-dimensional locally
homogeneous CR-germs $(M,o)$ is reduced to a classification problem of
certain CR-algebras, we begin in Section \ruf{Zehn} to carry out the
details of the proof.  It is subdivided into several sections, lemmata
and claims and will only be completed in Section \ruf{Sechzehn}.  Our
proof relies on a quite subtle analysis of Lie algebraic properties of
CR-algebras and uses basic structure theory of Lie algebras and Lie
groups.  The methods are quite general and can be adapted to handle
also higher dimensional cases.

A more detailed outline of the proof can be found in the first part of
Section \ruf{Zehn}.

\ifarx\vfill\eject\else\Partskip\fi

\centerline{\twelverm PART 1: Levi degenerate CR-manifolds}
\bigskip

\KAP{Zwei}{Preliminaries}

In the following let $E$ always be a complex vector space of finite
dimension and $M\subset E$ an immersed connected real-analytic
submanifold. In most cases $M$ will be locally closed in $E$. Due to
the canonical identifications $T_aE=E,$ for every $a\in M$ we consider
the tangent space $T_{a}M$ as an $\RR$-linear subspace of $E$. Then
$H_{a}M:=T_{a}M\cap iT_{a}M$ is the largest complex linear subspace of
$E$ contained in $T_{a}M$. The manifold $M$ is called a {\sl
CR-submanifold} if $\dim_{\CC}\!H_{a}M$ does not depend on $a\in
M$. This dimension is called the {\sl CR-dimension} of $M$ and
$H_{a}M$ is called the {\sl holomorphic tangent space} at $a$, compare
\Lit{BERO} as general reference for CR-manifolds. Given a further
real-analytic CR-submanifold $M'$ of a complex vector space $E'$, a
smooth mapping $g:M\to M'$ is called {\sl CR} if for every $a\in M$
the differential $dg_{a}:T_{a}M\to T_{ga}M'$ maps the corresponding
holomorphic tangent spaces in a complex-linear way to each
other. Keeping in mind the identification $T_aE=E$, a vector field on
$M$ is a mapping $f:M\to E$ with $f(a)\in T_{a}M$ for all $a\in
M$. For better distinction we also write $\xi=f(z)\dd z$ instead of
$f$ and $\xi_{a}$ instead of $f(a)$, compare the convention (2.1) in
\Lit{FEKA}.

An {\sl infinitesimal CR-transformation} on $M$ is by definition a
real-analytic vector field $f(z)\dd z$ on $M$ such that the
corresponding local flow consists of CR-transformations. Let us denote
by $\hol(M)$ the space of all such vector fields, which is a real Lie
algebra with respect to the usual bracket. For every $f(z)\dd
z\in\hol(M)$ and every $a\in M$ there exist an open neighbourhood
$U\subset M$ of $a$ with respect to the manifold topology on $M$, an
open neighbourhood $W$ of $a$ with respect to $E$ and
a holomorphic mapping $h:W\to E$ with $f(z)=h(z)$
for all $z\in U\cap W$, compare \Lit{ANHI} or 12.4.22 in
\Lit{BERO}.

Further, for every $a\in M$ we denote by $\hol(M,a)$ the Lie algebra
of all germs of infinitesimal CR-transformations defined on arbitrary
open neighbourhoods of $a$ with respect to the manifold topology of
$M$. For simplicity and without essential loss of generality we always
assume that the CR-submanifold $M$ is {\sl generic} in $E$, that is,
$E=T_{a}M\oplus iT_{a}M$ for all $a\in M$. This assumption allows us to
consider $\hol(M,a)$ in a canonical way as a real Lie subalgebra of
the complex Lie algebra $\hol(E,a)$, what we always do in the
following.  The CR-manifold $M$ is called {\sl holomorphically
nondegenerate} at $a$ if $\hol(M,a)$ is totally real in $\hol(E,a)$,
that is, if $\hol(M,a)\cap i\,\hol(M,a)=0$ in $\hol(E,a)$. This
condition together with the minimality of $M$ at $a$ implies $\dim
\hol(M,a)<\infty$, see 12.5.3 in \Lit{BERO}. Here, the CR-submanifold
$M\subset E$ is called {\sl minimal} at $a\in M$ if $T_{a}R=T_{a}M$
for every submanifold $R\subset M$ with $a\in R$ and $H_{z}M\subset
T_{z}R$ for all $z\in R$. By Proposition 15.5.1 in \Lit{BERO}, $\dim
\hol(M,a)<\infty$ implies that $M$ is holomorphically nondegenerate at
$a$. The CR-manifold $M$ is called {\sl locally homogeneous} at the
point $a\in M$ if there exists a Lie subalgebra $\7g\subset\hol(M,a)$
of finite dimension such that the canonical evaluation map $\7g\to
T_{a}M$ is surjective. In case $M$ is locally homogeneous at $a$ the
condition $\dim \hol(M,a)<\infty$ is equivalent to $M$ being
holomorphically nondegenerate and minimal at $a$.

By $\aut(M,a):=\{\xi\in\hol(M,a):\xi_{a}=0\}$ we denote the {\sl
isotropy subalgebra} at $a\in M$. Clearly, $\aut(M,a)$ has finite
codimension in $\hol(M,a)$. Furthermore, we denote by $\Aut(M,a)$ the
group of all {\sl germs} of real-analytic CR-isomorphisms $ h:W\to
\tilde W$ with $ h(a)=a$, where $W,\tilde W$ are arbitrary open
neighbourhoods of $a$ in $M$. It is known that every germ in
$\Aut(M,a)$ can be represented by a holomorphic map $U\to E$, where
$U$ is an open neighbourhood of $a$ in $E$, compare e.g. 1.7.13 in
\Lit{BERO}. Furthermore, $\Aut(M)$ denotes the group of all global
real-analytic CR-automorphisms $h:M\to M$ and $\Aut(M)_a$ its isotropy
subgroup at $a$. There is a canonical group monomorphism
$\Aut(M)_a\hookrightarrow\Aut(M,a)$ as well as an exponential map
$\exp:\aut(M,a)\to\Aut(M,a)$ for every $a\in M$.

By $\aff(M)\subset\hol(M)$ we denote the Lie subalgebra of all
(complex) affine infinitesimal CR-transformations on $M$. For every
$a\in M$ furthermore $\aff(M,a)\subset\hol(M,a)$ is the Lie subalgebra
of all affine germs. The canonical embedding
$\aff(M)\hookrightarrow\aff(M,a)$ is an isomorphism for every $a\in
M$.

Suppose that $g\colon U\to U'$ is a biholomorphic mapping between open
subsets $U,U'\subset E$. Then
$$
g_{*}\big(f(z)\dd z\big)=g'(g^{-1}z)\big(f(g^{-1}z)\big)\,\dd
z\Leqno{GP}
$$
defines a complex Lie algebra isomorphism
$g_{*}\colon\hol(U)\to\hol(U')$, where $g'(u)\in\End(E)$ for every
$u\in U$ is the derivative of $g$ at $u$. For real-analytic
CR-submanifolds $M,M'\subset E$ every CR-isomorphism
$g\colon(M,a)\to(M',a')$ of manifold germs induces a Lie algebra
isomorphism $g_{*}\colon\7g\to\7g'$, where $\7g\!:=\hol(M,a)$ and
$\7g'\!:=\hol(M',a')$. From \Ruf{GP} it is clear that $g_{*}$ extends
to a complex Lie algebra isomorphism $\7l\to\7l'$, where the sums
$\7l:=\7g\oplus i\7g\subset\hol(E,a)$ and
$\7l':=\7g'\oplus i\7g'\subset\hol(E,a')$ are not necessarily direct. In
particular, $g\mapsto g_{*}$ defines a group homomorphism
$\Aut(M,a)\to\Aut(\7g)$.

\medskip A basic invariant of a CR-manifold is the (vector valued)
{\sl Levi form}. Its definitions found in the literature may differ by
a constant factor. Here we choose the following one: It is
well-known that for every point $a$ in the CR-manifold $M$ there is a
well defined alternating $\RR$-bilinear map
$$
\omega_{a}:H_{a}M\times H_{a}M\;\to\; E/H_{a}M
$$
satisfying $\omega_{a}(\xi_{a},\zeta_{a})\equiv[\xi,\zeta]_{a}
\steil{mod}H_{a}M,$ where $\xi,\zeta$ are arbitrary smooth vector
fields on $M$ with $\xi_{z},\zeta_{z}\in H_{z}M$ for all $z\in M$. We
define the Levi form
$$
\5L_{a}:H_{a}M\times H_{a}M\;\to\; E/H_{a}M\Leqno{LE}
$$
to be twice the sesquilinear part of $\omega_{a}$. By {\sl
sesquilinear} we always mean `conjugate linear in the first and
complex linear in the second variable', that is,
$$
\5L_{a}(v,w)=\omega_{a}(v,w)+i\omega_{a}(iv,w)\qquad \hbox{for all }
v,w\in H_{a}M\,.
$$
In particular, the vectors $\5L_{a}(v,v),$ $v\in H_aM,$ are contained
in $iT_{a}M/H_{a}M$, which can be identified in a canonical way with
the normal space $E/T_{a}M$ to $M\subset E$ at $a$.  The following
remark follows immediately from the way the Levi form is defined:

\Remark{RM} Suppose $Z$ is a complex manifold, $\phi:Z\to M$ is a smooth
CR-mapping and $a=\phi(c)$ for some $c\in Z$. Then every vector $v\in
d\phi_{c}(T_{c}Z)\subset H_{a}M$ satisfies $\5L_{a}(v,v)=0$. In
general, $v$ is not contained in the {\sl Levi kernel}
$$
K_{a}M:=\{v\in H_{a}M:\5L_{a}(v,w)=0\steil{for all}w\in H_aM\}\,.
$$
\Formend The CR-manifold $M$ is called {\sl Levi nondegenerate} at $a$
if $K_{a}M=0$. Generalizing that, the notion of {\sl
$k$-nondegeneracy} for $M$ at $a$ has been introduced for every
integer $k\ge1$ (see \Lit{BAHR}, \Lit{BERO}). As shown in 11.5.1 of
\Lit{BERO} a real-analytic and connected CR-manifold $M$ is
holomorphically nondegenerate at $a$ (equivalently: at every $z\in M$)
if and only if there exists a $k\ge1$ such that $M$ is
$k$-nondegenerate at some point $b\in M$. For $k=1$ this notion is
equivalent to $M$ being Levi nondegenerate at $b\in M$.

\bigskip In the second part of our paper we also need a more general
notion of a (real-analytic) CR-manifold. This is a connected
real-analytic manifold $M$ together with a subbundle $HM\subset TM$
(called the `holomorphic subbundle') and a bundle endomorphism $J$ of
$HM$ with $J^2= - \id$ such that $(HM,J)$ is involutive, compare 7.4
in \Lit{BOGG}. By a theorem in \Lit{ANFR} there exists an embedding
$M\into Z$ into a complex manifold $Z$, such that $H_zM$ corresponds
to $T_{z}M\cap\,iT_{z}M$, where $T_zZ\to T_zZ$ is simply the
multiplication with the imaginary unit $i$ (here the real-analyticity
is necessary).  The bundle homomorphism $J:HM\to HM$ is then the
restriction of that multiplication with $i$ to $H_{z}M$ for every
$z\in M$.  For local considerations one can always assume that $Z$ is
(an open subset of) $\CC^n$.

\KAP{Drei}{Tube manifolds}

Let $V$ be a real vector space of finite dimension and $E:=V\,\oplus\,
iV$ its complexification. Let furthermore $F\subset V$ be a connected
real-analytic submanifold and $M:=F+iV\subset E$ the
corresponding {\sl tube manifold.} $M$ is a generic CR-submanifold of
$E$, invariant under all translations $z\mapsto z+iv$, $v\in V$. In
case $V'$ is another real vector space of finite dimension, $E'$ its
complexification, $F'\subset V'$ a real-analytic submanifold and
$\phi:V\to V'$ an affine mapping with $\phi(F)\subset F'$, then
clearly $\phi$ extends in a unique way to a complex-affine mapping
$E\to E'$ with $\phi(M)\subset M'$. However, it should be noted that
higher order real-analytic maps $\psi:F\to F'$ also extend locally to
holomorphic maps $\psi:U\to E'$, $U$ open in $E$. But in contrast to
the affine case we have in general $\psi(M\cap U)\not\subset M'$. We
may therefore ask how the CR-structure of $M$ is related to the real
affine structure of the submanifold $F\subset V$.

For every $a\in F$ let $T_{a}F\subset V$ be the {\sl tangent space}
and $N_{a}F:=V/T_{a}F$ the {\sl normal space} to $F$ at $a$. Then
$T_{a}M=T_{a}F\oplus iV$ for the corresponding tube manifold $M$, and
$N_{a}F$ can be canonically identified with the normal space
$N_{a}M=E/T_{a}M$ of $M$ in $E$. Define the map $\ell_{a}:T_{a}F\times
T_{a}F\to N_{a}F$ in the following way: For $v,w\in T_{a}F$ choose a
smooth map $f:V\to V$ with $f(a)=v$ and $f(x)\in T_{x}F$ for all $x\in
F$ (actually it suffices to choose such an $f$ only in a small
neighborhood of $a\,$). Then put
$$
\ell_{a}(v,w):=f'(a)(w)\;\mod\;T_{a}F\,,\Leqno{EL}
$$
where the linear operator $f'(a)\in\End(V)$ is the derivative of $f$
at $a$. One shows that $\ell_a$ does not depend on the choice of $f$
and is a symmetric bilinear map. We mention that if $V$ is provided
with a flat Riemannian metric and $N_aF$ is identified with
$T_aF^\perp$ then $\ell$ is nothing but the second fundamental form of
$F$ (see the subsection II.3.3 in \Lit{SAKA}). The form $\ell_{a}$ can
also be read off from local equations for $F$, more precisely, suppose
that $U\subset V$ is an open subset, $W$ is a real vector space and
$h:U\to W$ is a real-analytic submersion with $F=h^{-1}(0)$. Then the
derivative $h'(a):V\to W$ induces a linear isomorphism $N_{a}F\cong W$
and modulo this identification $\ell_{a}$ is nothing but the second
derivative $h''(a):V\times V\to W$ at $a$, restricted to $T_{a}F\times
T_{a}F$. By
$$
K_{a}F:=\{w\in T_{a}F:\ell_{a}(v,w)=0\steil{for all}v\in T_{a}F\}
$$
we denote the {\sl kernel} of $\ell_{a}$. The manifold $F$ is called
(affinely) {\sl nondegenerate} at $a$ if $K_{a}F=0$ holds. The
following statement follows directly from the definition of $\ell_a$:

\Lemma{CR} Suppose that $\phi\in\End(V)$ satisfies $\phi(x)\in T_{x}F$
for all $x\in F$. Then $\phi(a)\in K_{a}F$ if and only if
$\phi\big(T_{a}F\big)\subset T_{a}F$.\Formend

\noindent Lemma \ruf{CR} applies in particular for $\phi=\id$ in case
$F$ is a {\sl cone,} that is, $rF=F$ for all real $r>0$. More
generally, we call the submanifold $F\subset V$ {\sl conical} if $x\in
T_{x}F$ for all $x\in F$. Then $\RR a\subset K_{a}F$ holds for all
$a\in F$.

\medskip

In the remaining part of this section we explain how the CR-structure
of the tube manifold $M$ is related to the real objects $\ell_a,\,
TF,\, KF$ and $ K^r\!F$, to be defined below, which depend only on the
affine geometry of $F$. In general it needs some effort to check
whether a given CR-manifold $M$ is $k$-nondegenerate at a point $a\in
M$ (in sense of \Lit{BHUR}).  For affinely homogeneous tube manifolds,
however, there are simple criteria, see Propositions \ruf{UR} and
\ruf{US}. We start with some preparations. For every $a\in F\subset M$
$$
H_{a}M=T_{a}F\oplus iT_{a}F\,\subset\, E
\Leqno{HM}
$$
is the holomorphic
tangent space at $a$, and $E/H_{a}M$ can be canonically identified
with $N_{a}F\oplus iN_{a}F$. It is easily seen that the Levi form
$\5L_{a}$ of $M$ at $a$, compare \Ruf{LE}, is nothing but the
sesquilinear extension of the form $\ell_{a}$ from $T_{a}F\times
T_{a}F$ to $H_{a}M\times H_{a}M$. In particular,
$$
K_{a}M=K_{a}F\oplus iK_{a}F
$$
is the {\sl Levi kernel} of $M$ at $a$. In case the dimension of
$K_{a}F$ does not depend on $a\in F$ these spaces form a subbundle
$KF\subset TF$. In that case to every $v\in K_{a}F$ there exists a
smooth function $f:V\to V$ with $f(a)=v$ and $f(x)\in K_{x}F$ for all
$x\in F$, i.e., the tangent vector $v$ extends to a smooth section in
$KF$. In any case, let us define inductively linear subspaces
$K_{a}^{r}F$ of $T_{a}F$ as follows:

\Definition{HN} For every real-analytic submanifold $F\subset V$,
every $a\in F$ and every $r\in\NN$ put \0 $K_{a}^{0}F:=T_{a}F$ and
define \1 $K_{a}^{r+1}\!F$ to be the space of all vectors $v\in
K_{a}^{r}F$ such that there is a smooth mapping $f:V\to V$ with
$f'(a)(T_{a}F)\subset K_{a}^{r}F$, $f(a)=v$ and $f(x)\in K_{x}^{r}F$
for all $x\in F$.

\medskip\noindent It is clear that $K_{a}^{1}F=K_{a}F$ holds. Let us
call $F$ of {\sl uniform degeneracy} or {\sl uniformly degenerate} if
for every $r\in\NN$ the dimension of $K_{a}^{r}F$ does not depend on
$a\in F$. In that case it can be shown that for every $v\in
K_{a}^{r}F$ the outcome of the condition $f'(a)(T_{a}F)\subset
K_{a}^{r}F$ in (ii) does not depend on the choice of the smooth
mapping $f:V\to V$ satisfying $f(a)=v$ and $f(x)\in K_{x}^{r}F$ for
all $x\in F$. For instance, $F$ is of uniform degeneracy if $F$ is
locally affinely homogeneous, that is, if there exists a Lie algebra
$\7a$ of affine vector fields on $V$ such that every $\xi\in\7a$ is
tangent to $F$ and such that the canonical evaluation map $\7a\to
T_{a}F$ is surjective for every $a\in F$. Clearly, if $F$ is locally
affinely homogeneous in the above sense then the corresponding tube
manifold $M=F+iV$ is locally homogeneous as CR-manifold.

\medskip Recall our convention that every smooth map $f:V\to V$ is
considered as the smooth vector field $\xi=f(x)\dd x$ on $V$. Our
computations below are considerably simplified by the obvious fact
that every smooth vector field $\xi$ on V has a unique smooth
extension to $E$ that is invariant under all translations $z\mapsto
z+iv$, $v\in V$. In case $\xi$ is tangent to $F\subset V$ the
extension satisfies $\xi_{z}\in H_{z}M$ for all $z\in M$.

In case the submanifold $F\subset V$ is uniformly degenerate in a
neighbourhood of $a\in F$ we call $F\,$ {\sl affinely
$k$-nondegenerate} at $a$ if $K_{a}^{k}F=0$ and $k$ is minimal with
respect to this property. It can be seen that `affine
$k$-nondegeneracy' is invariant under affine coordinate changes. As a
consequence of \Lit{KAZT}, compare the last 5 lines in the Appendix
therein, we state:

\Proposition{UR} Suppose that $F$ is uniformly degenerate in a
neighbourhood of $a\in F$. Then the corresponding tube manifold
$M=F+iV$ is $k$-nondegenerate as CR-manifold at $a\in M$ if and
only if $F$ is affinely $k$-nondegenerate at $a$.\Formend

\Corollary{AF} Suppose $\dim F \ge2$ and $K_{x}F=\RR x$ for all $x\in
F$. Then $F$ is affinely $2$-nondegenerate at every point.

\Proof The map $f=\id$ has the property $f(x)\in K_xF$ for every $x\in
F$. Hence, the relation $f'(x)(T_{x}F)=T_{x}F\not\subset K_{x}F$
implies $x\notin K_{x}^{2}F$ and thus $K_{x}^{2}F=0$ as well as
$x\ne0$. In particular, $F$ is uniformly degenerate.\qed

\medskip For locally affinely homogeneous submanifolds $F\subset V$
the spaces $K_{a}^{r}F$ can easily be characterized. For each affine
vector field $\xi=h(x)\dd x$ on $V$ denote by
$\xi^{\lin}:=h-h(0)\in\End(V)$ the {\sl linear part} of $\xi$.

\Proposition{US} Suppose that $\7A$ is a linear space of affine vector
fields on $V$ such that every $\xi\in\7A$ is tangent to $F$ and the
canonical evaluation mapping $\7A\to T_{a}F$ is a linear
isomorphism. Then, given any $r\in\NN$, the vector $v\in K_{a}^{r}F$
is in $K_{a}^{r+1}\!F$ if and only if $\xi^{\lin}(v)\in K_{a}^{r}F$
for every $\xi\in\7A$.

\Proof By the implicit function theorem, there exist open
neighbourhoods $Y$ of $0\in\7A$ and $X$ of $a\in M$ such that
$g(y):=\exp(y)a$ defines a diffeomorphism $g:Y\to X$. Define the
smooth map $f:X\to V$ by $f(x)=\mu_{y}(v)$, where $\mu_{y}$ for
$y:=g^{-1}(x)$ is the linear part of the affine transformation
$\exp(y)$. Then $f(a)=v$ and $f(x)\in K_{x}^{r}F$ for every $x\in
X$. A simple computation shows that\medskip
\centerline{$f'(a)\big(g'(0)\xi\big)=\xi^{\lin}(v)\steil{for every}
\xi\in\7A$.}  \nline In view of Definition \ruf{HN}.ii this identity
implies the claim. \qed

\medskip It is easily seen that a necessary condition for $M$ being
minimal as CR-manifold is that $F$ is not contained in an affine
hyperplane of $V$. For later use in Proposition \ruf{PS} we state the
following sufficient condition.

\Proposition{PY} Suppose that $\7A$ has the same properties as in
Proposition \ruf{US} and denote by $\Lambda\subset\End(V)$ the {\rm
associative} real subalgebra generated by
$\{\xi^{\lin}:\xi\in\7A\}$. Then the tube manifold $M=F+iV$ is
minimal at $a$ if $V$ is the linear span of all vectors $\lambda(v)$
with $v\in T_{a}F$ and $\lambda\in\Lambda$.

\Proof Without loss of generality we assume that the canonical
evaluation mapping $\7A\to T_{x}F$ is a linear isomorphism for every
$x\in F$. We also assume that $V$ is the linear span of all
$\lambda(v)$ as above. Define inductively for every integer $r\ge1$
the subbundle $H^{r}M\subset TM$ in the following way: $H^{1}M:=HM$
and every $H^{r+1}_{z}M$, $z\in M$, is the linear span of $H^{r}_{z}M$
together with all vectors $[\xi,\eta]_{z}$, where $\xi,\eta$ are
arbitrary smooth sections in $H^{r}M$ over $M$. For the proof it is
enough to show
$T_{a}M=H^{\infty}_{a}M:=\bigcup_{r\ge1}\!T_{a}M$.\nline From
$T_{a}F\oplus iT_{a}F\,\subset\, H_{a}^{r}M\,\subset\, T_{a}F\oplus
iV$ we see that for every $1\le r\le\infty$ there is a unique linear
subspace $H^{r}_{a}F\subset V$ with $H^{r}_{a}M=T_{a}F\oplus
iH^{r}_{a}F$ and $H^{1}_{a}F=T_{a}F$. Therefore it is enough to show
$H^{\infty}_{a}F=V$. We claim that $\xi^{\lin}(H^{r}_{a}F)\subset
H^{r+1}_{a}F$ holds for all $\xi\in\7A$. To see this fix an arbitrary
$w\in H^{r}_{a}F$ and an arbitrary vector field $\xi\in\7A$. Choose a
smooth section $\eta$ over $F$ in the bundle $iH^{r}_{a}F$ with
$\eta_{a}=iw$ and extend $\eta$ as well as $\xi$ in the unique way to
smooth vector fields on $M$ that are invariant under all translations
$z\mapsto z+iv$ with $v\in V$. Then $\xi,\eta$ are sections in
$H^{r}M$, and $[\xi,\eta]_{a}=i\xi^{\lin}(w)\in H^{r+1}_{a}M$ implies
$\xi^{\lin}(w)\in H^{r+1}_{a}F$ as required. Now define inductively
the linear subspaces $W^{r}\subset V$ by $W^{1}:=H^{1}_{a}F=T_{a}F$
and $W^{r+1}$ as the linear span of $W^{r}$ together with all
$\xi^{\lin}(W^{r})$, $\xi\in\7A$. Then $V=\bigcup_{r\ge1}W^{r}$ by
assumption and $W^{r}\subset H^{r}_{a}F$ by induction gives $V\subset
H^{\infty}_{a}M\subset V$ as desired.\qed

\medskip

\Lemma{ZI} Suppose that $F\subset V$ is a submanifold such that for
every $c\in V$ with $c\ne0$ there exists a linear transformation
$g\in\GL(V)$ with $g(F)=F$ and $g(c)\ne c$ (this
condition is automatically satisfied if $F$ is a cone). Then for
$M=F+iV$ the CR-automorphism group $\Aut(M)$ has trivial center.

\Proof Let an element in the center of $\Aut(M)$ be given and let
$h:U\to E$ be its holomorphic extension to an appropriate connected
open neighbourhood $U$ of $M$. Since $h$ commutes with every
translation $z\mapsto z+iv$, $v\in V$, it is a translation itself:
Indeed, for $a\in F$ fixed and $c:=h(a)-a$ the translation
$\tau(z):=z+c$ coincides with $h$ on $a+iV$ and hence on $U$ by the
identity principle. For every $g\in\GL(V)\cap\Aut(M)$ the identity
$gh=hg$ implies $g(c)=c$. This forces $c=0$ by our
assumption, that is, $h(z)\equiv z$.\qed

\Proposition{UQ} Suppose that the homogeneous CR-manifold $M$ is
simply connected and that $\Aut(M)$ has trivial center. In case the
stability group $\Aut(M,a)$ is trivial for some (and hence every)
$a\in M$, the following properties hold: \0 Let $M'$ be an arbitrary
homogeneous CR-manifold and $D\subset M$, $D'\subset M'$ non-empty
domains. Then every real-analytic CR-isomorphism $h:D\to D'$ extends
to a real-analytic CR-\-isomorphism $M\to M'$. \1 Let $M'$ be an
arbitrary locally homogeneous CR-manifold and $D'\subset M'$ a domain
that is CR-equivalent to $M$. Then $D'=M'$.

\Proof ad (i): Fix a point $a\in D$. To every $g\in G:=\Aut(M)$
with $g(a)\in D$ there exists a transformation $g'\in G':=\Aut(M')$
with $hg(a)=g'h(a)$. Because of $\Aut(M',a')=\{\id\}$ the
transformation $g'$ is uniquely determined by $g$ and satisfies
$hg=g'h$ in a neighbourhood of $a$. Since the Lie group $G$ is simply
connected $g\mapsto g'$ extends to a group homomorphism $G\to G'$ and
$h$ extends to a CR-covering map $h:M\to M'$. The deck transformation
group $\Gamma:=\{g\in G:gh=h\}$ is in the center of $G$ and hence is
trivial by assumption. Therefore, $h:M\to M'$ is a
CR-isomorphism. \nline ad (ii): The proof is essentially the same
as for Proposition 6.3 in \Lit{FEKA}.\qed

The condition `locally homogeneous' in Proposition \ruf{UQ}.ii cannot
be omitted. A counterexample is given for every integer $k\ge3$
by the tube $M'\subset\CC^{3}$ over the cone
$$
F':=\{x\in\RR^{3}:x_{2}^{k}=x_{1}^{k-1}x_{3},\;x_{1}^{2}+
x_{2}^{2}>0\}\,.
$$ Then with $\;\RRp:=e^{\RR}$ the tube $M$ over
$F:=F'\cap(\RRp)^{3}$ is the Example \ruf{EX} below for $\theta=k$,
and $M,M'$ satisfy for $D'=M$ the assumption of Proposition
\ruf{UQ}.ii.

\ifarx\vfill\eject\fi

\KAP{Vier}{Tube manifolds over cones}

In this section we always assume that the submanifold $F\subset V$ is
conical (that is, $x\in T_xF$ for every $x\in F$) and that
$a\in F$ is a given point. Then, for $M:=F+iV$, the Lie algebra
$\7g:=\hol(M,a)$ contains the Euler vector field $\delta:=z\dd z$.
Denote by $\7P$ the complex Lie algebra of all polynomial holomorphic
vector fields $f(z)\dd z$ on $E$, that is, $f:E\to E$ is a polynomial
map. Then $\7P$ has the $\ZZ$-grading
$$
\7P=\bigoplus_{k\in\ZZ}\7P_{k}\,,\qquad[\7P_{k},\7P_{l}]\subset
\7P_{k+l}\,,\Leqno{HJ}
$$
where $\7P_{k}$ is the $k$-eigenspace of $\ad(\delta)$ in $\7P$. Then
$\7P_{k}$ is the subspace of all $({k+1})$-homogeneous vector fields
in $\7P$ if $k\ge-1$ and is $0$ otherwise. Define
$\7g_k:=\7g\cap\7P_k.$ Clearly, $\bigoplus_{k\in \ZZ}\7g_k$ is a
graded, in general proper subalgebra of $\7g$, and
$\aff(M,a)=\7g_{-1}\oplus\7g_{0}$.

\Proposition{DC} Retaining the above notation, suppose that
$\7g=\hol(M,a)$ has finite dimension. Then \0 $\7g\subset\7P$, that
is, every $f(z)\dd z\in\7g$ is a polynomial vector field on $E=V\oplus
iV$. Furthermore, $f(iV)\subset iV$.\vskip3pt \1 $\displaystyle
\7g=\bigoplus\limits_{k\ge-1}\7g_{k}\,,\quad[\7g_{k},\7g_{l}]\subset
\7g_{k+l}\Steil{and}\7g_{-1}=\{iv\dd z: v\in V\}\;$.\vskip2pt \1 For
every $z\in M$ the canonical map $\hol(M)\to\hol(M,z)$ is a Lie
algebra isomorphism. \1 $\7g_{k}=0$ for some $k\in\NN$ implies
$\7g_{j}=0$ for every $j\ge k$.

\Proof Consider $\7l:=\7g\oplus i\7g\subset\hol(E,a)$, which
contains the vector field $\eta:=(z-a)\dd z$. We first show
$\7l\subset \7P$: Fix an arbitrary $\xi:=f(z) \dd z\in\7l$. Then in a
certain neighbourhood of $a\in E$ there exists a unique expansion
$\xi=\sum_{k\in\NN}\xi_{k}$, where $\xi_{k}=p_{k}(z-a)\dd z$ for a
$k$-homogeneous polynomial map $p_{k}:E\to E$. It is easily verified
that the vector field $\ad(\eta)\xi\in\7l$ has the expansion
$\ad(\eta)\xi=\sum_{k\in\NN}(k{-}1)\xi_{k}$. Now assume that for
$d:=\dim \7l $ there exist indices $k_{0}<k_{1}<\dots<k_{d}$ such that
$\xi_{k_{l}}\ne0$ for $0\le l\le d$. Since the Vandermonde matrix
$\big((k_{l}-1)^{j}\big)$ in non-singular, we get that the vector
fields $(\ad\eta)^{j}\xi=\sum_{k\in\NN}(k{-}1)^{j}\xi_{k}$, $0\le j\le
d$, are linearly independent in $\7l$, a contradiction. This implies
$\xi\in\7P$ as claimed.\nline Since $\7g\subset\7P$ has finite
dimension, every $\eta\in\7g$ is a finite sum
$\eta=\sum_{k=-1}^{m}\eta_{k}$ with $\eta_{k}\in\7P_k$ and $m\in\NN$
not depending on $\eta$. For every polynomial $p\in\RR[X]$ then
$p(\ad\delta)\eta=\sum_{k=-1}^{m}p(k)\eta_{k}$ shows
$\eta_{k}\in\7g_{k}$ for all $k$, that is, $\7g=\bigoplus\7g_{k}$. The
identity $\7g_{-1}=\{iv\dd z: v\in V\}$ follows from the fact that
$\7g_{-1}$ is totally real in $\7P_{-1}$and this implies $f(iV)\subset
iV$ for all $f(z)\dd z\in\7g_{k}$ by
$[\7g_{-1},\7g_{k}]\subset\7g_{k-1}$ and induction on $k$. For the
proof of the last claim assume $\7g_{k}=0$ for some $k\ge0$. Then
$[\7P_{-1},\7g_{k+1}]\subset (\7g_{k}\oplus i\7g_{k})=0$ implies
$\7g_{k+1}\subset\7g_{-1}$ and hence $\7g_{k+1}=0$. \qed

\Corollary{ZU} In case $\7g=\hol(M,a)$ has finite dimension the
CR-manifold $M=F+iV$ is locally homogeneous at $a$ if and only
if $F$ is locally linearly homogeneous at $a\in F$.

\Proof From Proposition \ruf{DC} follows that
$W:=\{\xi_{a}:\xi\in\bigoplus_{k\in\NN}\7g_{2k}\}$ is a subspace of
$V$ while $\{\xi_{a}:\xi\in\bigoplus_{k\in\NN}\7g_{2k-1}\}=iV.$ Hence,
$M$ is locally homogeneous at $a$ if and only if $W=T_{a}F$. But for
every $k\in\NN$ and every $\xi\in\7g_{2k}$ the vector field
$\eta:=(\ad ia\dd z)^{2k}\xi$ is in $\7g_{0}$ and satisfies
$\eta_{a}=(-1)^{k}(2k+1)!\,\xi_{a}$, that is,
$W=\{\xi_{a}:\xi\in\7g_{0}\}$.\qed

Notice that the conclusion $\7g\subset\7P$ together with an eigenspace
decomposition as in Proposition \ruf{DC} can for $V=\RR^{n}$ be
obtained in the same way if instead of `$F$ conical' it is only
assumed for $M=F+i\RR^{n}$ that $\7g=\hol(M,a)$ contains a
vector field $\alpha_{1}z_{1}\dd{z_{1}}+
\alpha_{2}z_{2}\dd{z_{2}}+\dots+\alpha_{n}z_{n}\dd{z_{n}}$ with
$\alpha_{k}>0$ for all $k$, compare e.g. \Ruf{RU}.

\smallskip

\Proposition{ZB} Assume that $\7g:=\hol(M,a)$ has finite dimension and
that $F'\subset V$ is a further conical submanifold with tube manifold
$M'=F'+iV$ and $\7g':=\hol(M',a')$ for some $a'\in F'$. Assume
that the CR-germs $(M,a)$ and $(M',a')$ are isomorphic and let
$g,\tilde g:(M,a)\to(M',a')$ be arbitrary CR-isomorphisms. Then \0
$\dim\7g_{k}=\dim\7g_{k}'$ for all $k\in\ZZ$, where
$\7g^{}_{k},\7g'_{k}$ are given by the decomposition \ruf{DC}.ii.  \1
$g$ is represented by a rational transformation on $E$.  \1 In case
$\7g_{1}=0$, $\,g$ is represented by a linear transformation in
$\GL(V)\subset\GL(E)$ mapping every $K_{a}^{r}F$ onto
$K_{a'}^{r}F'$. \1 $g=\tilde g$ if and only if $g,\tilde g$ have the
same $d$-jet at $a$, where $d:=\min\{k\in\NN:\7g_k=0\}$.

\Proof Let $\7l:=\7g\oplus i\7g$ and $\7l_{k}:=\7l\cap\7P_{k}$ for all
$k$. The Lie algebra automorphism $\Psi:=\exp(\ad a\dd z)$ of $\7l$
maps every $f(z)\dd z$ to $f(z+a)\dd z$. For every $k$ denote by
$\7l^{k}\subset\7l$ the subspace of all vector fields that vanish of
order at least $k{+}1$ at $a$. Then $\Psi(\7l^{k})=\bigoplus_{j\ge
k}\7l_{j}$ implies $\dim\7l_{k}=\dim\7l^{k}/\7l^{k+1}$. As a
consequence, $\dim_{\RR}\7g_{k}=\dim_{\CC}\!\7l_{k}$ is a CR-invariant
of the germ $(M,a)$ for every $k$.\nline For the proof of (ii), (iii)
put $\7l':=\7g'+i\7g'$ and extend $g$ to a biholomorphic mapping
$g:U\to U'$ with $g(a)=a'$ and $g(U\cap M)=U'\cap M'$ for suitable
connected open neighbourhoods $U,U'$ of $a,a'\in E$. Consider
the induced Lie algebra isomorphism $ g_{*}:\7l\to\7l'$, compare
\Ruf{GP}. Its inverse $\Theta:=g_{*}^{-1}$ is given by
$$
\Theta\big(f(z)\dd z\re1\big)=g'(z)^{-1}\!f\big(g(z)\big)\re2\dd
z\,.\Leqno{ZA}
$$
Since $\7l$ consists of polynomial vector fields (Proposition
\ruf{DC}) there exist {\sl polynomial} maps $p:E\to E$ and $
q:E\to\End(E)$ such that
$$
\Theta\big(z\dd z\big)=p(z)\dd z\Steil{and}\Theta\big(e\dd
z\big)=\big(q(z)e\big)\dd z
$$
for all $e\in E$. Then \Ruf{ZA} implies $g'(z)^{-1}=q(z)$ and
$g'(z)^{-1}g(z)=p(z)$, that is,
$$
g(z)=q(z)^{-1}p(z)\Leqno{ZV}
$$
in a neighbourhood of $a\in E$ and, in particular, $g$ is
rational.\nline Now suppose $\7g_{1}=0$. Then also
$\7g_{k}=\7g'_{k}=0$ for all $k\ge1$ by \ruf{DC}.iv and (i).  Clearly
$\Theta(\7l'_{a'})=\7l_{a}$, where $\7l_{a}:=\{\xi\in\7l:\xi_{a}=0\}$
and similarly $\7l'_{a'}\subset\7l'$ are the isotropy subalgebras at
$a,a'$. Also, $\delta_{a}:=(z-a)\dd z$ is the unique element in
$\7l_{a}$ such that $\ad(\delta_{a})$ induces the negative identity on
the factor space $\7l/\7l_{a}$. Since $\delta_{a'}$ has the same
uniqueness property for $\7l'_{a'}$ we must have
$\Theta(\delta_{a'})=\delta_{a}$.  Since $\delta=z\dd z$ is the
$\7g$-component of $\delta_{a}$ as well as of $\delta_{a'}$ in
$\7g\oplus i\7g$ we actually get $\Theta(\delta)=\delta$, that is,
$p(z)\equiv z$.  Also $\7l_{-1}$ is $\Theta$-invariant, implying that
$q$ is constant. Therefore $g=q(a)^{-1}$.  \nline For the proof of
(iv) we assume without loss of generality that $M=M'$, $a=a'$ and that
$g,\id\in\Aut(M,a)$ have the same $d$-jet at $a$. This implies that
$(g(z)-z)$ vanishes of order $>d$ and $(g'(z)-\id)$ vanishes of order
$\ge d$ at $a$. Therefore also $(g'(z)^{-1}-\id)=(q(z)-\id)$ vanishes
of order $\ge d$ at $a$. Since $q$ is a polynomial of degree $\le d$,
there is a $d$-homogeneous polynomial $s:E\to\End(V)$ with
$q(z)=\id+s(z-a)$. Consider the vector field $\eta:=(z-a)\dd z\in\7l$
and define the holomorphic mappings $h,r:U\to E\,$ by
$$
h(z):=q(z)\big(g(z)-a\big)=(z-a)+r(z)\,.
$$
Then $\Theta\big(\eta)=h(z)\dd z\;\in\;\7g$ shows that $h$ and $r$ are
polynomials of degree $\le d$. But
$$
r(z)=\big(g(z)-z\big)+s(z-a)\big(g(z)-a\big)
$$
vanishes of order $>d$ at $a$, that is, $r=0$ and
$\Theta(\eta)=\eta$. This implies $\Theta(\7l_{-1})=\7l_{-1}$ since
$\7l_{-1}$ is the $(-1)$-eigenspace of $\ad(\eta)$. Therefore $g$ is
an affine transformation on $E$. From $g(a)=a$ and $g'(a)=\id$ we
finally get $g=\id$.\qed

\Corollary{PL} Let $M:=F+iV$, $M':=F'+iV$ with conical
submanifolds $F,F'\subset V$ and let $a\in F$, $a'\in F'$ be arbitrary
points. Assume furthermore that $\hol(M,a)=\aff(M,a)$ holds. Then the
following conditions are equivalent: \0 The manifold germs
$(M,a),\,(M',a')$ are CR-equivalent. \1 The manifold germs
$(M,a),\,(M',a')$ are affinely equivalent. \1 The manifold germs
$(F,a),\,(F',a')$ are linearly equivalent.\Formend

\bigskip\noindent Recall that $\aut(M,a)\subset\7g=\hol(M,a)$ is
defined as the isotropy subalgebra at $a$ and $\Aut(M,a)$ is the
CR-automorphism group of the manifold germ $(M,a)$, also called the
{\sl stability group} at $a\in M$.

\Proposition{LO} Let $F\subset V$ be conical and $M=F+iV$. In
case $\7g=\hol(M,a)$ has finite dimension, the following conditions
are equivalent: \0 $\7g_{1}=0$. \1 $\7g=\aff(M,a)$. \1 The tangential
representation $\;g\mapsto g'(a)$ induces a group monomorphism
$\Aut(M,a)\;\hookrightarrow\;\GL(V)$.  \smallskip\noindent Each of
these conditions is satisfied if $\aut(M,a)=0$.

\Proof Let $\7l:=\7g\oplus i\7g\subset\hol(E,a)$ and
$\7l_{k}:=\7g_{k}\oplus i\7g_{k}$ for all $k$. \nline\To12 Follows
from the last claim in Proposition \ruf{DC}. \nline\To23 By
Proposition \ruf{ZB}.iii every $g\in\Aut(M,a)$ is represented by a
linear transformation on $V$. \nline\To31 Let $\xi\in\7g_{1}$ be an
arbitrary vector field. Then there exists a unique symmetric bilinear
map $b:E\times E\to E$ with $\xi=b(z,z)\dd z$. Now $(\ad ia\dd
z)^{2}\xi=-2b(a,a)\dd z\in\7g$, that is, $\eta:=h(z)\dd z$ is in
$\aut(M,a)$, where $h(z):=b(z,z)-b(a,a)$. For every $t\in\RR$
therefore the transformation $\psi_{t}:=\exp(t\eta)\in\Aut(M,a)$ has
derivative $\psi_{t}'(a)=\exp(th'(a))\in\GL(E)$ in $a$. But
$\psi_{t}'(a)\in\GL(V)$ by \ruf{ZB}.iii and thus $2b(a,v)=h'(a)v\in V$
for all $v\in V$. On the other hand $b(a,v)\in iV$ by Lemma \ruf{ZI},
implying $\psi_{t}'(a)=\id$ for all $t\in\RR$. By the injectivity of
the tangential representation therefore $\eta_{t}$ does not depend on
$t$ and we get $\xi=0$. This proves (i) and thus the equivalence of
(i) -- (iii).\nline Suppose $\aut(M,a)=0$ and that there exists a
non-zero vector field $\xi\in\7g_{1}$. Then $\xi_{a}\in iV$ and there
exists an $\eta\in\7g_{-1}$ with $\xi-\eta\in\aut(M,a)$, a
contradiction.\qed

\Remark{} Notice that the condition (iii) in Proposition \ruf{LO}
states that the tangential representation takes its values in the
subgroup $\GL(V)\subset\GL(E)$. In general, the tangential
representation is not injective and also takes values outside
$\GL(V)$. The tube $\5M$ over the future light cone can serve as a
counterexample for both of these phenomena.\Formend

\Proposition{PS} Suppose that $M=F+iV$, $F\subset V$ a conical
submanifold, is locally homogeneous and that
$\7g=\7g_{-1}\oplus\7g_{0}$ for $\7g=\hol(M,a)$. Then the tangential
representation at $a$ induces a group isomorphism
$$
\Aut(M,a)\;\cong\;\{g\in\GL(V):g\7g_{0}g^{-1}=\7g_{0} \steil{and}g(a)=
a\}\,,
$$
where $\7g_{0}$ is considered in the canonical way as linear subspace
of $\,\End(V)$.

\Proof The assumptions imply that $H(a)\cap F$ is a neighbourhood of
$a$ in $F$ for $H:=\exp(\7g_{0})\subset\GL(V)$. Let $g\in\GL(V)$ be an
arbitrary linear transformation with $g(a)=a$ and
$g\7g_{0}g^{-1}=\7g_{0}$. Then $gHg^{-1}=H$ and hence
$gH(a)=Hg(a)=H(a)$, that is $g\in\Aut(M,a)$. Conversely, by \ruf{ZB}.iii
every element of $\Aut(M,a)$ can be represented by some
$g\in\GL(V)$ with $g(a)=a$. The Lie algebra automorphism
$\Theta=g_{*}$ of $\7g$ leaves $\delta$ and thus also $\7g_{0}$
invariant. As a consequence, $\Theta(\phi)=g\phi g^{-1}$ for every
$\phi\in\7g_{0}$, that is, $g\7g_{0}g^{-1}=\7g_{0}$.\qed

\medskip The real Lie algebra structure of $\hol(M,a)$ is a
CR-invariant for the manifold germ $(M,a)$. For certain classes of
conical tube manifolds this gives a complete invariant:

\Proposition{PQ} Let $F,F'\subset V$ be conical submanifolds for which
the corresponding tubes $M:=F+iV$, $M':=F'+iV$ are locally
homogeneous CR-manifolds. Let $a\in F$, $a'\in F'$ be arbitrary points
and assume that for $\7g:=\hol(M,a)$ the spaces $\7g_{k}$ occurring in
the gradation \ruf{DC}.ii satisfy $\7g_{k}=[\7g_{0},\7g_{0}]=0$. Then
the following conditions are equivalent: \0 The germs
$(M,a),\,(M',a')$ are CR-equivalent. \1 $\hol(M,a)$, $\,\hol(M',a')$
are isomorphic as real Lie algebras.

\Proof Only the implication \To21 is not obvious. Suppose that
$\Theta:\7g\to\7g':=\hol(M',a')$ is a Lie algebra isomorphism. We use
the same symbol for the complex linear extension $\Theta:\7l\to
\7l'$. \nline Our first step is to show that $\Theta$ can be modified
in such a way that it respects the gradations. To begin with,
$[\7g',\7g']$ is abelian since our assumption implies that
$[\7g,\7g]=\7g_{-1}$ has this property. With $\7g'=\bigoplus \7g_k'$
being the gradation for $\7g'$ as in \ruf{DC}.ii assume that there
exists a minimal integer $k\ge1$ with $\7g'_{k}\ne0$. Then
$\7g'_{-1}=[\delta,\7g'_{-1}]$ and $\7c:=[\7g'_{-1},\7g'_{k}]\ne0$
imply $\7g'_{-1},\7c\subset[\7g',\7g']$ together with
$[\7g'_{-1},\7c]\ne0$, a contradiction. Therefore
$\7g'=\7g'_{-1}\oplus\7g'_{0}$ with $[\7g',\7g']=\7g'_{-1}$, and as a
consequence,
$\Theta(\7g_{-1})=\Theta([\7g,\7g])=[\7g',\7g']=\7g_{-1}'$. Since
$\ad:\7g'_{0}\to \gl(\7g'_{-1})$ is injective, $\delta+\7g'_{-1}$ is
precisely the set of all $\xi\in\7g'$ such that $\ad\xi$ induces on
$\7g'_{-1}$ the negative identity. Therefore, there exists
$\eta\in\7g'_{-1}$ with $\Theta(\delta)=\delta-\eta$. Replacing
$\Theta$ by $\exp(\ad \eta)\Theta=(\id+\ad \eta)\Theta,$ we get
$\Theta(\delta)=\delta$ and finally $\Theta(\7g_{k})=\7g'_{k}$ for all
$k$.\nline There exists a linear operator
$\theta\in\GL(V)\subset\GL(E)$ with $\Theta(e\dd z)=\theta(e)\dd z$
for all $e\in E=V\oplus iV$. We claim that $H'=\theta H\theta^{-1}$
holds for the abelian subgroups $H:=\exp \7g_0$ and $H':=\exp\7g_0'$
of $\GL(V)$. Indeed, application of $\Theta$ to $\big[e\dd
z,\lambda(z)\dd z\big]=\lambda(e)\dd z\;$ yields
$$
\tilde\lambda \theta=\theta\lambda\Steil{ for all}\lambda(z)\dd
z\in\7g_{0}\steil{and}\tilde\lambda(z)\dd z:=\Theta\big(\lambda(z)\dd
z\big)\in\7g'_{0}\,.
$$
We may therefore assume (possibly after replacing $F$ by $\theta F$
and $a$ by $\theta a$) that $H=H'$, $F=H(a)$ and $F'=H(a')$. Denote by
$\Lambda\subset\End(V)$ the associative subalgebra generated by all
$\lambda\in\End(V)$ with $\lambda(z)\dd z\in\7g_{0}$. Then $\Lambda$
is abelian, contains the identity of $\End(V)$ and
$H\subset\Lambda$. Since $\dim\7g<\infty$, the CR-manifold $F+
iV$ is minimal and consequently $F$ cannot be contained in a
hyperplane of $V$. This implies $\Lambda(a)=V$ and thus the existence
of a $g\in\Lambda$ with $a'=g(a)$. From $gH=Hg$ we get
$F'=g(F)$. Since $g(F)$ also cannot be contained in a hyperplane of
$V$ finally $g\in\GL(V)$ follows.\qed

In Propositions \ruf{UH} an d \ruf{HU} a large class of linearly
homogeneous conical submanifolds $F\subset V$ is given for which the
corresponding tubes $M$ satisfy the condition
$\7g_{1}=[\7g_{0},\7g_{0}]=0$ in Proposition \ruf{PQ}.

We note that we do not know a single example with $\dim\7g<\infty$ and
$\dim\7g_{k}>\dim\7g_{-k}$ for some $k\in\NN$.  We also do not know
any pair $M,M'$ of holomorphically nondegenerate conical tube
manifolds, that are locally CR-equivalent but are not locally affinely
equivalent. For Levi nondegenerate tube manifolds (which necessarily
cannot be conical) such examples can be found in \Lit{DAYA}. In
\Lit{ISMI} even two affinely homogeneous examples are contained which
are locally affinely non-equivalent but whose associated tube
manifolds are locally CR-equivalent.

\KAP{Fuenf}{Some examples}

In this section we present two classes of examples.  To our knowledge,
the only known example of a homogeneous $k$-nondegenerate CR-manifold
with $k\ge3$ occurs in \Lit{FELS} for the case $k=3$: That is a
hypersurface $M$ in a $7$-dimensional compact complex manifold, on
which the simple Lie group $\SO(3,4)^{0}$ acts by biholomorphic
transformations with orbit $M$. In our first example we give for
arbitrary CR-codimenion $c\ge1$ a minimal homogeneous
$3$-nondegenerate as well as a minimal homogeneous $4$-nondegenerate
CR-manifold.  The second class of examples deals with tubes $M$ over
cones of the form $\{x\in(\RRp\!)^{n}:\sum_{j\le p}
x_{j}^{\alpha}=\sum_{j>p} x_{j}^{\alpha}\}$ with $1\le p<n$ and
$\alpha\ne0,1$. Using results from the preceding section, we
explicitly determine all Lie algebras $\7g=\hol(M,a)$ for certain
$\alpha$. Among these are all hyperquadrics (that is $\alpha=2$) of
signature $(p,q)$ with $q:=n-p$, where $\7g$ turns out to be
isomorphic to $\so(p{+}1,q{+}1)$.

The first example of CR-manifolds introduced below consists of tubes
$M=\Gamma(a)+iV$ over certain group orbits $\Gamma(a)$, where the
connected group $\Gamma:=\{g\in\GL(2,\RR):\det(g)>0\}$ acts linearly
on a real vector space $V$. In that way we obtain homogeneous
$k$-nondegenerate CR-manifolds for $k\in\{2,3,4\}$. For dimension
reasons it is impossible to construct CR-manifolds of higher
nondegeneracy, employing the group $\Gamma=\GL(2,\RR)^{0}$.  We do not
know how to construct $k$-nondegenerate homogeneous tube manifolds with
$k\ge5$ (should these exist) with suitable other groups, either.

\Example{JP} For fixed integers $k\in\{2,3,4\}$ and $c\ge1$ let
$V\subset\RR[u_{1},u_{2}]$ be the subspace of all homogeneous
polynomials of degree $m:=k+c-1$. Then the group $\Gamma$ (see above)
acts irreducibly on $V$ by $p\mapsto p\circ g^{-1}$ for all
$g\in\Gamma$ and has the subgroup $\{g\in\RR\id:g^{m}=\id\}$ as kernel
of ineffectivity. For $a:=\sum^{k-2}_{j=0}u_{1}^{j}u_{2}^{m-j}\in V$
the orbit $F=\5F^{k,c}:=\Gamma(a)$ is a connected conical submanifold
of dimension $k$ in $V$.  With Proposition \ruf{US} it is easily seen
that
$$
K^{r}_{a}F=\sum^{k-1-r}_{j=0}\RR\,u_{1}^{j}u_{2}^{m-j}\steil{for
all}r\ge0\,.
$$
In particular, the tube manifold $\5M^{k,c}:=F+iV$ is a
$k$-nondegenerate homogeneous CR-manifold of CR-dimension $k$ and
CR-codimension $c$. The manifolds $\5M^{2,c}$ will be discussed in
more detail in Section \ruf{Sieben}. In particular, $\5M^{2,1}$
is linearly equivalent to the future light cone tube $\5M$.  For
easier handling let us identify $\RR^{m+1}$ with $V$ via
$(x_{0},x_{1},\dots,x_{m})
\;\leftrightsquigarrow\;\sum^{m}_{j=0}x_{j}u_{1}^{j}u_{2}^{m-j}\;.$
Since $\Gamma$ acts on $\5M^{k,c}$ by (linear) CR-automorphisms, the
linear part $\7g_0$ of $\7g=\hol(\5M^{k,c},a)$ (compare \ruf{DC}.ii)
contains a copy of $\gl(2,\RR)$. More explicitly, this subalgebra is
spanned by the vector fields
$$
\eqalign{ \zeta_{1}:&=\sum^{m}_{j=0}jz_{j}\dd{z_{j}}\,,\mskip53.6mu
    \xi^{-1,1}:=\sum^{m-1}_{j=0}(m-j)z_{j+1}\dd{z_{j}}
\cr
\xi^{1,-1}:&=\sum^{m}_{j=1}jz_{j-1}\dd{z_{j}} \,,\mskip59mu
   \zeta_{2}:=\sum^{m}_{j=0}(m-j)z_{j}\dd{z_{j}}
\;.\cr}\Leqno{AQ}
$$
In particular, the vector fields $\xi^{1,-1},\,\xi^{-1,1}$ and
$\xi^{0,0}:=\zeta_{1}-\zeta_{2}=\sum^{m}_{j=0}(2j-m)z_{j}\dd{z_{j}}$
span a copy of $\7{sl}(2,\RR)$ in $\7g_{0}$. In case $k=2$ a
straightforward computation shows that
$\Gamma_{a}=\Big\{\Matrix\alpha\beta0\epsilon\in\Gamma:
\epsilon^{m}=1\Big\}$ is the isotropy subgroup at $a\in\5M^{2,c}$.
Hence,
$\Big\{\Matrix\alpha\beta{-\beta}\alpha\in\Gamma\Big\}\cong\CC^{*}$
acts transitively\vadjust{\vskip4pt} on $F$, that is, $\5M^{2,c}$ is
diffeomorphic to $\RR^{n}\times\CC^{*}$, where $n:=m+1=\dim V$.

\medskip\noindent It seems quite hard to find explicit global
equations for $\5M^{k,c}$ in case $k>2$. For $k=2$ see \Ruf{IP}. Only
for $k=3$ and $c=1$ we have the following: Consider in $\RR^{4}$ the
algebraic hypersurface given by the homogeneous equation
$$
S:=\{x\in\RR^{3}:x^{2}_{0}x_{3}^2+4x_{0}x_{2}^{3}-6x_{0}x_{1}x_{2}x_{3}-
3x_{1}^{2}x_{2}^{2}+4x_{1}^{3}x_{3}=0\}\,.\Leqno{PI}
$$ It is obvious that $S$ contains the point $a=(1,1,0,0)$. Since the
application of every vector field in $\Ruf{AQ}$ to the defining
function of $S$ gives a multiple of this function, $S$ is invariant
under the group $\Gamma$. Therefore the $3$-nondegenerate homogeneous
CR-submanifold $\5M^{3,1}=\Gamma(a)+i\RR^{4}$ of $\CC^{4}$ is an
open piece of the tube $S+i\RR^{4}$.

\bigskip\Example{EB} Fix integers $p\ge q\ge1$ with $n:=p+q\ge3$ and a
real number $\alpha$ with $\alpha^2\ne\alpha$. Then
$$
F=F_{p,q}^{\alpha}:=\Big\{x\in(\RRp)^{n}:\sum_{j=1}^{p}
x_{j}^{\alpha}=\sum_{j=p+1}^{n}
x_{j}^{\alpha} \Big\}\,,
$$
is a hypersurface in $V:=\RR^{n}$. Furthermore, $F$ is a cone and
therefore $\dim K_{a}F\ge1$ for every $a\in F$. On the other hand, the
second derivative at $a$ of the defining equation for $F$ gives a
non-degenerate symmetric bilinear form on $V\times V$, whose
restriction to $T_{a}F\times T_{a}F$ then has a kernel of dimension
$\le1$. Therefore $\dim K_{a}F=1$ for every $a\in F$ and by Corollary
\ruf{AF} the CR-manifold $M=M_{p,q}^{\alpha}:=F_{p,q}^{\alpha}+
iV$ is everywhere $2$-nondegenerate (compare Example 4.2.5 in
\Lit{EBFT} for the special case $n=\alpha=3$). Since $M$ as
hypersurface is also minimal, $\7g=\hol(M,a)$ has finite dimension.

For the special case $\alpha=2$ and $q=1$ the above cone
$F=F_{n-1,1}^{2}$ is an open piece of the future light cone
$$
\{x\in\RR^{n}:x_{1}^{2}+\dots+x_{n-1}^{2}=x_{n}^{2},\;x_{n}>0\}
$$
in $n$-dimensional space-time, which is affinely homogeneous. In
\Lit{KAZT} it has been shown that for the corresponding tube manifold
$M$ the Lie algebra $\7g=\hol(M,a)$ is isomorphic to $\so(n,2)$ for
every $a\in M$. In case $q>1$ the following result seems to be new:

\noindent{\fett Case $\alpha=2$:} Consider on $\CC^{n}$ the
symmetric bilinear form $\langle
z|w\rangle:=\sum\epsilon_{j}z_{j}w_{j}$. Then $F$ is an open piece of
the hypersurface $\{x\in\RR^{n}:x\ne0\,,\;\langle
x|x\rangle=0\}$, on which the reductive group $\RR^{*}\cd\O(p,q)$ acts
transitively. Therefore, $\7s_{0}^{}:=\RR\delta\oplus\so(p,q)$ is
contained in $\7g_{0}$. One checks that
$$
\7s:=\7g_{-1}\oplus\7s_{0}^{}\oplus\7s_{1}^{}\,,\qquad\7s_{1}^{}:=
\big\{\big(2i\langle c|z\rangle z-i\langle z|z\rangle c\big)\dd
z:c\in\RR^{n}\big\}
$$
is a Lie subalgebra of $\7g$. The radical $\7r$ of $\7s$ is
$\ad(\delta)$-invariant and hence of the form
$\7r=\7r_{-1}\oplus\7r_{0}\oplus\7r_{1}$ for
$\7r_{k}^{}:=\7r\cap\7g_{k}$. From $\so(p,q)$ semisimple we
conclude $\7r_{0}\subset\RR\delta$. But $\delta$ cannot be in $\7r$
since otherwise $\7g_{-1}\subset\7r$ would give the false statement
$[\7g_{-1},\7s_{1}^{}]\subset\RR\delta$. Therefore $\7r_{0}^{}=0$, and
$[\7g_{-1},\7r_{1}^{}]=[\7r_{-1}^{},\7s_{1}^{}]=0$ implies
$\7r=0$. Now Proposition 3.8 in \Lit{KAZT} implies $\7g=\7s$, and, in
particular, that $\7g$ has dimension ${n+2\choose2}$. In fact, it can
be seen that $\7g$ is isomorphic to $\so(p+1,q+1)$.

\medskip\noindent{\fett Case $\alpha$ an integer $\ge3$:} Then $F$
is an open piece of the real-analytic submanifold
$$
S:=\Big\{x\in\RR^{n}:x\ne0\Steil{and}\sum_{j=1}^{n}\epsilon_{j}
x_{j}^{\alpha}= 0\Big\}\Leqno{NP}
$$
which is connected in case $q>1$ and has two connected components
otherwise. For every $x\in\RR^{n}$ let $d(x)\in\NN$ be the cardinality
of the set $\{j:x_{j}=0\}$. It is easily seen that $\dim
K_{x}S=1+d(x)$ holds for every $x\in S$. Now consider the group
$$
\GL(F):=\{g\in\GL(V):g(F)=F\}\,.
$$
Every $g\in\GL(F)$ leaves $S$ and hence also
$H:=\{x\in\RR^{n}:d(x)>0\}$ invariant, that is, $g$ is the product of
a diagonal with a permutation matrix. Inspecting the action of
$\GL(F)$ on $\{c\in\overline F:d(c)=n-2\}$ we see that $\GL(F)$ as
group is generated by $\RRp\!\cd\id$ and certain coordinate
permutations. As a consequence, $\7g_{0}=\RR\delta$. Now suppose that
there exists a non-zero vector field $\xi\in\7g_{1}$. Then
$\xi=q(z,z)\dd z$ for some symmetric bilinear map
$q:\CC^{n}\times\CC^{n}\to\CC^{n}$ with $q(c,z)\dd z\in\CC\delta$ for
every $c\in\CC^{n}$. Because of $q\ne0$ symmetric therefore any two
vectors in $\CC^{n}$ must be linearly dependent, which contradicts
$n\ge3$.  Therefore $\7g_{1}=0$ and hence $\7g_{k}=0$ for all $k\ge1$
by Proposition \ruf{ZB}.iv. In particular, $\dim\7g=n+1<\dim M$ for
the tube manifold $M=F+iV$, that is, $M$ is not locally
homogeneous. For $n=3$ this gives an alternative proof for Proposition
6.36 in \Lit{EBEN}.

\Proposition{} Let $a,a'\in F=F^{\alpha}_{p,q}$ be arbitrary
points. Then in case $3\le\alpha\in\NN$ the CR-manifold germs $(M,a)$
and $(M,a')$ are CR-equivalent if and only if $a'\in\GL(F)(a)$.

\Proof Suppose that $g:(M,a)\to(M,a')$ is an isomorphism of
CR-manifold germs. From $\7g'=\7g_{-1}'\oplus\,\RR\delta$ for
$\7g':=\hol(M',a')$ Proposition \ruf{ZB}.iii implies that $g$ is
represented by a linear transformation in $\GL(V)$ that we also denote
by $g$. But then $g(F)\subset S$ with $S$ defined in \Ruf{NP}. Because
$g(F)$ has empty intersection with $H$ we actually have $g(F)\subset
F$. Replacing $g$ by its inverse we get the opposite inclusion, that
is $g\in\GL(F)$.\qed

\KAP{Sechs}{Levi degenerate CR-manifolds associated with an
endomorphism}

The lowest CR-dimension for which there exist homogeneous CR-manifolds
that are Levi degenerate but not holomorphically degenerate is
$2$. The construction recipe below will give, up to local affine
equivalence, all affinely homogeneous conical tube submanifolds of
$\CC^{n}$ with CR-dimension $2$. Indeed, it is based on the following
simple observation: Suppose that $F\subset V:=\RR^{n}$ is a conical
locally linearly homogeneous submanifold of dimension $2$. Denote by
$\7a$ the Lie algebra of all linear vector fields on $V$ that are
tangent to $F$. Then, fixing a point $a\in F$, there exists a
$\phi\in\End(V)$ with $\phi(x)\dd x\in\7a$ and $T_{a}F=\RR
a\oplus\RR\phi(a)$. Therefore, the orbit $H(a)$ under the subgroup
$H:=\{\exp(r\id+t\,\phi):r,t\in\RR\}$ is an (immersed) surface in $V$
having the same germ at $a$ as $F$.

\Joker{CO}{Construction recipe} Throughout this section let $1<d<n$ be
arbitrary integers and $V$ a real vector space of dimension $n$.  Let
furthermore $\phi\in\End(V)$ be a fixed endomorphism and
$\7h\subset\End(V)$ the linear span of all powers $\phi^{k}$ for
$k=0,1,\dots,d-1$. Then $H:=\exp(\7h)\subset\GL(V)$ is an abelian
subgroup, and for given $a\in V$ the orbit $F:=H(a)$ is a cone and an
immersed submanifold of $V$ (not necessarily locally closed in case
$n\ge4$). Furthermore, the tube $M=F+iV\subset E$ is an immersed
CR-submanifold of $E$.

\Joker{}{Cyclic endomorphisms and vectors} A vector $a\in V$ is called
{\sl cyclic} with respect to $\phi\in \End(V)$ if the $\phi^k(a)$,
$k\ge 0$, span $V$. This is equivalent to
$a,\phi(a),\ldots,\phi^{n-1}(a)$ being a basis of $V$. We call
$\phi\in\End(V)$ {\sl cyclic} if it has a cyclic vector and denote
by $\Cyc(V)\subset\End(V)$ the subset of all cyclic endomorphisms. If
$a,b\in V$ both are cyclic vectors of $\phi$ then there exists a
transformation $g\in\RR[\phi]\subset\End(V)$ with $b=g(a)$. But $g$
commutes with every element of the group $H=\exp(\7h)$ and hence maps
the orbit $H(a)$ onto $H(b)$. In particular, the CR-isomorphism type
of $M=H(a)+iV$ only depends on $\phi$ and $d$, but not on the
choice of the cyclic vector $a$. To emphasize this dependence we also
write $M^{\phi,d}$ for $M$ and $F^{\phi,d}$ for $H(a)$, but only if
$\phi$ is cyclic.  We tacitly assume that a choice for a cyclic
vector $a$ has been made. In case $d=2$ we even write $M^{\phi}$ and
$F^{\phi}$ instead of $M^{\phi,2}$ and $F^{\phi,2},$ respectively.
The following proposition shows the relevance of these manifolds in
our discussion.

\Proposition{PR} For $F=H(a)$ and $M=F+iV$ as in \ruf{CO} the
following conditions are equivalent: \0 $\hol(M,a)$ has finite
dimension. \1 $a$ is a cyclic vector of $\phi$.\par\noindent If these
conditions are satisfied, $M$ is a minimal $2$-nondegenerate
homogeneous CR-manifold with CR-dimension $d$ and Levi kernel
$K_{a}M=\CC a$.

\Proof \To12 Condition (i) together with the homogeneity of $M$
implies that $M$ is minimal. Let $W\subset V$ be the linear span of
all vectors $\phi^{k}(a)$, $k\ge0$. Then $H\subset\RR[\phi]$ implies
$H(a)\subset W$ and hence $W=V$ by the minimality of $M$. Therefore,
$a$ is a cyclic vector, and the $\phi^{k}(a)$, $0\le k<d$, form a
basis of the tangent space $T_{a}F$. In particular, $F$ has dimension
$d$, which is also the CR-dimension of $M$.  \nline \To21 Suppose that
$a$ is a cyclic vector of $\phi$.  Lemma \ruf{CR} gives $\RR a\subset
K_{a}F$ since $F$ is a cone in $V$. For the proof of the opposite
inclusion fix an arbitrary $w\in T_{a}F$ with $w\notin\RR a$. Then
$w=\sum_0^m c_j \phi^j(a)$ with $c_m\ne0$ for some $1\le m<d$ and
$\phi^{d-m}(w)\notin T_{a}F$ shows $w\notin K_{a}F$ by Proposition
\ruf{US}. Therefore, $M$ is $2$-nondegenerate by Corollary \ruf{AF}
and $K_{a}M=\CC a$. It remains to show that $M$ is minimal at $a$. But
this immediately follows from Proposition \ruf{PY}.\qed

For the manifolds $M=M^{\phi,d}$ it is possible to compute the Lie
algebras $\7g=\hol(M,a)$ in `most cases'. Clearly, the Lie algebra
$\7g_{0}$ contains the $d$-dimensional Lie algebra $\7h$ (as before,
$\7g_{0}$ is canonically identified with a linear subspace of
$\End(V)\,$). In the Propositions \ruf{UH} and \ruf{HU} we will show
that the equality $\7g=\7g_{-1}\oplus\7h$ holds for all $\phi$ in
`general position'.

\noindent Suppose that $\phi\in\End(V)$ has the cyclic vector $a\in V$
and that the integer $d$ satisfies $1<d<n$. Let
$$
\chi^{}_{\phi}:=\prod_{j=1}^{s}(X-\lambda_{j})^{n_{j}}\;\in\;\RR[X]\,,\quad
n_{j}\ge1\,,\Leqno{CP}
$$ with mutually distinct eigenvalues
$\lambda_{1},\dots,\lambda_{s}\in\CC$ be the characteristic polynomial
of $\phi$.  Let furthermore $\alpha_{1},\dots,\alpha_{n}$ be the
family of all roots of $\chi_{\phi}$, that is, each $\lambda_{j}$
occurs $n_{j}$-times in this string.  As usual,
$\alpha_{1},\dots,\alpha_{n}$ is called an {\sl arithmetic
progression} in $\CC$ if there exists a $\beta\in\CC$ such that, after
a suitable permutation, $\alpha_{j}=\alpha_{1}+(j-1)\beta$ for all
$1\le j\le n$.

\Proposition{UH} Let $M=M^{\phi,d}$ and $\7g:=\hol(M,a)$ for a cyclic
vector $a\in V\cap M$ of $\phi$. Then $\7g_{0}=\7h$ holds if one of the
following conditions is satisfied. \0 $d=2$ and the characteristic
roots $\alpha_{1},\dots,\alpha_{n}$ of $\phi$ do not form an {\sl
arithmetic progression}. \1 $d\ge3$ and $s>d$.\Formend

\noindent For the proof we need several preparations.  To simplify the
notation at various places let us introduce
$$
S:=\{1,\dots,s\},\quad m:=d-1\Steil{and}m_{k}:=n_{k}-1\steil{for
all}k\in S\,.\Leqno{UK}
$$
For every $K\subset S$ define furthermore
$$
\Delta_{K}:= \big\{\beta_{jk}:j\in
S,k\in K\big\}\steil{with}\beta_{jk}:=(\lambda_{j}^{}-\lambda_{k}^{},
\lambda_{j}^{2}- \lambda_{k}^{2},
\dots,\lambda_{j}^{m}-\lambda_{k}^{m})\in\CC^{m}\Leqno{JG}
$$
and denote by $\6K$ the set of all non-empty subsets $K\subset S$ such
that there exist subsets $P\subset Q\subset (S\setminus K)$ with the
following two properties, where the maximum over the empty set here is
defined to be $0\,$: \0 $\sum_{k\in K}n_{k}\;+\;\sum_{q\in Q}n_q
\;\;\ge\;\; d\;+\;\max\re2\{n_{p}-1:p\in P\}\,.$ \1 To every
$\beta\in\Delta_{Q}\!\setminus\Delta_{K}$ there exist uniquely
determined $j\in S$ and $q\in Q$ with $\beta=\beta_{jq}$ such that
$n_{j}<n_{q}$ and $q\in P$.

\noindent Notice that $\6K$ contains every subset $K\subset S$ with
$\;\sum_{k\in K}n_{k}\ge d\;$ (just take $Q=\emptyset$).  In the
following Lemma we show that (i) or (ii) in Proposition \ruf{UH} will
follow from a more general technical condition that will allow us to
give a uniform proof of \ruf{UH} for all $d$.

\Lemma{UL} Suppose that $ s\ge d$ and $\bigcap_{K\in\6K}\!\Delta_{K}
\,=\,\{0\}$ holds. Then one of the conditions (i) and (ii) in
Proposition \ruf{UH} is satisfied.

\Proof {\fett $d=2$:} Assume that $\alpha_{1},\dots,\alpha_{n}$ is not
an arithmetic progression and that there exists a non-zero
$\beta\in\Cap_{K\in\6K}\Delta_{K}$. We claim that
$\lambda_{1},\dots,\lambda_{s}$ is an arithmetic progression:
Otherwise $s\ge3$ and there exists a $k$ with $1<k<s$ such that,
without loss of generality, $\lambda_{j}=\lambda_1+(j-1)\beta$ for all
$1\le j\le k$ and $\lambda_1-\beta\ne\lambda_{r}\ne \lambda_k+\beta$
for all $r>k$. For every $j>k$ the set $\Delta_{\{k,j\}}$ contains the
number $\beta$, that is, there is an $r\in S$ with
$\lambda_{r}=\lambda_{k}+\beta$ or
$\lambda_{r}=\lambda_{j}+\beta$. The first possibility violates our
assumptions. In the second case necessarily $r>k$ must hold since
$r\le k$ would imply $\lambda_{j}=\lambda_1-\beta$ and thus the second
possibility cannot be true for {\it all}\/ $j$ with $k<j\le s$. This
proves that $\lambda_1,...\,,\lambda_s$ is an arithmetic progression
and also $s<n$ by assumption. In particular, $n_j>1$ for some $j\in
S$. Next we claim that $\{1\}\in\6K$ and $\{s\}\in\6K$: To see that
$\{1\}\in\6K,$ let $k\in \{1,...\,, s\}$ be minimal with
$n_{k}>1$. With $P:=Q:=\{k\}\setminus\{1\}$ condition \Ruf{JG}.ii is
fulfilled and the first part of the claim follows. A similar argument
proves also $\{s\}\in\6K$.  But then
$\Delta_{\{1\}}\cap\Delta_{\{s\}}=\{0\}$ gives a contradiction.\nline
{\fett $d\ge3$:} Suppose that $0\ne \beta_{jk}\in
\Cap_{K\in\6K}\Delta_{K}$ and that $s>d$. Then
$L:=S\setminus\{k\}\in\6K$ (with $Q=\emptyset$) and there are $\ell\in
L$, $r\in S$ with $\beta_{jk}=\beta_{r\ell}$. Since $d\ge 3,$ the
equation $\beta_{jk}=\beta_{r\ell}$ implies $k=\ell$, a
contradiction.\phantom{XXXX}\qed

\medskip\noindent{\fett Proof of
Proposition \UH:} It is enough to assume that the assumption of Lemma
\ruf{UL} is satisfied. Consider the decomposition
$E=E_{1}\oplus\cdots\oplus E_{s}$ with $E_{k}$ the kernel of
$(\phi-\lambda_{k})^{n}$ for every $k\in S$. Denote by
$\pi_{k}:E\to E_{k}$ the canonical projection and by
$\epsilon_{k}:E_{k}\to E$ the canonical injection. Then
$a_{k}:=\pi_{k}(a)$ is a cyclic vector for
$\phi_{k}:=\pi_{k}\phi\re1\epsilon_{k}\in\End(E_{k})$. Furthermore, if
we put $a^{j}_{k}:=(\phi_{k}-\lambda_{k})^{j}(a_{k})$ for all $j\ge0$,
then $a^{0}_{k},\dots,a^{m_{k}}_{k}$ is a basis of $E_{k}$. With
$m=d-1$ as defined above let
$\Phi:=(\phi^{1},\dots,\phi^{m})\in\End(V)^{m}$.

\noindent For every real (or complex) vector space $W$ and all tuples
$t=(t_{1},\dots,t_{m})\in\RR^{m}$, $w=(w_{1},\dots,w_{m})\in W^{m}$
let us write as shorthand $t\cd w:=\sum_{j}t_{j}w_{j}$. For every
$t\in\RR^{m}$ the point $e^{t\cd\Phi}(a)$ is contained in $F=M\cap V$.

\noindent Now fix an arbitrary $\mu\in\7g_{0}\subset\End(V)$. Since
the vector field $\mu$ is tangent to $F$, to every
$t\in\RR^{m}$ there exist real coefficients
$r_{0},r_{1},\dots,r_{m}$ with
$$
\mu e^{t\cd\Phi}(a)\;=\;Re^{t\cd\Phi}(a)\Steil{for}R:=
\sum^{m}_{\ell=0}r^{}_{\li2\ell}\phi^{\ell}\,.\Leqno{PT}
$$
Actually, every $r_{\ell}$ has to be considered as a real valued function
on $\RR^{m}$. Put
$$
\mu_{k}:=\pi_{k}\mu\Steil{and}N_{k}:=\big((\phi_{k}-\lambda_{k}),
\dots,(\phi^{m}_{k}-\lambda_{k}^{m})\big)\;\in\;\End(E)^{m}
$$
for all $k\in S$. Applying $\pi_{k}$ to \Ruf{PT}
gives
$$
\eqalign{ R\,e^{t\cd N_{k}}(a_{k})\;&=\;e^{t\cd N_{k}}R(a_{k})=\;
\sum_{j=1}^{s}e^{t\cd\beta_{jk}} \mu_{k}e^{t\cd N_{j}}(a_{j})
\,\quad\hbox{with}\cr
R(a_{k})&=\;\sum^{m_{k}}_{j=0}\rho_{k,j}a^{j}_{k}\Steil{for}
\rho_{k,j}(t):=\sum^{m}_{\ell=j}\textstyle{\ell\choose
j}\lambda^{\ell-j}_{k}r^{}_{\ell}(t)\;\;.\cr }\Leqno{JA}
$$
For every subset $B\subset\CC^{m}$ denote by
$\6F(B)\subset\5C(\RR^{m},\CC)$ the smallest linear subspace
containing all functions $h(t)e^{t\cd\beta}$ with $\beta\in B$ and
$h\in\CC[t]=\CC[t_{1},\dots,t_{m}]$ a polynomial function on
$\RR^{m}$.  Then it is well known that every $f\in\6F(B)$ has a unique
representation $f=\sum_{\beta\in B}f^{[\beta]}e^{t\cd\beta}$ with
$f^{[\beta]}\in\CC[t]$ and $\{\beta\in B:f^{[\beta]}\ne0\}$ finite.
Since $t\cd N_{j}$ is nilpotent, i.e., $e^{\pm t\cd N_j}$ are
polynomial, as a consequence of the identities \Ruf{JA} we get for
every $K\subset S$
$$
\rho_{k,j}\;\in\;\6F(\Delta_{K})\Steil{for every}k\in K\steil{and}0\le
j\le m_{k}\,.\Leqno{JC}
$$
Denote by $\6B$ the set of all subsets $B\subset\CC^{m}$ with
$r_{\ell}\in\6F(B)$ for all $\ell$.

\noindent{\bf Claim 1:} \quad {\sl $\bigcap_{K\in
\6K}\Delta_K\in\6B$, that is, $r_{\ell}\in\CC[t]$ for all
$\ell$.}

\ProofC Fix an arbitrary $K\in\6K$ and let $P\subset Q$ be as in the
definition of $\6K$. Since $\6B$ is closed under intersections it is
enough to show $\Delta_{K}\in\6B$. Assume to the contrary that this is
not true.  Consider the linear system of equations for the $r_{\ell}$,
compare also \ruf{HA},
$$
\rho_{k,j}=\sum^{m}_{\ell=j}\textstyle{\ell\choose
j}\lambda^{\ell-j}_{k}r^{}_{\ell}\;\in\;\6F(\Delta_{K\cup
Q})\Steil{for all}k\in K\cup Q\steil{and}0\le j\le
m_{k}\,.\Leqno{KZ}
$$
The coefficient matrix is of generalized Vandermonde type and hence
has rank $d=m+1$ since by the definition of $\6K$ the number of
equations is at least $d$. This implies
$r_{\ell}\in^{}\6F(\Delta_{K}\cup\Delta_{Q})$ for all $\ell$.  Since
by assumption not all $r_{\ell}$ are in $\6F(\Delta_{K})$ there is a
$\beta\in\Delta_{Q}\!\setminus\Delta_{K}$ such that the
$\beta$-components $r^{[\beta]}_{\ell}\in\CC[t]$ do not vanish for all
$\ell$ simultaneously. By the definition of $\6K$ there are uniquely
determined $p\in S$ and $q\in P$ with $\beta=\beta_{pq}$.  Define
$L:=K\cup Q\setminus \{q\}$.  Since $\beta\not\in\Delta_L$ we get from
\Ruf{JA} the following linear system
$$
\rho^{[\beta]}_{k,j}=\sum^{m}_{\ell=j}\textstyle{\ell\choose
j}\lambda^{\ell-j}_{k}r^{[\beta]}_{\ell}\;=0\,,\steil{where}k
\steil{runs through}L\steil{and}0\le j\le m_{k}\,\,.\Leqno{KX}
$$
By the very definition of $K$, $P$ and $Q$ it follows that the above
linear system consists of at least $d-1$ equations. Consequently we
can write it in the form:
$$
\sum^{m-1}_{\ell=j}\textstyle{\ell\choose
j}\lambda^{\ell-j}_{k}r^{[\beta]}_{\ell}\;\in\CC
r^{[\beta]}_{m}\Steil{for all}k\in L \steil{and}\;0\le j\le m_{k}\;.
$$
Since its coefficient matrix is of generalized Vandermonde type, every
$r^{[\beta]}_{\ell}$ is a complex multiple of $r^{[\beta]}_{m}$ and,
in particular, $r^{[\beta]}_{m}\ne0$. We claim
$\rho^{[\beta]}_{q,0}\ne0$. Indeed, otherwise we could add the
equation $\rho^{[\beta]}_{q,0}=0$ to the linear system \Ruf{KX}, which
then has a coefficient matrix of generalized Vandermonde type with
rank $d$ contradicting $r^{[\beta]}_{m}\ne0$.  Denote by $D$ the
degree of
$\rho^{[\beta]}_{q,0}=\sum_{\ell}\lambda^{\ell}_{q}r^{[\beta]}_{\ell}$. Then
$D=\deg \rho^{[\beta]}_{q,0}=\deg r_{m}^{[\beta]}\ge 0$ and all
$r_{\ell}$ and $\rho^{[\beta]}_{q,j}$ have degree $\le D$. The
equations \Ruf{JA} imply (replace $k$ by $q$, form the
$\beta$-components for $\beta:=\beta_{pq}$ and carry out the
multiplication with $e^{t\cd N_{q}}$)
$$
\mu_{q}e^{t\cd N_{p}}\re1(a_{p})\;=\;e^{t\cd
N_{q}}\sum^{m_{q}}_{j=0}\rho^{[\beta]}_{q,j}a^{j}_{q}=\sum^{m_{q}}_{j=0}
\big(f_{j}\rho^{[\beta]}_{q,0}+g_{j}\big)a^{j}_{q}\Leqno{KI}
$$
for certain polynomials $f_{j},g_{j}\in\CC[t]$ with
$\deg(g_{j})<\deg(f_{j})=j$. Comparing degrees on both sides in \Ruf{KI}
and keeping in mind $n_{p}<n_{q}$ (by the definition of $\6K$) we get
$$
D+m_{q}=\deg(f_{m_{q}}\rho^{[\beta]}_{q,m_{q}}+g_{m_{q}})\le
m_{p}<m_{q}\,.
$$
This contradicts $r_m^{[\beta]}\ne 0$ and Claim 1 is proved.\qd

\noindent{\bf Claim 2: \sl Every $r_{\ell}$ is a constant polynomial.}

\ProofC Fix a $k\in S$ with $D:=\deg(\rho_{k,0})=\max_{j\in
S}\deg(\rho_{j,0})$. With $s\ge d$ a Vandermonde argument applied to
the linear system $\rho_{j,0}=\sum_{\ell}\lambda^{\ell}_{j}r_{\ell}$,
$\,j\in S$, gives that every $r_{\ell}$ has degree $\le D$. As in
\Ruf{KI} we have
$$
\mu_{k}e^{t\cd N_{k}}\re1(a_{k})\;=\;e^{t\cd
N_{k}}\sum^{m_{k}}_{j=0}\rho_{k,j}a^{j}_{k}=\sum^{m_{k}}_{j=0}
\big(f_{j}\rho_{k,0}+g_{j}\big)a^{j}_{k}\Leqno{KY}
$$ for polynomials $f_{j},g_{j}\in\CC[t]$ with
$\deg(g_{j})<\deg(f_{j})=j$.  All coefficient polynomials in \Ruf{KY}
in front of the $a^{j}_{k}$ have degree $\le m_{k}$, that is
$D+m_{k}\le m_{k}$ and hence $D\le0$. This proves Claim 2.\qd
\nline The proof of \ruf{UH} now is complete: Indeed, since
$F$ contains a basis of $V$, the endomorphism $\mu$ is uniquely
determined by the function tuple $(r_{\ell})$, that is, $\dim\7g_{0}\le
d=\dim \7h $.\qed

\smallskip \Remark{HA} For given tuples
$\lambda_{1},\dots,\lambda_{s}\in\CC$ and $\,n_{1},\dots,n_{s}\in\NN$
with $n_{k}\ge1$ and $n:=\sum\!n_{k}$ let $L:=\{(k,j):1\le k\le
s,\,0\le j<n_{k}\}$ endowed with the lexicographic order. Then it can
be seen that the following $n\Times n$-matrix (every entry with $\ell<
j$ is zero)
$$
\Big(\textstyle{\ell\choose
j}\lambda^{\ell-j}_{k}\Big)_{0\le\ell<n,\re1(k,j)\in L}\Steil{has
determinant}\displaystyle\prod\limits_{p<q}
(\lambda_{q}-\lambda_{p})^{n_{p}n_{q}}\;.
$$

Notice from the proof of \ruf{UH} that the condition
$\bigcap_{K\in\6K}\!\Delta_{K} \,=\,\{0\}$ guarantees that every
$\mu\in\7g_{0}$ leaves every generalized eigenspace $E_{k}$ of $\phi$
invariant, while $s\ge d$ guarantees that every such $\mu$ actually is
in $\7h$. In the next section we will see that the condition (i) in
Proposition \ruf{UH} for $d=2$ is optimal, compare Proposition \ruf{ZR}.

\Proposition{HU} Let $M=M^{\phi,d}$ and assume that $\7g_{0}=\7h$ for
$\7g=\hol(M,a)$. Then $\7g=\aff(M,a)$ and $\aut(M,a)=0$.

\Proof For the proof of $\7g=\aff(M,a)$ it is enough to show
$\7g_{1}=0$ by Proposition \ruf{DC}.iv. This is more easily done in
the more general complex setting, compare the following Lemma
\ruf{UF}. Finally, counting dimensions yields $\aut(M,a)=0$.\qed

\medskip\noindent It remains to show the next Lemma. As before,we
identify the spaces $\End(E)$ and $\7P_{0}$, see \Ruf{HJ}.\vskip-30pt

\Lemma{UF} Let $\phi\in\Cyc(E)$ be an arbitrary cyclic
endomorphisms. For given integer $d<n=\dim E$ let furthermore $\7h$ be
the complex linear span of all powers $\phi^{j}$, $0\le j<d$. Then
$$
\big\{\xi\in\7P_{1}:[\7P_{-1},\xi]\subset\7h\big\}\;=0\;\,.
$$
\Proof We identify $E$ with $\CC^{n}$ in such a way that the
matrix $\Phi$ of $\phi$ is in Jordan normal form, more precisely:
$\Phi$ is a block diagonal matrix with Jordan blocks
$J_{1},\dots,J_{s}$, where each block $J_{l}$ is lower triangular, has
the eigenvalue $\lambda_{l}$ on its main diagonal and is of size
$n_{l}\Times n_{l}$ for some $n_{l}\ge1$. We also introduce an
equivalence relation on $\{1,\dots,n\}$ in the following way: Put
$j\sim k$ if the $j^{\th}$ row and the $k^{\th}$ column in $\Phi$
intersect in one of the Jordan blocks. \nline Now suppose that there
exists a non-zero vector field $\xi\in \7P_{1}$ with
$[\7P_{-1},\xi]\subset\7h$. This $\xi$ has a unique representation
$$
\xi=\sum_{j,k,p=1}^{n}c_{p}^{jk}z_{j}z_{k}\dd{z_{p}}\Steil{with}
c_{p}^{jk}=c_{p}^{kj}\in\CC\,.
$$
For every $j\le n$ the vector field $\dd z_{j}$ is contained in
$\7P_{-1}$. Therefore
$$
\xi_{j}:={1\over2}\big[\dd{z_{j}},\xi\big]=\sum_{k,p=1}^{n}
c^{jk}_{p}z_{k} \dd{z_{p}}
$$
is contained in $\7h$ implying $c_{p}^{jk}=0$ if $k\not\sim p$ and, by
symmetry, $c_{p}^{jk}=0$ if $j\not\sim p$. This implies
$$
\xi_{j}=\sum_{p\sim j}\sum_{k\sim j}c^{jk}_{p}z_{k}
\dd{z_{p}}\,.\Leqno{VZ}
$$ By assumption $\xi\ne0$ and therefore $\xi_{r}\ne0$ for some $r\le
n$. Without loss of generality we may assume $r\sim1$. Replacing
$\phi$ by $(\phi-\lambda_{1})$ we also may assume $\lambda_{1}=0$.
Put $b:=n_{1}$, that is, the Jordan block $J_{1}$ has size $b\times b$
and is nilpotent. Define for all $1\le p,k\le b$
$$
\eta_{p}:=z_{1}\dd{z_{p}}+z_{2}\dd{z_{p+1}}+\dots+z_{b+1-p}\dd{z_{b}}
\,,
$$
$$
\7D_{k}:=\sum_{j=k}^{b}\CC\re1\eta_{j}\Steil{and}\psi_{k}:=
\phi^{k-1}\prod_{l=2}^{s}(\phi-\lambda_{l})^{n_{l}}\;
\in\;\End(E)\,.
$$
From \Ruf{VZ} we know
$\xi_{j}\in\7D_{1}\cap\7h$ for all $j\le b$. In particular,
$$ \xi_{j}=\sum_{p=1}^{b}c^{j1}_{p}\eta_{p}\,,\qquad c_{p}^{jk}=
\cases{c^{j1}_{p+1-k}&if $k\le p$\cr0&otherwise\cr}\qquad\hbox{and hence}
$$
$$
\xi_{j}=c^{11}_{1}\eta^{}_{j}+c^{11}_{2} \eta^{}_{j+1}+
\dots+c^{11}_{b+1-j}\eta^{}_{b}
$$ for all $j$ due to the symmetry of $c^{jk}_{p}$ in the upper
indices.  As a consequence there exists a minimal $q\le b$ with
$q\ge1$ and $c^{11}_{q}\ne0$. But then
$\xi_{b+1-q}=c^{11}_{q}\eta^{}_{b}$ implies
$\eta^{}_{b}\in\7D_{b}\cap\7h$. We show that this cannot be
true:\nline Since $\psi_{1}$ as polynomial in $\phi$ has constant term
$\lambda_{2}\lambda_{3}\dots\lambda_{s}\ne0$, we get that $\psi_{k}$
spans $\7D_{k}$ over $\7D_{k+1}$ for every $k\ge1$. This implies
$\eta_{b}=\sum_{j=1}^{n}e_{j}\phi^{j-1}$ for suitable real
coefficients $e_{j}$ with $e_{n}\ne0$ and thus
$\eta_{b}\notin\7h$.\phantom{XXX}\qed

\bigskip For the application of Proposition \ruf{UQ} to manifolds of
the type $M=M^{\phi,d}$ it is necessary to know when $M$ is simply
connected and when $\Aut(M,a)$ is the trivial group. \nline A
sufficient condition for $M^{\phi,d}$ (and $F^{\phi,d}=H\Kl a$) to be
simply connected is the following: {\sl There exist eigenvalues
$\lambda_{1},\dots,\lambda_{d}$ of $\phi$ such that ${\det(A+\overline
A)}\ne0$, where $A=(\lambda_{j}^{k-1})_{1\le j,k\le d}$ is the
corresponding Vandermonde matrix}. \nline To get a partial answer to
the second question consider to given $\epsilon\in\RR^{*}$ a $g\in\GL(V)$
with $g(a)=a$ and $g\phi g^{-1}=\epsilon\phi$. From
$g\phi^{k}g^{-1}=(\epsilon\phi)^{k}$ we get
$g(\phi^{k}(a))=\epsilon^{k}\phi^{k}(a)$ for all $k\ge0$, that is, $g$
is uniquely determined in $\GL(V)$ by the above assumptions since
$\phi$ is cyclic. A further consequence is
$g\re1\exp(t\phi^{k})\re2g^{-1}=\exp(t\epsilon^{k}\phi^{k})$ for all
$t\in\RR$ and $k\ge0$. This means $gHg^{-1}=H$ for the group
$H=\exp(\7h)$. But then $g(F)=F$ for $F=H(a)$ and consequently
$g\in\Aut(M,a)$.

Notice that we always may assume without loss of generality that
$\phi\in\Cyc(V)$ has trace $0$ (otherwise replace $\phi$ by
$(\phi-c\,\id)\in\End(V)$ for $c:=n^{-1}\tr(\phi)$, since this
procedure does not change the algebra $\7h)$.

\Lemma{KR} Suppose that $2d\le n+1$ for $n=\dim V $ and that $\phi$
is trace-free. Suppose in addition that $\aut(M,a)=0$ holds for the
cyclic vector $a\in V$ and $M:=M^{\phi,d}$. Then
$$\Aut(M,a)=\big\{g\in\GL(V):g(a)=a\steil{and}g\phi
g^{-1}=\pm\phi\big\}\,.$$ In particular, $\Aut(M,a)$ always has order
$\le 2$ and is trivial, for instance, if the spectrum of $\phi$ in
$\CC$ is not symmetric with respect to the origin of $\CC$.

\Proof By Proposition \ruf{LO} the assumption $\aut(M,a)=0$ implies
$\7g_{1}=0$ and $\7g_{0}=\7h$. As a consequence of Proposition
\ruf{PS} therefore
$\Aut(M,a)=\big\{g\in\GL(V):g(a)=a\steil{and}g\7hg^{-1}= \7h\big\}$
holds. Let $g\in\Aut(M,a)$ be an arbitrary automorphism. Then $g\phi
g^{-1}=\sum^{m}_{j=0}c_{j}\phi^{j}$ for real coefficients $c_{j}$ and
$m:=d-1$. We show by induction on $k$ that $c_{j}=0$ holds for all
$j>m/k$ and all $1\le k\le m$. For $k=1$ this is obvious. So fix a
$k>1$ with $k\le m$. By induction hypothesis $g\phi^{k}
g^{-1}=\big(\sum c_{j}\phi^{j}\big)^{k}
=\sum^{2m}_{\ell=0}e_{\ell}\phi^{\ell}\in\7h$ with real coefficients
$e_{\ell}$. Since $2m<n$ by assumption, we must have $e_{\ell}=0$ for
all $\ell>m$, that is $c_{j}=0$ for all $j>m/k$. For $k=m$ this
implies $c_{j}=0$ for all $j>1$. Taking traces finally gives $g\phi
g^{-1}=\epsilon\phi$ for $\epsilon:=c_{1}$. By the above example
$\phi$ cannot be nilpotent, that is, $\phi$ has non-zero spectrum in
$\CC$. Since this spectrum is invariant under multiplication with
$\epsilon$ necessarily $\epsilon=\pm1$ holds.\qed

\KAP{Sieben}{Homogeneous 2-nondegenerate manifolds of CR-dimension 2}

In this section we specialize to homogeneous tube manifolds $M=F+iV$
in $E=V\oplus iV$ of CR-dimension $2$, that is, where $F\subset V$ is
a surface of dimension $2$. We begin with manifolds of type
$M^{\phi}=M^{\phi,2}$ that are obtained by the construction recipe
\ruf{CO}. Since Propositions \ruf{UH} and \ruf{HU} of the preceding
section do not cover the case where the characteristic roots
$\alpha_{1},\dots,\alpha_{n}$ of $\phi\in\Cyc(V)$, $n=\dim V$, form an
arithmetic progression, let us discuss this case first:\nline Possibly
after replacing $\phi$ by $\phi-r\id$ with an appropriately chosen
constant $r$ we may assume without loss of generality that $\phi$ is
trace-free. Since multiplication of $\phi$ by any non-zero real number
does not change the algebra $\7h=\RR\id\oplus\RR\phi$, there are
essentially 3 different cases for $\phi$ with characteristic roots
forming an arithmetic progression -- either $\phi$ is nilpotent or
$\phi$ has pairwise different characteristic roots in $\RR$ or in
$i\RR$.  Consider the manifold $M:=\5M^{2,n-2}=\5F^{2,n-2}+iV$ from
Example \ruf{JP}. As already remarked in \ruf{JP} the conformal
subgroup $H^{\i}:=\RRp\cd\SO(2)\subset\GL(2,\RR)$ acts transitively on
$\5F^{2,n-2}$. The corresponding Lie algebra is
$\7h_{\i}:=\RR\id\oplus\RR\phi^{}_{\i}$, where
$\phi_{\i}:=\xi^{1,-1}-\xi^{-1,1}\in\7g_{0}$ is a trace-free
semisimple endomorphism with eigenvalues in $i\RR$, see \Ruf{AQ} for
the notation. Next consider the subgroup
$H^{\r}:=\RRp\cd\SO(1,1)\subset\GL(2,\RR)$ with Lie algebra
$\7h^{\r}:=\RR\id\oplus\RR\phi^{}_{\r}$, where
$\phi_{\r}:=\xi^{1,-1}+\xi^{-1,1}\in\7g_{0}$ is a trace-free
semisimple endomorphism with real eigenvalues. In particular,
$H^{\r}(a)$ is open in $H^{\i}(a)$. Finally, consider the solvable
subgroup
$$
H^{\z}:=\Big\{\Matrix\alpha0\beta\alpha\in\GL(2,\RR):\alpha>0\Big\}
\steil{with Lie algebra}\7h^{\z}:=\RR\id\oplus\RR\phi^{}_{\z}\,,
$$
where $\phi_{\z}:=\xi^{1,-1}$ only has zero eigenvalues. Again, the
orbit $H^{\z}(a)$ is open in $H^{\i}(a)$. This implies that the
manifolds $M^{\phi_{\r}}$ and $M^{\phi_{\z}}$ are open subsets of
$M^{\phi_{i}}=\5M^{2,n-2}$ and hence that all three of them are
locally CR-equivalent. Since $\5M^{2,n-2}$ is minimal as CR-manifold
the endomorphisms $\phi_{\i},\phi_{\r}$ and $\phi_{\z}$ have $a$ as
cyclic vector by Proposition \ruf{PR} (what also easily can be
verified). By construction, the characteristic roots of the
endomorphisms $\phi_{\i}$, $\phi_{\r}$ and $\phi_{\z}$ form an
arithmetic progression and represent the three types with imaginary,
real and zero characteristic roots.  The estimate $\dim \7g_{0}\ge
4>\dim \7h^\psi$ for $\psi=\phi_{\i},\phi_r,\phi_n$ together with
Propositions \ruf{UH} and \ruf{HU} implies that the characteristic
roots of all 3 endomorphisms form an arithmetic progression.  Summing
up we have proved  the following result.

\Proposition{IO} Let $\phi,\phi'\in\End(V)$ be endomorphisms with cyclic
vectors $a,a'\in V$. Assume that for both endomorphisms the families
of characteristic roots $\alpha_{1},\dots,\alpha_{n}$ and
$\alpha'_{1},\dots,\alpha'_{n}$ form an arithmetic progression. Then
the germs $(M^{\phi},a)$ and $(M^{\phi'},a')$ are
CR-equivalent.\Formend

\medskip We can use Proposition \ruf{IO} to get explicit global
equations for every $M^{\phi}$ where the characteristic roots of
$\phi\in\Cyc(V)$ form an arithmetic progression: Indeed, we may take
$\phi:=\xi^{1,-1}$ from \Ruf{AQ} on $\CC^{m+1}$ with coordinates
$(z_{0},z_{1},\dots,z_{m})$ and $a:=(1,0,\dots,0)$. Then $\phi$ is
nilpotent and
$S:=^{}\exp(\RR\phi)(a)=\{(1,t,t^{2},\dots,t^{m}):t\in\RR\}$. Consequently,
$M^{\phi}=F+iV$ is the tube over the cone $F$ generated by $S$ (that
is $F=\RRp\!\cd S$). As a consequence, $F$ is an open piece of the
algebraic surface given by the following explicit system of quadratic
equations on $\RR^{m+1}$ with coordinates $(x_{0},x_{1},\dots,x_{m})$
$$
x_{0}x_{j+1}=x_{1}x_{j}\Steil{for}0<j<m\,.\Leqno{IP}
$$
This can be reformulated also in the following slightly different
form: Let $C:=\{(t,t^{2},\dots,t^{n}):t\in\RR\}$ in $\RR^{n}$ be the
{\sl twisted $n$-ic} (also called twisted cubic, quartic etc, see
\Lit{HART} for interesting properties of these curves). Then the cone
$\RRp\!\cd C$ generated by $C$ is a nonsingular surface outside the
origin and the corresponding tube manifold is locally CR-equivalent to
$M^{\phi}$ with $\phi$ as in \ruf{IO}. The twisted $n$-ic will also
show up in another type of examples, compare \ruf{KU}.

\medskip Next we extend Propositions \ruf{UH} and \ruf{HU} to the case
$d=2$ where the characteristic roots of $\phi$ do form an arithmetic
progression.  Recall that for the light cone tube $\5M\cong\5M^{2,1}$
the Lie algebra $\7g=\hol(\5M,a)$ is isomorphic to $\so(2,3)$, compare
\Lit{KAZT}, \Lit{FEKA}. In particular, $\7g_{0}\cong\gl(2,\RR)$ and
$\dim\7g_{1}=3$ in this case.

\Proposition{ZR} Assume that the characteristic roots of
$\phi\in\Cyc(V)$ form an arithmetic progression. Then for
$M=M^{\phi}$, $\,n=\dim V$ and $\7g=\hol(M,a)$ the following
properties hold. \0 $\7g_{0}$ is isomorphic to $\gl(2,\RR)$ and hence
has dimension $4$. \1 $\7g=\aff(M,a)$ in case $n\ge4$. In particular,
$\dim\7g=n+4$ in this case.

\Proof ad (i): We may assume that $\lambda_{j}=(1-n)/2+(j-1)$ for
all $j\in S$ using the notation in \Ruf{UK} and \Ruf{JG}. Then
$\bigcap_{K\in\6K}\Delta_{K}=\{-1,0,1\}$. Solving \Ruf{PT} for
$r_{0},r_{1}$ gives that all pairs
$$
\eqalign{r_{0}&=(n-1)(ue^{-t}+v_{0}-we^{t})\cr
r_{1}&=2(ue^{-t}+v_{1}+we^{t})\,,\cr}\Leqno{AP}
$$ $u,v_{0},v_{1},w\in\RR$ arbitrary, form the solution space. This
implies $\dim\7g_{0}=4$. Since $M$ is locally CR-equivalent to
$\5M^{2,n-2}$ the Lie algebra $\7g_{0}$ contains a copy of
$\gl(2,\RR)$, that is $\7g_{0}\cong\gl(2,\RR)$.\nline ad (ii): Let
$n\ge4$. We may assume that for $m=n-1$ the Lie algebra $\7g_{0}$ is
the linear span of the vector fields \Ruf{AQ}. For every
$\nu\in\ZZ^{2}$ let
$\7g^{\nu}:=\{\xi\in\7g:[\zeta_{j},\xi]=\nu_{j}\xi\steil{for}j=1,2\}\,.$
Then $\7g^{\nu}\subset\7g_{k}$ for $k=(\nu^{}_{1}+\nu^{}_{2})/m$ and
$$
\7g=\bigoplus_{\nu\in\ZZ^{2}}\7g^{\nu}\Steil{with}[\7g^{\nu},
\7g^{\mu}]\subset\7g^{\nu+\mu}\,,
$$ compare also (3.5) in \Lit{FEKA}.
Clearly $i\dd{z_{k}}\in\7g^{-k,k-m}\steil{for all}0\le k\le m$.
Because of Proposition \ruf{DC}.iv it is enough to show $\7g_{1}=0$.
Assume to the contrary that there exists a non-zero $\xi\in\7g_{1}$.
Then we may assume without loss of generality that $\xi\in\7g^{k,m-k}$
for some $k\in\ZZ$. Let $c$ be the cardinality of $\{0\le j\le
m:[\dd{z_{j}},\xi]\ne0\}$. From
$\7g_{0}=\7g^{-1,1}\oplus\7g^{0,0}\oplus\7g^{1,-1}$ and
$[i\dd{z_{j}},\xi]\in\7g^{k-j,j-k}$ we see $c\le3$. Assume $c=3$,
which implies $0<k<m$. From $[\dd{z_{j}},\xi]\ne0$ for $j=k\pm1$ and
the special form of $\xi^{-1,1}$, $\xi^{1,-1}$ in \Ruf{AQ} we see that
$\xi$ must depend on all $n$ variables, a contradiction to $n>3$.  But
$c\le2$ also gives a contradiction since in all spaces $\7g^{-1,1}$,
$\7g^{0,0}$, $\7g^{1,-1}$ every non-zero vector field must depend on
at least $n-2$ variables.\qed

\medskip Recall that $\Aut(M)$ is the group of all global
CR-automorphisms and $\Aff(M)$ is the subgroup of all of affine
transformations of $M$.

\Proposition{GO} Let $\phi\in\End(V)$ be a trace-free cyclic
endomorphism. Then the groups $\Aut(M)$ and $\Aff(M)$
coincide. Furthermore, with $n=\dim V$ the following dimension
estimates hold.\0 $\dim\Aut(M)=n+4$ if the characteristic roots of
$\phi$ are pairwise distinct and form an arithmetic progression in
$i\RR$.\1 $\dim\Aut(M)=n+3$ if $\phi$ is nilpotent.\1
$\dim\Aut(M)=n+2$ in all other cases.

\Proof Let $F:=M\cap V$ and denote by $\7a\subset\7g=\hol(M,a)$ the
Lie algebra of $\Aut(M)$. Since $\7a$ contains the Euler vector field
we have $\7a=\7a_{-1}\oplus\7a_{0}\oplus \7a_{1}$ for
$\7a_{j}:=\7a\cap\7g_{j}$. We first determine $\dim\7a_{0}$. Because
of $\7h:=\RR\id\oplus\RR\phi\subset\7a_{0}$ we have
$\dim\7a_{0}\ge2$.\nline In case (i) $M$ is CR-equivalent to
$\5M^{2,n-2}$, compare Example \ruf{JP}. Therefore $\GL(2,\RR)$ acts
transitively on $F$ and $\dim\7a_{0}=4$ by Proposition \ruf{ZR}.\nline
Next consider the case (ii), that is, $\phi$ is nilpotent. Then
$\7a_{0}$ consists of all $\xi\in\7g_{0}\subset\End(V)$ with
$\xi(c)\in\RR c$ for all $c\in\partial F:=\overline F\setminus F$,
where $\overline F$ is the closure of $F$ in $V$. We may assume that
$a=(1,0,\dots,0)\in\RR^{n}$ and $\phi=\xi^{1,-1}$ in the notation of
\Ruf{AQ}. This implies
$F=\{e^{s}(1,t,t^{2},\dots,t^{n-1})\in\RR^{n}:s,t\in\RR\}$ and hence
$\partial F=\RR c$ for $c:=(0,\dots,0,1)$. Therefore, $\7a_{0}$ is the
linear span of $\zeta_{1}$, $\zeta_{2}$ and $\xi^{1,-1}$ and
$\dim\7a_{0}=3$ in this situation.\nline Next consider the case
$V=V_{1}\oplus V_{2}\oplus\dots\oplus V_{n}$ where every $V_{j}$ is
the $\big((1-n)/2+(j-1)\big)$ eigenspace of $\phi$, that is, the
characteristic roots of $\phi$ form an arithmetic progression in
$\RR$. Here $\partial F=V_{1}\cup V_{n}$ is easily verified.  The
vector fields in $\7g_{0}$ are characterized by the function tuples
$(r_{0},r_{1})$ in \Ruf{AP}. The condition $\xi(V_{j})\subset V_{j}$
for $j=1,n$ implies that $r_{0},r_{1}$ are constant for every
$\xi\in\7a_{0}$, that is $\7a_{0}=\7h$. On the other hand, if the
characteristic roots of $\phi$ do not form an arithmetic progression,
then also $\7g_{0}=\7h$ by Proposition \ruf{UH}.i, that is,
$\dim\7a_{0}=2$ always holds in case (iii).\nline Next we show
$\7a_{1}=0$ in all cases: For $n\ge4$ this follows from $\7g_{1}=0$,
see Proposition \ruf{ZR}.ii. In case (iii) we have $\7a_{0}=\7h$ and
the claim follows with Lemma \ruf{UF}. Therefore we only have to
consider the cases (i) and (ii) for $n=3$. In case (i) $M$ is the
future light cone tube $\5M$ and $\Aut(\5M)=\Aff(\5M)$ follows as a
special case of Proposition 6.9 in \Lit{KAZT}. In case (ii) $M$ is a
proper domain in $\5M$: We realize $\5M=\5F+i\RR^{3}$ in
$\CC^{3}$ with coordinates $(z_{0},z_{1},z_{2})$ as
$\5F=\{x\in\RR^{3}:x_{0}x_{2}=x^{2}_{1}\steil{and}x_{0}+x_{2}>0\}$.
Then $\hol(\5M)$ is the linear span of the vector fields (3.5) and
(3.7) in \Lit{FEKA}. We may assume without loss of generality that
$\phi=\xi^{-1,1}=2z_{1}\dd{z_{0}}+z_{2}\dd{z_{1}}$. This implies
$\5M\setminus M=\RRp\!\cd c+i\RR^{3}$ for $c:=(1,0,0)\in\5F$.
We know already that $\7a$ is the linear span of the vector fields
$\zeta_{1}$, $\zeta_{1}$ and $\xi^{-1,1}$. From Figure 1 and (3.7) in
\Lit{FEKA} we therefore derive that either $\7a_{1}=0$ or
$\7a_{1}=\RR\xi^{0,2}$ for
$\xi^{0,2}:=iz^{2}_{1}\dd{z_{0}}+iz_{1}z_{2}
\dd{z_{1}}+iz^{2}_{2}\dd{z_{2}}$. The latter possibility cannot occur
since $\xi^{0,2}$ is not tangent to $\5M\setminus M$: Check, for
instance, the point $(1,i,0)$. This proves $\7a=\7g_{-1}\oplus\7a_{0}$
in all cases. As in the proof of Proposition \ruf{ZB}.iii it is shown
that this implies $\Aut(M)=\Aff(M)$. The above dimension estimates for
$\7a_{0}$ imply the dimension estimates in (i) -- (iii).\qed

\medskip Next we solve the local as well as the global CR-equivalence
problem for all manifolds $M^{\phi}$. Because of Proposition \ruf{IO}
in the local situation only the case has to be considered where the
characteristic roots of $\phi$ do not form an arithmetic
progression. Recall that without loss of generality we always may
assume that $\phi$ is trace-free.

\Proposition{KO} Let $\phi,\phi'\in\End(V)$ be trace-free cyclic
endomorphisms with characteristic roots $\alpha_{1},\dots,\alpha_{n}$
and $\alpha'_{1},\dots,\alpha'_{n}$ respectively. Suppose that the
$\alpha_{1},\dots,\alpha_{n}$ do not form an arithmetic
progression. Then for given cyclic vectors $a,a'\in V$ and
corresponding $M=M^{\phi}$, $M'=M^{\phi'}$ the Lie algebras
$\hol(M,a)$ and $\aff(M,a)$ coincide. Furthermore, the following
conditions are equivalent. \0 The Lie algebras $\hol(M,a)$ and
$\hol(M',a)$ are isomorphic . \1 The germs $(M,a)$ and $(M',a')$ are
CR-equivalent. \1 $g\phi'g^{-1}=r\phi$ for suitable $g\in\GL(V)$ and
$r\in\RR^{*}$. \1 There exists a permutation $\pi\in\7S_{n}$ and an
$r\in\RR^{*}$ with $\alpha'_{j}=r\alpha_{\pi(j)}$ for all $j$.

\Proof $\7g:=\hol(M,a)=\aff(M,a)$ and $\dim\7g=n+2$ follows from
Propositions \ruf{UH} and \ruf{HU}.\nline \To12 With $\7g$ also
$\7g':=\hol(M',a')$ has dimension $n+2$. Therefore
also $\alpha'_{1},\dots,\alpha'_{n}$ do not form an
arithmetic progression, see Proposition \ruf{ZR}. This implies
$\7g_{0}=\7h$, $\7g'_{0}=\7h'$ and (ii) follows with Proposition
\ruf{PQ}. \nline\To23 Let $g$ be a CR-isomorphism
$(M,a)\to(M',a')$. Then $g\in\GL(V)$ as a consequence of Proposition
\ruf{ZB}.iii. and clearly $g\7h g^{-1}=\7h'$. Since
$\RR\phi'\subset\7h'$ is precisely the subset of all trace-free
endomorphisms (iv) follows. \nline The remaining implications are easy
to check and left to the reader.\qed

\ifarx\eject\fi

\Proposition{BP} Let $\phi,\phi'\in\End(V)$ be trace-free cyclic
endomorphisms and $M:=M^{\phi}$, $M':=M^{\phi'}$. Then the following
conditions are equivalent. \0 The groups $\Aff(M)$ and $\Aff(M')$ are
isomorphic. \1 $M$ and $M'$ are globally CR-equivalent. \1
$g\phi'g^{-1}=r\phi$ for suitable $g\in\GL(V)$ and $r\in\RR^{*}$.
\Proof \To13 Suppose that (i) holds. Because of Proposition \ruf{GO}
we may assume $\dim\Aut(M)=n+2$ without loss of generality.  Then
(iii) follows from Proposition \ruf{KO} if at least for one of the
$\phi,\phi'$ the characteristic roots do not form an arithmetic
progression. In the remaining cases the claim follows from Proposition
\ruf{GO} since for both endomorphisms the characteristic roots form an
arithmetic progression in $\RR$ and are pairwise distinct.\nline\To32
$\Rightarrow$ (i) is obvious since $\Aut(M)=\Aff(M)$ and
$\Aut(M')=\Aff(M')$ by Proposition \ruf{GO}.\qed

\ifarx\bigskip\else\medskip\fi
\Joker{MO}{Some moduli spaces} For fixed $n=\dim V$ let $\6M$ be
the space of all global CR-equivalence classes $[M^{\phi}]$ of
manifolds $M^{\phi}$ with $\phi\in\Cyc(V)$, that is, every
$[M^{\phi}]$ is the set of all $M^{\psi}$ that are globally
CR-equivalent to $M^{\phi}$. Furthermore put
$$\Phi:=\{\phi\in\Cyc(V):\tr(\phi)=0\}\Steil{and}\6M^{*}:=
\big\{[M^{\phi}]\in\6M: \phi\in\Phi,\,\phi^{n}\ne0\big\}\,.$$ The
reductive group $\RR^{*}\Times\GL(V)$ acts on $\End(V)$ by
$\phi\mapsto rg\phi g^{-1}$ for every $(r,g)\in\RR^{*}\Times\GL(V)$
and leaves the cone $\Phi$ invariant. By Proposition \ruf{BP} $\6M$
can be identified as set with the quotient
$\Phi/(\RR^{*}\Times\GL(V))$.  This quotient can be built in several
steps: For every $j$ let $\sigma_{j}(\phi)\in\RR$ be the $j^{\th}$
elementary symmetric function in $n$ variables evaluated on the
characteristic roots of $\phi$, that is,
$$X^{n}+\sum^{n}_{j=2}(-1)^{j}\sigma_{j}(\phi)X^{n-j}\;\in\;\RR[X]$$
is the characteristic polynomial of $\phi\in\Phi$. Let $W:=\RR^{n-1}$
with coordinates $(x_{2},\dots,x_{n})$ and denote by $\sigma:\Phi\to
W$ the mapping given by
$\phi\mapsto(\sigma_{2}(\phi),\dots,\sigma_{n}(\phi))$. Since every
real polynomial factors in a product of linear and quadratic real
polynomials, the map $\sigma$ is surjective and $\6M$ can be
canonically identified as set with the quotient $W/\RR^{*}$, where
$\RR^{*}$ acts on $W$ by
$(x_{2},x_{3},\dots,x_{n})\mapsto(t^{2}x_{2},t^{3}x_{3},\dots,t^{n}x_{n})$
for every $t\in\RR^{*}$. The subgroup $\{\pm1\}\subset\RR^{*}$ leaves
the sphere $S^{n-2}=\{x\in W:\sum x_{j}^{2}=1\}$ invariant and
$\6M^{*}$ can be identified with the quotient
$Q^{n-2}:=S^{n-2}/\{\pm1\}$.  In general, $Q^{n-2}$ can be stratified
into a finite number of manifolds. For instance, $Q^{1}$ is a compact
line segment and $Q^{2}$ is homeomorphic to the sphere $S^{2}$. At
this point a word of caution is necessary: We do not give a topology
to $\6M$, the topology on $Q^{n-2}$ only serves for the readers
imagination.

\ifarx\bigskip\fi
Instead of $\6M$ we can also consider the space of {\sl local }
CR-equivalence classes for manifolds of type $M^{\phi}$. By our
results this space is of the form $\6M/\Sim\,$, where the equivalence
relation $\Sim$ on $\6M$ just identifies the $3$ equivalence classes
$[M^{\phi}]\in\6M$ such that the characteristic roots of $\phi$ form
an arithmetic progression. Clearly, $\6M/\Sim=\6M^{*}/\Sim$ can be
obtained by identifying two points in $Q^{n-2}$. In the spacial case
$n=3$ the endpoints of the line segment $Q^{1}$ have to be identified,
that is, $\6M$ can be thought of in this case as the circle
$\overline{\RR}:=\RR\cup\infty$ (without topology) where the point
$\infty$ corresponds to the class represented by the future light cone
tube $\5M$. To be more specific, call in case $n=3$ for every
$\phi\in\Phi$
$$
\mu\big(M^{\phi}\big):=-{\sigma_{2}(\phi)^{3}\over
\sigma_{3}(\phi)^{2}} \;\;\in\;\;\overline{\RR}\Leqno{KW}
$$
the {\sl modulus} of the CR-manifold $M^{\phi}$ (with $t/0:=\infty$
for all $t\in\RR$). It is clear that $M^{\phi}$ and $M^{\psi}$ in case
$n=3$ are locally CR-equivalent if and only if they have the same
modulus. A special meaning has the modulus $\mu_{0}:=27/4$:
For real moduli $>\mu_{0}$ the endomorphism $\phi$ has $3$ distinct
real eigenvalues while in case of real moduli $<\mu_{0}$ the
endomorphism $\phi$ has one real and two purely imaginary eigenvalues.

\ifarx\vfill\eject\fi
\bigskip\medskip\Joker{KU}{Another type of examples} For $k=3$ and
$c\ge1$ let $V$, $\;\Gamma=\GL(2,\RR)^{0}$ and
$a=u^{m}_{2}+u_{1}u^{m-1}_{2}\in V$ be as in Example
\ruf{JP}. Consider the subgroup
$$
\Sigma:=\Big\{\Matrix\alpha0\beta1:\alpha\in\RRp,\beta\in\RR\Big\}
\;\subset\;\Gamma\,.\Leqno{LP}
$$
Then $F:=\Sigma(a)$ is a homogeneous surface in
$V$. For $o:=u_{2}^{m}$ the orbit $C:=\Sigma(o)$ is a homogeneous
curve. Identify $V$ and $\RR^{m+1}$ with coordinates
$(x_{0},x_{1},\dots,x_{m})$ as in Example \ruf{JP}. Then
$a=(1,1,0,\dots,0),\,o=(1,0,\dots,0)\in\RR^{m+1}$ and the Lie algebra
of $\Sigma$ corresponds to the linear span of the two vector fields
$\zeta_{1}$ and $\xi^{1,-1}$ in \Ruf{AQ}. In particular
$$
C=\{(1,t,t^{2},\dots,t^{m}):t\in\RR\}
$$
and $T_{o}C=\RR\cd b$ with $b:=(0,1,0,\dots,0)$ for the tangent space
at $o\in C$. On the other hand, the affine {half}line $o+\RRp\!\cd b$
is contained in $F$. The geometric meaning of this is the following:
The development $S:=\bigcup_{c\in C}(c+T_{c}C)$ of the curve $C$ is
divided by $C$ in two $\Sigma$-orbits, one of which is $F$ (compare
\Lit{EAEZ}, \p45 for the special case $m=3$). Now identify the
$\Sigma$-invariant hyperplane $W:=\{x\in V:x_{0}=1\}$ with $\RR^{m}$
by 'dropping the coordinate' $x_{0}$. Then $C$ becomes the twisted
$m$-ic $\{(t,t^{2},\dots,t^{m}):t\in\RR\}$, $o$ becomes the origin and
$a$ the first basis vector $(1,0,\dots,0)$ in $\RR^{m}$. In the
coordinates of $\RR^{m}$ the vector fields $\zeta_{1}$ and
$\xi^{1,-1}$ are affine and have the forms
$$
\zeta_{1}=\sum^{m}_{j=1}jz_{j}\dd{z_{j}}\Steil{and}
\xi^{1,-1}=\dd{z_{1}}+\sum^{m}_{j=2}jz_{j-1}\dd{z_{j}}\,.\Leqno{HQ}
$$
With Proposition \ruf{US} it is easily verified that
$K^{r}_{a}F=\{x\in\RR^{m}:x_{j}=0\steil{if}j+r>2\}$ for all
$r\ge0$. Since the twisted $m$-ic is not contained in any hyperplane
of $\RR^{m}$ we therefore get that the tube $M:=C+i\RR^{m}$ is a
homogeneous minimal $2$-nondegenerate submanifold of $\CC^{m}$ with
CR-dimension $2$ and CR-codimension $m-2$. Notice that for the cone
$\RRp\!\cd F$ generated by $F$ in $V$ the tube $\RRp\!\cd
F+iV$ is an open piece of $\5M^{3,n-3}$, $n:=m+1$, and hence is
$3$-nondegenerate.\nline Denote for every integer $j$ by $\7g^ {(j)}$
the $k$-eigenspace of $\ad(\zeta_{1})$ in $\7g:=\hol(M,a)$. Then every
$\xi\in\7g^{(j)}$ is a complex linear combination of monomial vector
fields $z_{1}^{\nu_{1}}z_{2}^{\nu_{2}}\cdots
z_{m}^{\nu_{m}}\dd{z_{p}}$ with
$\nu_{1}+2\nu_{2}+\dots+m\nu_{m}=j+p$. As in the proof of Proposition
\ruf{DC}, it is shown that $\7g$ has the $\ZZ$-gradation
$$
\7g=\bigoplus_{j\ge-m}\7g^{(j)}\,.\Leqno{RU}
$$ It can be seen that $\7g$ is the linear span of all $i\dd{z_{j}}$,
$1\le j\le m$, as well as $\zeta_{1}$ and $\xi^{1,-1}$. In particular,
$\7g$ is a solvable Lie algebra of dimension $m+2$, coincides with
$\aff(M,a)$ and has commutator subgroup of dimension $m+1$. A proof
will be sketched for the special case $m=3$ in Example \ruf{EV}.

\KAP{Acht}{Homogeneous 2-nondegenerate CR-manifolds in dimension 5}
\medskip

In this section we specialize the examples of the previous section to
the case $V=\RR^{3}$. We start with manifolds of type
$M=M^{\phi}=F^{\phi}+i\RR^{3}$ in $\CC^{3}$, Then the local
CR-equivalence classes of these manifolds are parameterized by the
modules $\mu(M)\in\overline{\RR}$, compare \Ruf{KW}. The $\phi$
occurring in the following examples not necessarily are trace-free but
easily could be transformed to.  As defined in the previous section
let $\mu_{0}:=27/4$.

\Example{EI}{\fett ($\mu=\infty$) } Let
$F:=\{x\in\RR^{3}:x_{1}^{2}+x_{2}^{2}=x^{2}_{3},\,x_{3}>0\}$ be the
future light cone. This surface occurs as $F^{\phi}$ for
$\phi:=x_{2}\dd{x_{1}}-x_{1}\dd{x_{2}}$ having spectrum $\{\pm i,0\}$.

\Example{EY}{\fett ($\mu<\mu_{0}$) } For $\omega>0$ let
$F\subset\RR^{3}$ be the orbit of $(1,0,1)$ under the group of all
linear transformations $x\mapsto r(\cos t\,x_{1}-\sin t\,x_{2},\sin
t\,x_{1}+\cos t\,x_{2},e^{\omega t}x_{3})$, $r\in\RRp,t\in\RR$.\nline
With $r:=(x_{1}^{2}+x_{2}^{2})^{1/2}$, the manifold $F$ is given in
$\{x\in\RR^{3}:r>0\}$ by the explicit equations
$$
x_{3}=r\exp\big(\omega\cos^{-1}(x_{1}/r)\big)=r\exp\big
(\omega\sin^{-1}(x_{2}/r)\big)\,,
$$
where locally always one of these suffices. A suitable choice is
$\phi=x_{1}\dd{x_{2}}-x_{2}\dd{x_{1}}+\omega x_{3}\dd{x_{3}}$ with
spectrum $\{\pm i,\omega\}$.

\medskip

\Example{EZ}{\fett ($\mu=\mu_{0}$) } Let $F\subset\RR^{3}$ be the
orbit of $(1,0,1)$ under the group of all linear transformations
$x\mapsto r(x_{1},x_{2}+tx_{1},e^{t}x_{3})$ with
$r\in\RRp,t\in\RR$, that is,
$$
F=\big\{x\in\RR^{3}:x_{1}>0,\;x_{3}^{}= x_{1}^{}\re1e^{ x_{2}/x_{1}}
\big\}\,.
$$
Here $\phi=x_{1}\dd{x_{2}}+x_{3}\dd{x_{3}}$ has characteristic roots
$0,0,1$.

\medskip

\Example{EX}{\fett ($\mu>\mu_{0})$ } For $\theta>2$ let
$$
F:=\{x\in(\RRp)^{3}:x_{3}=x_{1}(x_{2}/x_{1})^{\theta}\}\,.
$$
Here $\phi=x_{2}\dd{x_{2}}+\omega x_{3}\dd{x_{3}}$ has eigenvalues
$\{0,1,\omega\}$.

\bigskip

The following is Example \ruf{KU} specialized to $m=3$.

\Example{EV} Let $\Sigma$ be the group generated by the following two
one-parameter groups
$$
x\mapsto(e^{t}x_{1},e^{2t}x_{2},e^{3t}x_{3})\,,\qquad
x\mapsto(x_{1}+t,x_{2}+2tx_{1}+t^{2},x_{3}+3tx_{2}+3t^{2}x_{1}+
t^{3})\Leqno{NM}
$$
of affine transformations on $\RR^{3}$. Then $\Sigma$ is isomorphic to
the group defined in \Ruf{LP}. For $a:=(1,0,0)$ the orbit
$F:=\Gamma(a)$ is
$$
F\;=\;\{(t,t^{2},t^{3})+r(1,2t,3t^{2})\in\RR^{3}:r\in\RRp,t\in\RR\}\,.
$$
The tube $M:=F+i\RR^{3}$ is an affinely homogeneous
$2$-nondegenerate CR-manifold. The Lie algebra of $\Sigma$ is spanned
by the affine vector fields
$$
\zeta_{1}:=z_{1}\dd{z_{1}}+2z_{2}\dd{z_{2}}+3z_{3}\dd{z_{3}}
\Steil{and}\xi^{1,-1}:=\dd{z_{1}}+z_{1}\dd{z_{2}}+z_{2}\dd{z_{3}}
\,,
$$
compare also \Ruf{HQ}. The Lie algebra $\7g:=\hol(M,a)$ is of finite
dimension and has the gradation \Ruf{RU} for $m=3$, where $\7g^ {(k)}$
is the $k$-eigenspace of $\ad(\zeta_{1})$.  We claim that $\7g$ has
dimension $5$ and coincides with $\aff(M,a)$. The proof consists of
several elementary steps which we only sketch here. To begin with,
define $\7a^{(k)}\subset\7g^{(k)}$ by
$$
\7a^{(-3)}:=\RR\, i\dd{z_{3}}\,,\; \7a^{(-2)}:=\RR
i\dd{z_{2}}\,,\;\7a^{(-1)}:=\RR
i\dd{z_{1}}\oplus\RR\phi\,,\;\7a^{(0)}:=\RR\zeta
$$
and $\7a^{(k)}:=0$ for all other $k$. By induction on $k$ it is seen
that $\7g^{(k)}=\7a^{(k)}$ holds for all $k$: For $k\le-3$ this is
obvious. For $k=-2$ suppose there exists a
$\xi\in\7g^{(-2)}\setminus\7a^{(-2)}$. Then
$\xi=\alpha\dd{z_{2}}+\beta z_{1}\dd{z_{3}}$ for some
$\alpha,\beta\in\CC$. From $[\xi,\7a^{(-1)}]\subset\7a^{(-3)}$ we get
$\beta\in\RR$ and then $\xi_{a}\in T_{a}M$ implies $\beta=0$. But then
$\xi\in\7a^{(-2)}$ since $\7a^{(-2)}\!\oplus
i\7a^{(-2)}\not\subset\7g^{(-2)}$, a contradiction. For $k\ge-1$ the
procedure is as follows. Suppose there exists a
$\xi\in\7g^{(k)}\setminus\7a^{(k)}$. Then write $\xi$ as complex
linear combination of monomial vector fields as mentioned above and
subtract from $\xi$ a suitable element of $\7a^{(k)}$ thus killing as
many coefficients in front of monomial terms of the form
$f(z)\dd{z_{1}}$ as possible. By induction hypothesis
$[\7g^{(k)},\7a^{(j)}]\subset\7a^{(k+j)}$ holds for all $j<0$ and
gives $\xi=0$, a contradiction.  This proves the claim and also the
first part of the following statement.

\Lemma{LS} Let $M:=F+i\RR^{3}$ for $F=H(a)\subset\RR^{3}$ as in
Example \ruf{EV}. Then $\7g=\hol(M,a)$ is a solvable Lie algebra of
dimension $5$ with commutator algebra $[\7g,\7g]$ of dimension
$4$. Furthermore, $\Aut(M,a)=\{\id\}$.

\Proof Fix $a:=(1,0,0)\in F$ and write $\xi_j:=i\ddz{j}$ for
$j=1,2,3$. Let $h\in \Aut(M,a)$ and $\Theta$ the induced Lie algebra
automorphism of $\;\7l:=\7g\oplus i\7g$. With $\7g$ the subspaces
$$
\7n:=[\7g,\7g]=\bigoplus_{k<0}\7g^{(k)}\;,\qquad
[\7n,\7n]=\RR\xi_{2}\oplus
\RR\xi_{3}\Steil{and}[\7n,[\7n,\7n]]=\RR\xi_{3}
$$
are stable under $\Theta$ and hence also
$$
\7n\cap\;K_aM = \RR\xi_1\Steil{and}[\7n,\7n]\cap H_aM= \RR \xi_2\,,
$$
where $\7g$ and the tangent space $T_aM$ are identified via the
evaluation map. In particular, $\Theta(\7P_{-1})=\7P_{-1}$ and
$h(z)=g(z)+c$ for a diagonal matrix $g\in\GL(3,\RR)$ with $c=g(a)-a$,
compare \Ruf{GP}. Taking commutators of $h$ with all elements in the
second $1$-parameter group of \Ruf{NM} and then taking the derivative
by $t$ at $t=0$ gives $c\dd z\in\7g$. From $\7g\cap\;i\7g=0$ we
conclude $c=0$ and thus $g_{11}=1$ for the diagonal matrix $g$. Now
$\psi$ is the unique vector field in $\7g$ that has the value $\dd z$
at both points $0,a$, implying $\Theta(\psi)=\psi$. Therefore $g$ is
the unit matrix and $h$ is the identity in $\Aut(M,a)$.\qed

\smallskip For every $M=F+i\RR^{3}$ with $F\subset\RR^{3}$ a
cone from Example \ruf{EI} -- \ruf{EX} the commutator of $\hol(M,a)$
has either dimension $10$ (Example \ruf{EI}) or dimension $3$ (all the
others). As a consequence of Proposition \ruf{KO} and Lemma \ruf{LS}
we therefore get:

\Proposition{UA} The CR-manifolds $M=F+i\RR^{3}$ with
$F\subset\RR^{3}$ occurring in the examples \ruf{EI} -- \ruf{EV}, are
all homogeneous and $2$-nondegenerate. Furthermore, they are mutually
locally CR-inequivalent.\Formend

\medskip With an argument from \Lit{FEKA} together with 2.5.10 in
\Lit{HORM} a holomorphic extension property for {\sl global}
continuous CR-functions $f$ on $M=F\oplus i\RR^{3}\subset\CC^{3}$, $F$
one of the cones from examples \ruf{EI} -- \ruf{EV}, can be obtained:
Every such $f$ has a unique continuous extension to the convex hull
$\hat{M}$ of $M$ in $\CC^{3}$ that is holomorphic on the interior of
$\hat{M}$ with respect $\CC^{3}$. Since $M$ is completely contained in
the interior of $\hat{M}$ in case $F$ belongs to Example \ruf{EY},
every global continuous CR-function on such an $M$ is real-analytic.

\bigskip For the tube $\5M$ over the future light cone (that is
Example \ruf{EI}) there exist many (even simply-connected) homogeneous
CR-manifolds that are all locally CR-equivalent to $\5M$ but are
mutually non-diffeomorphic, compare \Lit{KAZT}. In contrast to this,
using already Theorem II from Section \ruf{Neun}, we can state the
following global result:

\Proposition{UB} Let $M$ be a homogeneous $2$-nondegenerate
CR-manifold that is not locally CR-equiva\-lent to the tube $\5M$ over
the future light cone. Then $M$ is simply connected and $\Aut(M)$ is a
solvable Lie group of dimension 5 acting transitively and freely on
$M$. For every $a\in M$ the stability group $\Aut(M,a)$ is trivial and
every homogeneous real-analytic CR-manifold $M'$, that is locally
CR-equivalent to $M$, is already globally CR-equivalent to $M$.

\Proof By Theorem II $M=F+iV$ for  $F$ as in one of the
examples \ruf{EY} -- \ruf{EV}. It is easily checked that $F$ and hence
$M$ is simply connected. In case $F$ is a cone, the claim follows with
Lemma \ruf{KR}. Therefore we may assume that $F$ is the
submanifold of Example \ruf{EV}. But then $\Aut(M,a)$ is the trivial
group by Lemma \ruf{LS} and $\Aut(M)$ has trivial center by
Proposition \ruf{ZI}. But then the claim follows from Proposition
\ruf{UQ}.~ \qed

\bigskip The affinely homogeneous surfaces $F\subset\RR^{3}$ of
examples \ruf{EI} -- \ruf{EV} occur already in \Lit{EAEZ} \p43. There
the surfaces are presented in their affine normal forms: Our Example
\ruf{EI} ($\mu=\infty$) corresponds to P4, the Examples \ruf{EY} --
\ruf{EX} ($\mu\in\RR$) to P3 and Example \ruf{EV} to P1. The remaining
degenerate types in \Lit{EAEZ}, types P2$^{\pm}$ and P5, do not show
up among our examples since the associated tube manifolds are
holomorphically degenerate.

Let us consider the type P3 in \Lit{EAEZ} \p43 a little bit closer.
This is the family of local surfaces in $\RR^{3}$ given by the local
equations in affine normal form
$$x_{3}=x_{1}^{2}+x_{1}^{2}x_{2}+x_{1}^{2}x_{2}^{2}+x_{1}^{5}+
x_{1}^{2}x_{2}^{3}+4x_{1}^{5}x_{2}+x_{1}^{2}x_{2}^{4}+ax_{1}^{7}
+10x_{1}^{5}x_{2}^{2}+x_{1}^{2}x_{2}^{5}+{\rm O}(8),\Leqno{OQ}$$ where
$a\in\RR$ is an arbitrary parameter and ${\rm O}(8)$ for every fixed
$a$ is a convergent power series in $x=(x_{1},x_{2},x_{3})$ vanishing
of order $\ge8$ at the origin and uniquely determined by the
requirement, that \Ruf{OQ} defines near the origin of $\RR^{3}$ a
locally affinely homogeneous surface. Different values of $a\in\RR$
give locally affinely inequivalent surfaces and there tubes in
$\CC^{3}$ correspond in a 1-1-way to our Examples \ruf{EY} --
\ruf{EX}. It is not difficult to see that the modulus $\mu$ of every
such surface is related to the parameter $a$ in \Ruf{OQ} by the
formula $100\mu=(28a)^{3}$.

\medskip In \Lit{DOKR}, \Lit{EAEZ} all locally affinely homogeneous
surfaces in $\RR^{3}$ have been classified up to local affine
equivalence. Inspecting the degenerate surfaces in these classifications
gives together with our results the following

\Proposition{IZ} Let $F$ be a locally affinely homogeneous surface in
$\RR^{3}$ and assume that the corresponding tube $M:=F+i\RR^{3}$
in $\CC^{3}$ is $2$-nondegenerate. Then \0 $M$ is locally
CR-equivalent to a manifold occurring in the Examples \ruf{EI} --
\ruf{EV}, and \1 for every further locally affinely homogeneous
surface $F'$ in $\RR^{3}$ the corresponding tubes $M$ and
$M':=F'+i\RR^{3}$ are locally CR-equivalent if and only if $F$
and $F'$ are locally affinely equivalent.\Formend

\Partskip

\ifarx\bigskip\fi

\centerline{\twelverm PART 2: The classification}

\KAP{Neun}{Lie-theoretic characterization of locally homogeneous
CR-manifolds}

In this part of the paper we  classify
all homogeneous 5-dimensional 2-nondegenerate CR-manifolds up to
local CR-equivalence, that is, we
carry out the proof of

\medskip\noindent
{\bf Theorem II.} {\sl
   Let $M$ be a locally homogeneous $2$-nondegenerate
real-analytic CR-manifold of dimension $5$. Then $M$ is locally
CR-equivalent to a tube $F+ i\RR^{3}\;\subset\;\CC^{3}$, where
$F\subset\RR^{3}$ is one of the affinely homogeneous surfaces
occurring in the Examples \ruf{EI} -- \ruf{EV}.}

\smallskip\noindent We call (in accordance with Section \ruf{Zwei}) a
real-analytic CR-manifold $M$ {\sl locally homogeneous} at a point
$o\in M$ if there exists a Lie subalgebra $\7g\subset\hol(M,o)$ of
finite dimension such that the canonical evaluation map $\7g\to
T_{o}M$ is surjective, that is, such that the tangent vectors
$\xi_{o}$, $\xi\in\7g$, span the tangent space $T_{o}M$. If this is
the case we also call the corresponding CR-germ $(M,o)$ {\sl locally
homogeneous}. If a particular locally transitive $\7g\subset\hol(M,o)$
has been fixed we also say that the germ $(M,o)$ is $\7g$-homogeneous.

The proof of Theorem II relies on a natural equivalence (see
\Lit{FELS}, Prop. 4.1) between the category of CR-manifold germs with
a locally transitive Lie algebra action and a certain purely
algebraically defined category. Before we briefly outline the main
steps of our proof, we recall the notion of a CR-algebra, taken from
\Lit{MENA}, and introduce some notation.

\Definition{FJ} A CR-algebra is a pair $(\7g,\7q)$, where $\7g$ is a
real Lie algebra of {\sl finite} dimension and $\7q$ is a complex Lie
subalgebra of the complexification $\7l:=\7g\oplus i\7g$. The
CR-algebra $(\7g,\7q)$ is called {\sl effective} if $0$ is the only
ideal of $\7g$ contained in $\7g\cap\7q$. \Formend

\Remarks{XA} \0 In \Lit{MENA} also the case is allowed where $\7g$ has
infinite dimension, but $\7q$ has to have finite codimension in $\7l$.
In this part of the paper however, only finite-dimensional Lie
algebras occur. \1 The CR-algebras form in an obvious way a category:
A morphism $(\7g,\7q)\to(\7g',\7q')$ of CR-algebras is a Lie algebra
homomorphism $\7g\to\7g'$ in the usual sense whose complex linear
extension $\7l\to\7l'$ maps $\7q$ to $\7q'$. Unfortunately, the
resulting notion of isomorphism between CR-algebras is too strict for
our purposes. We therefore mainly work with the coarser notion of
geometric equivalence between CR-algebras to be introduced later. \1
The geometric situation behind the notion of a CR-algebra is the
following: Let $Z$ be a complex manifold homogeneous under a complex
Lie group $L$, let $o\in Z$ be a point with isotropy subgroup
$Q\subset L$ at $o$ and $G\subset L$ a connected real form of $L$,
that is, the connected identity component of the fixed point set
$L^\sigma$ for an involutive antiholomorphic automorphism $\sigma$ of
$L$. Then each $G$--orbit $M$ in $Z$ is a generic (immersed)
CR-submanifold. Let $\7g,\7q$ be the Lie algebras of $G$ and $Q$,
respectively. Then $(\7g,\7q)$ is a CR-algebra that completely
describes the CR-germ $(M,o)$ together with the local action of $G$
near $o$. \1 In general, however, not every CR-algebra $(\7g,\7q)$ can
be obtained from a global situation as described in (iii), only from a
more general local setting.  Nevertheless, the set \CRG of all
CR-equivalence classes of locally homogeneous CR-germs and the set
\CRA of all geometric equivalence classes of CR-algebras stand in a
canonical 1-1-correspondence, compare also Section 4 in
\Lit{FELS}. For convenience of the reader, we briefly describe below
this correspondence in both directions. As a reference for a general
discussion of `local' and `infinitesimal' actions we refer to the
original paper of Palais, \Lit{PALA}.\Formend

\medskip\noindent In the following let $M$ always be a locally
homogeneous real-analytic CR-manifold with base point $o\in M$ that
{\sl locally} can be embedded in some $\CC^{n}$ (equivalently, $M$
is an (abstract) real-analytic {\sl involutive} CR-manifold as defined
at the end of Section \ruf{Zwei}). Each such CR-manifold
can globally and generically be embedded into a complex manifold $Z$
\Lit{ANFR}. Since we are only interested in the local structure of $M$
at $o$ and therefore mostly deal with CR-germs $(M,o)$, we may assume
without loss of generality that $M$ is embedded in a complex vector
space $E\cong \CC^n$ as a locally closed generic CR-submanifold.

Next we describe the interplay between locally homogeneous CR-germs
and CR-algebras more closely. In particular, we give two canonical
constructions that induce the 1-1-correspondence between the sets \CRG
and \CRA mentioned in \ruf{XA}.iv and also allow the precise
definition of `geometric equivalence' for CR-algebras:

Let $(M,o)$ be a CR-germ and $\7g\subset\hol(M,o)$ a locally
transitive Lie subalgebra of finite dimension. Then an effective
CR-algebra $(\7g,\7q)$ can be associated in the following way: To
begin with, realize $\hol(M,o)$ in the canonical way as real Lie
subalgebra of the complex Lie algebra $\hol(E,o)$. This is possible,
since we assumed $M$ to be generic in $E$ and we can use Proposition
12.4.22 in \Lit{BERO}. As in Definition \ruf{FJ} we always denote by
$\7l=\7g^{\CC}=\7g\oplus i\7g$ the formal complexification of
$\7g$. Now $\Xi(\xi+ i\eta):=\xi+J\eta$ defines a Lie algebra
homomorphism $\Xi:\7l\to\hol(E,o)$, where $J$ denotes the complex
structure tensor $J:TE\to TE$. We will not make a notational
distinction between the complex structures in $\7l$ or $TZ$, and write
`$i$' for it. The homomorphism $\Xi$ is in general not injective (it
is, if $M$ is holomorphically nondegenerate). Let $\7q\subset\7l$ be
the $\Xi$-preimage of the isotropy subalgebra
$\{\xi\in\hol(E,o):\xi_{o}=0\}$. Then the CR-algebra $(\7g,\7q)$ is
called a {\sl CR-algebra associated} to the locally homogeneous
CR-germ $(M,o)$. It is obvious that $\7g\cap\7q$ is nothing but the
isotropy subalgebra $\7g_{o}=\{\xi\in\7g:\xi_{o}=0\}$ and the tangent
space $T_{o}M$ can be canonically identified with $\7g/\7g_{o}$. Also
the holomorphic tangent space $H_{o}M$ and the partial complex
structure $J:H_oM\to H_oM$ (equivalently: the decomposition
$H^{1,0}_oM\oplus H^{0,1}_oM$ of the complexification
$H^{\CC}_oM=H_oM\otimes \CC$) can be read off the CR-algebra
$(\7g,\7q)$: Let $\sigma$ be the conjugate linear involution of $\7l$
with $\7l^{\sigma}:=\Fix(\sigma)=\7g$. Then it is easily verified that
$\7H:=(\7q+\sigma\7q)^{\sigma}$ coincides with $\{\xi\in\7g:\xi_{o}\in
iT_{o}M\}$, that is, $\7H/\7g_{o}$ is canonically isomorphic to
$H_{o}M$ (the capital letter for $\7H$ is chosen to indicate that
$\7H$ in general is not a Lie algebra, only a linear
subspace). Further, $H^{0,1}_oM=\qu{\7q}{\7q\!\cap\!\sigma\7q}$ and
$H^{1,0}_oM=\qu{\sigma\7q}{\7q\!\cap\! \sigma\7q}$. In \Lit{FELS} it
has been shown that the geometric properties of the CR-structure of
the CR-germ $(M,o)$ like minimality, $k$-nondegeneracy, holomorphic
degeneracy can be read off every CR-algebra $(\7g,\7q)$ associated
with $(M,o)$. The facts relevant for our classification will be
discussed below.

There is also a canonical way to associate a locally homogeneous
CR-manifold germ $(M,o)$ to a given CR-algebra $(\7g,\7q)$ (not
necessarily effective): Choose a complex Lie group $L$ with Lie
algebra $\7l=\7g\oplus i\7g$ and a complex linear subspace
$E\subset\7l$ with $\7l=\7q\oplus E$. Then there exist open
neighbourhoods $U$ of $0\in E$, $V$ of $0\in\7q$ and $R$ of $\id\in L$
such that $(u,v)\mapsto\exp(u)\exp(v)$ defines a biholomorphic mapping
$\phi:U\times V\to R$. Choose an open neighbourhood $P$ of $0\in\7g$
with $\exp(P)\subset R$. Then, in our particular situation, with
$\pi:U\times V\to U$ being the canonical projection, the mapping
$\psi:=\pi\circ\phi^{-1}\circ\exp:P\to U$ has constant rank. Without
loss of generality we therefore may assume that $M:=\psi(P)$ is a
connected real-analytic submanifold of $E$ containing the origin $0\in
E$. The Lie algebra $\7l$ can be identified with the Lie algebra of
all right-invariant vector fields on $L$ and every $\xi\in\7l$ can be
projected along $\pi\circ\phi^{-1}:R\to U$ to a holomorphic vector
field $\tilde\xi\in\hol(U)$. Thus the real subalgebra
$\tilde\7g:=\{\tilde\xi:\xi\in\7g\}\subset\hol(U)$ is a homomorphic
image of $\7g$ and spans at every $x\in M$ the tangent space
$T_{x}M$. In particular, $M$ is a generic CR-submanifold of $E$ and
$\tilde\7g$ is a locally transitive subalgebra. It is not difficult so
see that $\tilde\7g$ is obtained from $\7g$ by factoring out the
kernel of ineffectivity, more precisely, let $\7j$ be the largest
ideal in $\7g$ with $\7j\subset \7g\cap\7q $. Then $\tilde\7g$ is
isomorphic to $\7g/\7j$. If there exists a Lie group $L$ with Lie
algebra $\7l$ such that the subalgebra $\7q$ corresponds to a {\sl
closed} complex subgroup $Q\subset L$ then we may take $L/Q$ for $U$
and the $G$--orbit through $[Q]\in L/Q$ for $M.$ We call $(M,o)$ the
CR-{\sl germ associated} to the CR-algebra $(\7g,\7q)$.

\Definition{} The CR-algebras $(\7g,\7q)$ and $(\7g',\7q')$ are called
{\sl geometrically equivalent} if the associated CR-germs are
CR-equivalent.

Notice that CR-algebras are always geometrically equivalent if they
are isomorphic in the categorical sense of \ruf{XA}.ii, but not
conversely in general. Notice also that every CR-algebra is
geometrically equivalent to an effective one. In the following
paragraphs \ruf{FA}--\ruf{FD}, we fix

\vskip5pt\centerline{\sl for the rest of the paper}

\noindent the basic setup and notation, which are mainly taken from
\Lit{FELS}:

\Notation{FA} Given a CR-algebra $(\7g,\7q)$, let $(M,o)$ be the
associated CR-germ. Write $\7l:=\7g^{\CC}=\7g\oplus i\7g$ for the
complexification and $\sigma:\7l\to \7l$ for the complex conjugation
with $\7l^\sigma=\7g$. Then \0 $\7g_o : =\7g\cap \7q$ is called the
{\sl real isotropy subalgebra}. Define
$\7l_o:=\7q^{(\infty)}:=\7q\cap\sigma\7q$ and note that $\7l_o$ is
the complexification $(\7g_o)^{\CC}$ of $\7g_{o}$.  \1 $\7g/\7g_{o}$
and $\7H/\7g_{o}\subset\7g/\7g_{o}$ for
$\7H:=(\7q+\sigma\7q)^\sigma\subset\7g$ are called the {\sl real} and
the {\sl holomorphic tangent space} respectively.  \1 The descending
chain $\7q^{(0)}\supset\7q^{(1)}\supset\7q^{(2)}\supset\cdots\supset
\7q^{(\infty)}$ of complex subalgebras is inductively defined by
$\7q^{(0)}:=\7q$, $\7q^{(\infty)}:=\7q\cap \sigma\7q$ and
$\7q^{(k+1)}:=\{w\in\7q^{(k)}:[w,\sigma\7q]\subset\7q^{(k)}\li1+
\sigma\7q\}$, for $k\in \NN$.

\noindent If $(M,o)$ is the manifold germ associated to the CR-algebra
in \ruf{FA} then the holomorphic tangent space
$\7H/\7g_{o}\subset\7g/\7g_{o}$ in the sense of (ii) can be
canonically identified with the holomorphic tangent space $H_{o}M$ in
the geometric sense. As shown in \Lit{FELS}, the mapping $\7q\to\7H$,
$w\mapsto w+\sigma w$, induces a complex linear isomorphism
$\7q/\7q^{(\infty)}\cong\7H/\7g_{o}\,$. The former quotient is
canonically isomorphic to $H^{0,1}_oM$ (similarly, $H^{1,0}_oM\cong
\sigma\7q/\7q^{(\infty)}$). The Levi kernel $K_{o}M$ and its higher
order analogues $K^{r}_{o}M$ can be considered as complex linear
subspaces of $\7q/\7q^{(\infty)}$. We will make extensive use of some
of the main results of \Lit{FELS}, such as Theorem 5.10:

\Joker{CA}{Algebraic characterization of $k$-nondegeneracy} For every
$r\ge0$ the space $\7q^{(r)}$ is a Lie subalgebra of $\7q$ and the
$r^{{\th}}$ Levi kernel $K^{r}_{o}M$ is isomorphic to
$\7q^{(r)}/\7q^{(\infty)}$. In particular, for every $k\ge1$ the
locally homogeneous CR-manifold $M$ is $k$-nondegenerate if and only
if

\vskip5pt\centerline{$\7q^{(k-1)}\ne\;\7q^{(k)}=\;\7q^{(\infty)}$.}\Formend

\medskip\noindent To handle the $5$-dimensional case we also introduce
the following abbreviations \1 $\7f:=\7q^{(1)}=\{v\in
\7q:[v,\sigma\7q]\subset \7q+\sigma\7q\}$ and
$\7F:=(\7f+\sigma\7f)^\sigma$.

\noindent Then $\7g_o\subset\7F\subset \7H\subset\7g$ are real
subspaces stable under $\ad(\7g_o)$ and $\7l_o\subset\7f\subset\7q$
are complex subalgebras. In Lemma 5.9. of \Lit{FELS} it has been shown
that actually $\7F$ is a Lie algebra and coincides with
$N_{\7g}(\7H)\cap \7H$ (here, given $W\subset \7g$, $N_\7g(W):=\{v\in
\7g:[v,W]\subset W\}$).

\smallskip Summarizing the above discussion, the following result is
the key for our classification.

\Proposition{MC} Let $(M,o)$ be an $\7g$-homogeneous CR-germ and let
$(\7g,\7q)$ be the corresponding CR-algebra.  Then $M$ is
5-dimensional, minimal and 2-nondegenerate if and only if
$$
\codim_{\7g}(\7g_{0})=5\Steil{and}
\7q\ne\7f\ne\7q^{(2)}=\7q^{(\infty)}:=\7q\cap\sigma\7q\;.
$$

\ifarx\eject\fi

{\openup2pt\Lemma{FC} The Lie algebraic terms in \Ruf{MC} are
equivalent to the following set of conditions: \o
$\dim\qu{\7F}{\7g_o}=\dim \qu{\7H}{\7F}=2,\quad
\dim_{\RR}\qu{\7g}{\7H}=1=
\dim_{\CC}\qu{\7f}{\7l_o}=\dim_{\CC}\qu{\7q}{\7f}$

\l $[\7g_o,\7F]\subset\7F,\quad[\7F,\7F]\subset\7F,\quad
[\7F,\7H]\subset\7H$

\l $[\7q,\sigma\7q]\not\subset\7q\,{+}\kern1pt\sigma\7q$\hfill\hbox
to 7truecm{\Klein\sl $M$ is not Levi flat\hss}

\l $[\7f,\sigma\7q]\,\subset \7q\,{+}\kern1pt\sigma\7q,\quad
\7f\ne\7l_o$\hfill\hbox to 7truecm{\Klein\sl $M$ is Levi
degenerate\hss}

\l $[\7f,\sigma\7q]\,\not\subset\,\7f\,{+}\,\sigma\7q\,$\hfill\hbox to
7truecm{\Klein\sl $M$ is 2-nondegenerate.\hss}
}
\Formend

\smallskip\noindent We will frequently use the fact that the condition
`$\,[\7F,\7H]\subset \7F\,$' (instead of $\subset \7H$) violates
condition \25.

\medskip By the canonical bijection between the classes \CRG and \CRA
mentioned in \ruf{XA}.iv and made precise above our classification
problem is transferred to the classification of certain effective
CR-algebras up to geometric equivalence. Unfortunately, to a given
locally homogeneous CR-germ $(M,o)$ there may be associated many
CR-algebras $(\7g,\7q)$ for which the $\7g$'s are non-isomorphic. For
instance, if $M=\5M$ is the tube over the future light cone, then to
$(\5M,o)$ there are associated CR-algebras $(\7g,\7q)$ with
$\7g\cong\so(2,3)$, $\so(1,3)$, $\so(2,2)$ together with a bunch of
other Lie algebras that are not semisimple, see \Lit{FEKA} for
explicit realizations. Therefore, the best we can do in the following
is to consider only CR-algebras $(\7g,\7q)$ such that $\dim\7g$ is
minimal in the geometric equivalence class of $(\7g,\7q)$.  But, also
then, we are still left in the example $(\5M,o)$ with several
non-isomorphic solvable Lie algebras of dimension 5 as well as one
non-solvable Lie algebra of dimension $5$ with $2$-dimensional
non-abelian radical and Levi part $\cong\7{sl}(2,\RR)$.

\smallskip \Joker{FD}{Fundamental Assumption} In the following, every
CR-algebra $(\7g,\7q)$ under consideration is assumed to satisfy the
condition \Ruf{MC} (equivalently \21 -- \25\hskip1pt) as well as the
following additional condition. \l For every CR-algebra $(\7g',\7q')$,
which is geometrically CR-equivalent to $(\7g,\7q)$, the dimension
estimate $\dim\7g'\ge\dim\7g$ holds.

\smallskip\noindent Condition \26 implies, in particular, the
following two conditions.

\item{\26$_{1}$} $(\7g,\7q)$ is effective.

\item{\26$_{2}$} There is no proper subalgebra $\7g'\subset \7g$ with
$\7g'+\7g_o=\7g$ for $\7g_o=\7g\cap\7q$.

\noindent Indeed, $(\7g',\7q')$ with $\7q':=(\7g'+i\7g')\cap\7q$ is a
CR-algebra that is geometrically equivalent to $(\7g,\7q)$ in case
$\7g'+\7g_o=\7g$.

\bigskip\Joker{ST}{Basic structure theory} In the following we
frequently use standard facts concerning reductive Lie algebras and
parabolic subalgebras, see \Lit{KNAP}, Chap.VI-VII, as a general
reference. To fix our notation, let $\7s$ be a complex reductive Lie
algebra and $\7t$ a Cartan subalgebra of $\7s$. Denote by
$\Phi=\Phi(\7s,\7t)\subset \7t^*$ the corresponding root system and by
$\Pi\subset\Phi$ the subset of simple roots. A subalgebra
$\7r\subset\7s$ is called {\sl parabolic} if it contains a maximal
solvable subalgebra (also called a {\sl Borel subalgebra}) $\7b$ of
$\7s$. Conjugacy classes of parabolic subalgebras in $\7s$ are
parameterized by the subsets of $\Pi$: For every $\6P\subset \Pi$ and
$\langle\kern-2pt\langle \6P \rangle\kern-2pt\rangle:=\Phi\cap
\bigoplus_{\alpha\in \6P} \ZZ \alpha$ the corresponding parabolic
subalgebra is defined by
$$
\7r=\7r_\6P:=\7r^{\red}\oplus
\7r^{\nil}\Steil{with}\7r^{\red}:=\7t\;\oplus\!
{\textstyle\bigoplus\limits_{\alpha\in\langle\!\langle
\6P\rangle\!\rangle}}\7s_\alpha\Steil{and}\7r^{\nil}:=\!
{\textstyle\bigoplus\limits_{\alpha\notin\langle\!\langle
\6P\rangle\!\rangle}}\7s_\alpha \;\;.\Leqno{PA}
$$
The case of a reductive {\sl real} Lie algebra $\7s$ is a little bit
more sophisticated. In contrast to the complex situation there may
exist several conjugacy classes of Cartan subalgebras. Among these the
class most suitable for our purposes consists of the so-called {\sl
maximally split Cartan subalgebras} $\7t\subset\7s$, defined as
follows: Select a Cartan decomposition $\7s=\7k\oplus\7p$ and let
$\7a$ be a maximal abelian subalgebra of $\7p$. Consider the
centralizer $\7m:=C_\7k(\7a)$ and put $\7t:=\7a\oplus \7t_\7m$, where
$\7t_\7m\subset \7m$ is a maximal abelian subalgebra. The conjugacy
classes of real parabolic subalgebras in $\7s$ are parameterized by
the subsets of the simple roots $\Pi\subset \Phi(\7s,\7a)$ in the
restricted root system $\Phi(\7s,\7a)\subset \7a^*$. Each parabolic
subalgebra $\7r$ in $\7s$ has the decomposition
$\7r=\7r^{{\red}}\ltimes \7r^{\nil}$ into the reductive and nilpotent
parts (once a maximally split Cartan subalgebra $\7t\subset\7r$ is
selected, the reductive factor with $\7r^{\red}\supset \7t$ is
unique). Following a common convention the roots $\lambda$ occurring
in the root space decomposition of the nilpotent ideal
$\7r^{\nil}=\bigoplus_{\lambda\in\Lambda} \7s_\lambda$ are negative
and we write $\7r^{-\n}:=\7r^{\nil}$, $\;\7r^{\r}:=\7r^{\red}$ and
$\Phi^{-\n}:=\Phi(\7r^{\nil},\7a)=\Lambda$. Each parabolic subalgebra
$\7r$ containing $\7t$ determines the decomposition
$$
\7s=\7r^{\n}\oplus \7r^{\r}\oplus \7r^{-\n}\Steil{with}
\7r^{\n}:=\!\!{\textstyle\bigoplus\limits_{\lambda\in\Phi(\7r^
{\nil},\7a)}}\li{12}\7s_{-\lambda}\;,\Steil{and we
put}\7r^{\opp}:=\7r^{\r}\ltimes \7r^{\n}\;.\Leqno{FG}
$$
\vskip-6pt\noindent We call $\rk(\7s):=\dim \7t$ the {\sl rank} of
$\7s$ and $\rk_{\RR}(\7s):=\dim \7a$ with $\7a\subset \7p$ as above
the {\sl real rank} of $\7s.$

\KAP{Zehn}{The Lie algebra $\7g$ has small semisimple part}

In the preceding section we have explained how 2-nondegeneracy can be
expressed in pure Lie algebraic terms. In this section we start with
the actual proof of Theorem II. Since this proof will be quite
involved, we subdivide it into several sections, lemmata and claims
until the final step is completed in Section \ruf{Sechzehn}. For the
convenience of the reader we briefly outline the main steps.

As explained above the classification can be reduced to the
determination of all CR-algebras satisfying the fundamental assumption
\ruf{FD}.  Once all possible CR-algebras are known, we have to
identify the underlying CR-germs. In general, it may happen that
algebraically non-equivalent CR-algebras give rise to equivalent
CR-germs. For this last part of the proof we use results from Section
\ruf{Acht}.

Our proceeding will be to show that the assumption \ruf{FD} severely
restricts the possibilities for $(\7g,\7q)$.  This will be achieved by
a detailed structural study of the Lie algebras $\7g$ occurring in
$(\7g,\7q)$. Recall that every Lie algebra $\7h$ has a Levi-\Malcev
decomposition $\7h=\7h^{\s}\ltimes\rad(\7h)$, where $\7h^{\s}$ is
semisimple and is uniquely determined up to an inner automorphism of
$\7h$. Furthermore, $\rad(\7h)$ is the radical of $\7h$, i.e., the
unique maximal solvable ideal in $\7h$. In the first part of the proof
we investigate the various possibilities for $\7g^{\s}$ where
$(\7g,\7q)$ satisfies certain conditions stated in the previous
section. To be precise: {\sl For the rest of the paper the fundamental
assumption \ruf{FD} remains in force for {\rm all} CR-algebras
$(\7g,\7q)$ under consideration}. In particular, every $(\7g,\7q)$ is
effective (condition \26$_{1}$, and therefore $\7g$ can be considered
as a transitive Lie subalgebra of $\hol(M,o)$). Furthermore, condition
\26$_{2}$ states that there is no proper Lie subalgebra of $\7g$ that
also is transitive on $(M,o)$.

Let an arbitrary CR-algebra $(\7g,\7q)$ subject to \ruf{FD} be given
and let $\7g^{\s}\ltimes \rad(\7g)$ be a Levi-\Malcev decomposition
of $\7g$.  In this and in the following few sections we assume that
$\7g^{\s}\ne 0$ and investigate which simple factors can occur in
$\7g^{\s}$. Thereby we use the following notation: We fix a simple
ideal $\7s$ in $\7g^{\s}$ and denote by $\7s'$ the corresponding
complementary ideal, i.e.,
$$
\7g=\7g^{\s}\ltimes \rad{\7g}\Steil{and} \7g^{\s}=\7s\times \7s'\;.
\Leqno{SS}
$$

In this section we show that $\7g^{\s}$ only can contain simple
factors isomorphic to one of the Lie algebras $\so(2,3),\;
\so(1,3),\;\su(2)$ and $\7{sl}(2,\RR)$.  This result is obtained by
analyzing which simple real Lie algebras can contain proper
subalgebras of very low codimensions.  In the next sections we exclude
further possibilities for $\7s$: In Section \ruf{Elf} we show that the
factors $\so(2,3)$ and $\so(1,3)$ cannot occur in $\7g^{\s}$ as it
will turn out that their existence in the Levi factor would violate
the minimality assumption \26. Nevertheless notice that there exist
CR-algebras $(\so(2,3),\7q)$ and $(\so(1,3),\7q)$ satisfying \21 --
\25 and \26$_1$. All the underlying CR-germs of such CR-algebras are
locally CR-equivalent to the light cone tube $\5M$. In Section
\ruf{Zwoelf} we find first examples of CR-algebras which satisfy
\ruf{FD}. In these cases $\7g^{\s}$ contains a simple factor $\7s\cong
\7{sl}(2,\RR)$, and then necessarily $\7g\cong \7{sl}(2,\RR)\times
\7r$, where $\7r$ is a 2-dimensional non-abelian Lie algebra. Also all
such CR-algebras give rise only to CR-germ locally CR-equivalent to
the tube over the light cone.  In Section \ruf{Dreizehn} we finally
eliminate the possibility $\7s\cong \su(2)$ for the simple factor
$\7s$ in $\7g^{\s}$.  At that stage of the proof we will have proved
the following dichotomy: Let $(\7g,\7q)$ be an arbitrary CR-algebra
which obeys \ruf{FD}. If $\7g^{\s}\ne 0$ then $\7g^{\s}\cong
\7{sl}(2,\RR)$ and furthermore the underlying CR-germ is locally
CR-equivalent to the light cone tube $\5M$.  Or, $\7g^{\s}=0$, i.e.,
$\7g$ is solvable.

For that reason, from Section \ruf{Vierzehn} on we only consider
CR-algebras $(\7g,\7q)$ with $\7g$ solvable and show that then
necessarily $\dim \7g=5=\dim M$. In Section \ruf{Fuenfzehn} we look
closer at the nilcenter $\7z$ of $\7g$ and find out that $\dim
\7z\in \{1,3\}$. The Main Lemma \ruf{GS}, which might be of interest
for itself, gives a sufficient condition for a CR-algebra to be
associated to a tube $F+i\RR^{3}\subset\CC^{3}$ over an affinely
homogeneous surface $F\subset\RR^{3}$. In the last section we show by
ad-hoc methods that indeed every 5-dimensional solvable Lie algebra
$\7g$ occurring in the CR-algebra $(\7g,\7q)$ under consideration
fulfills the assumption of the Main Lemma.  Hence, due to the
aforementioned assertions and since (up to local affine equivalence)
all affinely homogeneous surfaces in $\RR^{3}$ occur among the
examples \ruf{EI} -- \ruf{EV}, this completes the proof of the
classification theorem.

\medskip

We now begin with the proof of Theorem II:

\Lemma{LA} Let $(\7g,\7q)$ be a CR-algebra subject to \Ruf{FD}.
Then the simple Lie subalgebra $\7s\subset \7g^{\s}$ can only
be isomorphic to
$\;\so(2,3)\;,\;\so(1,3)\,,\,\7{sl}(2,\RR)\;$ or $\;\su(2)$.

\Proof The proof is carried out in several reduction steps. To begin
with, we write as shorthand
$$
\7h_o:=\7g_o\cap \7h\,,\quad\7h_\7F:=\7F\cap
\7h\steil{and}\7h_\7H:=\7H\cap\7h
$$
for every subalgebra $\7h\subset\7g$. Notice that $\7h_o\subset
\7h_{\7F}$ are subalgebras and $\7h_{\7H}$ is a linear
$\ad(\7h_{\7F})$--stable subspace of $\7h.$ \nline Consider the
subalgebras $\7s_o\subset \7s_\7F$ of $\7s$. The case
$\7s_o=\7s_\7F=\7s$, that is $\7s\subset\7g_o$, can be ruled out since
then $\7s'\oplus\rad(\7g)$ would be a proper locally transitive
subalgebra of $\7g$, contradicting assumption \26$_{2}$. Therefore, at
least one of the inclusions $\7s_{o}\subset \7s_\7F\subset\7s$ is
proper. Consequently, there is always a proper subalgebra of
codimension $\le3$ in $\7s$. Indeed, in case $\7s_\7F\ne\7s$ the
subalgebra $\7s_\7F$ has this property and in case $\7s_\7F=\7s$ the
proper subalgebra $\7s_{o}$ has codimension $\le2$. Hence, there
exists a maximal proper subalgebra $\7h$ of $\7s$ with either
$\7s_{\7F}\subset\7h$ or $\7s_o\subset \7h.$ Such a maximal subalgebra
$\7h$ has codimension $\le3$ in $\7s$. Due to \Lit{BOUR}, Chap.~VIII,
\S 10, Cor.~1, every maximal proper subalgebra of $\7s$ is either
reductive or parabolic. In the following claims we list all simple Lie
algebras $\7s$ which admit proper maximal subalgebras of such low
codimensions. We discuss the reductive and parabolic case separately.

\noindent{\fett Claim 1: } \sl Let $\7k$ be a simple real algebra and
$\7h\subset \7k$ reductive with
$0<\codim_\7s\7h\le 3$. Then $\7k$ can only be isomorphic to
$\7{sl}(2,\RR),\;\su(2)$ or $\so(1,3)$.

\ProofC Let $\7h=\7h_1\times \cdots \times \7h_p\times \7z$ be the
decomposition of the reductive subalgebra $\7h$ into the simple
factors $\7h_j$ and the center $\7z.$ Weyl's Theorem implies the
existence of an $\ad(\7h)$-stable complement $\7v\subset \7s.$ Let
$\rho:\7h\to \gl(\7v)$ be the induced adjoint representation. Every
restriction $\rho_j:\7h_j\to \gl(\7v)$ must be faithful since
otherwise $\7h_j$ would be an ideal in $\7s$. The crucial condition
here is $\dim \7v\le 3$ which, in turn, implies that each $\7h_j$ is
isomorphic either to $\7{sl}(2,\RR)$, $\so(3)$ or $\7{sl}(3,\RR)$. As
a consequence
$$
\dim \7h\le \cases{8k&$r=2k$\cr8k+3&$r=2k+1$\cr}
\Steil{for}r:={\rk}(\7h)\le{\rk}(\7s)\,.
$$
On the other hand, a
glance at the classification of simple Lie algebras shows
$$
\dim \7s \ge \cases{2k^2+4k=\dim_{\RR} \7{sl}(k+1,\CC)& $r=2k$\cr
4k^2+8k+3=\dim\7{sl}(2k+2,\RR) &$r=2k+1\,.$}
$$
Putting both inequalities together and bearing in mind $\dim\7s-\dim
\7h\le 3$ shows that the rank of $\7s$ can only be one of the numbers
$1,2,4$. Since $\7{sl}(2,\RR)$ and $\su(2)$ are the only simple real
Lie algebras of rank $1$ we may assume that $\7s$ has rank $2$ or $4$.
The case $\rk(\7s)=4$ can be ruled out in the following way: Consider
first the situation where $\7s$ is {\sl of complex type,} that is,
$\7s$ is the underlying real Lie algebra of a complex simple Lie
algebra $\7c$ of (complex) rank $2$. Then $\7c$ is either
$\7{sl}(3,\CC)$ or $\so(5,\CC)$. But in both cases a proper reductive
real subalgebra has at least (real) codimension 4 (in fact, 8). If
$\7s$ is {\sl of real type} then the above estimates give
$\dim\7h\le16$ and $\dim\7s\ge24=\dim \7{sl}(5,\RR)$. For the remaining
case $\rk(\7s)=2$ either $\7s\cong \so(1,3)$, which is in the list of
the claim, or $\7s$ is isomorphic to a real form of $\7{sl}(3,\CC)$ or
$\so(5,\CC)$. In both cases every proper reductive complex subalgebra
has at least codimension 4. This proves Claim 1.\qd

\noindent{\fett Claim 2:} \sl Let $\7k$ is a simple real Lie algebra
and $\,\7h\subset \7k$ parabolic with $0<\codim_\7s\7h\le 3$. Then
$\7s$ is isomorphic to $\so(1,3)$, $\7{sl}(4,\RR)$, $\su(2)$ or to a
non-compact real form of $\7{sl}(3,\CC) $ or $\so(5,\CC)$.

\ProofC Since every parabolic subalgebra of a compact Lie algebra is
trivial, apart from $\7s\cong \su(2)$ we only have to consider the
case where $\7s$ is non-compact. The estimate
$\rk(\7s)\le\dim\7s-\dim\7h$, compare \Ruf{ST}, implies
$\rk(\7s)\le3$. We work out the various cases separately.

\noindent {\fett $\rk\7s=3\,$:} The complexification $\7s^{\CC}$ of
$\7s$ is one of the Lie algebras $\7{sl}(4,\CC)$, $\7{so}(7,\CC)$ or
$\7{sp}(3,\CC)$. The latter two can be immediately ruled out since a
glance at the corresponding Satake diagrams shows that every proper
parabolic subalgebra of them is at least of codimension 4. In the
remaining case $\7s^{\CC}\cong \7{sl}(4,\CC)$ only the normal real
form $\7{sl}(4,\RR)$ has a parabolic subalgebra of codimension 3. (A
glance at the Satake diagrams for the remaining non-compact real forms
of $\7{sl}(4,\CC)$ excludes further possibilities).

\noindent {\fett $\rk\7s=2\,$ or $1$:} The rank conditions imply that
$\7s$ is the underlying real Lie algebra of $\7{sl}(2,\CC)$ or a real
form of $\7{sl}(3,\CC)$ or $\so(5,\CC).$ The Lie algebra $\so(1,3)$ as
well as every non-compact real form of $\7{sl}(3,\CC)$ or $\so(5,\CC)$
contains a parabolic subalgebra of codimension less or equal 3. Since
$\dim \7{sl}(2,\RR)=\dim\su(2)=3$, all simple Lie algebras of rank one
also contain parabolic subalgebras of the required codimension.

\noindent In order to further reduce the list of possibilities for
$\7s$ we have to look at the particular real forms obtained in Claim 1
and 2 but not in the list of Lemma \ruf{LA} more closely:

\noindent {\fett Elimination of $\7s\cong \7{sl}(4,\RR)$:} Up
to an automorphism of $\7{sl}(4,\RR)$ there is only one such parabolic
subalgebra $\7h$ of codimension $3$. Let $\7t=\7a\subset \7h$ be the
split Cartan subalgebra and $\7h=\7h^{\r}\ltimes \7h^{-\n}$ the
decomposition as in \Ruf{FG}. Note that here the reductive factor
$\7h^{\r}\cong \gl(3,\RR)$ acts irreducibly on $\7h^{\pm\n}\cong
\RR^3$. The only possibility for the flag $\7s_o\subset \7s_\7F\subset
\7s_\7H\subset \7s$, which cannot be trivially excluded, is when
$\7s_o\subset \7s_\7F=\7h\subsetneq \7s_\7H\subsetneq \7s$ and
$\codim_\7h(\7s_o)\le 2.$ This cannot be true since
$\dim\qqu{\7s_\7H}{ \7s_\7F}=2$, $[\7s_\7F:\7s_\7H]\subset \7s_\7H$
(condition \22) but $\7h=\7s_\7F$ acts irreducibly on the
3-dimensional space $\qqu{\7s}{\7s_\7F}\cong \7h^{\n}.$

\noindent{\fett Elimination of $\7s\cong\7{su}(1,2)\,$:} In this case
of real rank 1, there exists up to conjugacy only one parabolic
subalgebra $\7h$: This is the minimal one $\7h=\7h_\emptyset
=\7m\oplus \7a\oplus\7n$ which is solvable with $\codim_{\7s}\7h
=3$. The reductive part $\7h^{\r}=\7m\oplus \7a=:\7t$ of $\7h$ is a
maximally split Cartan subalgebra. As before, the only situation which
cannot be trivially disposed of is when $\7s_o\subset \7s_\7F=\7h$ and
$c:=\codim_\7h(\7s_o)\in \{0,1,2\}.$ Recall that $\7h$ yields the
decomposition $\7s=\7h^{\n}\oplus \7t\oplus \7h^{-\n}$. If $c=0,$ that
is $\7s_o=\7h$, then $\7h^{\n}\oplus\rad(\7g)$ would be a proper
locally transitive subalgebra of $\7g$ in contradiction to assumption
\26$_{2}$. If $c=1$ then either $\7t\cap \7s_o\ne \7t$, and then
$\7h^{{\opp}}\oplus\rad(\7g)$ would be a proper locally transitive
subalgebra of $\7g$ (again contradicting \26$_{2}$) or $\7t\subset
\7s_o.$ But taking into account the particular structure of
$\7h^{-\n}=\7s_{-\lambda}\oplus \7s_{-2\lambda}$ this also does not
occur since otherwise $\7s_o\cap \7h^{-\n}$ would be a $\7t$-stable
2-dimensional subalgebra, which is impossible since $\ad(\7t)$ acts
irreducibly on the 2-dimensional root spaces $\7s_{\pm\lambda}$ and
$[\7s_{\lambda},\7s_{\lambda}]=\7s_{2\lambda}.$ \nline It remains the
case $c=2,$ that is, $\7s=\7g\cong \su(1,2)$ is locally
transitive. Then $\7l=\7s^{\CC}\cong\7{sl}(3,\CC)$ and $\7q$ is a
subalgebra of complex dimension $5$. Consequently $\7q$ is contained
in a maximal (6-dimensional) parabolic subalgebra
$\7h\cong\7{gl}(2,\CC)\ltimes\CC^2$ of $\7l$, that is, either $\7q$
coincides with a Borel subalgebra $\7b$ or is conjugate to a
subalgebra $\7j\cong \7{sl}(2,\CC)\ltimes \CC^2$. In both cases the
subgroups $Q$ of $L=\SL(3,\CC)$ corresponding to $\7b$ or $\7j$ are
closed and the underlying CR-germ $(M,o)$ is locally CR-equivalent to
an $\SU(1,2)$-orbit either in $L/B\cong\FF(\CC^3)$, the complex
manifold of full flags in $\CC^3$, or in the $\CC^*$-principal bundle
$L/J$ over $\PP_2(\CC)$. A direct check shows that in both cases there
do not exist 2-nondegenerate $\SU(1,2)$-orbits.

\noindent{\fett Elimination of $ \7s\cong\7{sl}(3,\RR)$.} Then
$3\le\dim\7s_o\le6$ holds except for the trivial case $\7s_o=\7s$.  If
$\dim\7s_o=3$ and hence $\7s$ is locally transitive, then as in the
previous situation each CR-germ $(M,o)$ associated with a CR-algebra
$(\7{sl}(3,\RR),\7q)$ is locally CR-equivalent to an
$\SL(3,\RR)$-orbit either in $\SL(3,\CC)/B$ or $\SL(3,\CC)/J$ (as
discussed above, with $\SL(3,\RR)$ in place of $\SU(1,2)$). Again,
none of these orbits is 2-nondegenerate. It remains the case $\dim
\7s_o\ge 4$. But then there always exists a parabolic (proper)
subalgebra $\7h\subset\7s$ with $\7s_o+\7h=\7s$, that is,
$\7h\oplus\rad(\7g)$ is a proper locally transitive subalgebra of
$\7g$ excluding the case $\7s^{\CC}\cong\7{sl}(3,\CC)$.

\noindent{\fett Elimination of $ \7s\cong\so(1,4)$.} There exists up
to conjugacy a unique parabolic subalgebra
$\7h=\7h^{\r}\ltimes\7h^{-\n}\subset\7s$ of codimension $3$, and we
have to investigate the cases $\7s\supsetneq \7s_\7H\supsetneq
\7s_\7F=\7h\supset \7s_o$ only. A close look at the minimal (and
maximal proper) parabolic subalgebra $\7h=\7m\oplus \7a\oplus
\7s_{\lambda}$ shows that $\7m\cong \so(3)$ and $\7m$ acts irreducibly
on the 3-dimensional nilpotent ideal $\7h^{-\n}=\7s_{-\lambda}$.
Consequently $\7h$ acts irreducibly on $\7s/\7h =\7s/\7s_\7F\cong
\RR^3.$ This leads to a contradiction as
$[\7s_{\scriptscriptstyle\7F},\7s_{\scriptscriptstyle\7H}]\subset
\7s_{\scriptscriptstyle\7H}$, i.e.,
$\qqu{\7s_{\scriptscriptstyle\7H}}{\7s_{\scriptscriptstyle\7F}}$ would
be a 2-dimensional stable subspace. The proof of Lemma \ruf{LA} is now
complete.\qed

\medskip
In Theorem 6 of \Lit{BELO} it is claimed that for every
$2$-nondegenerate real-analytic hypersurface $M\subset\CC^{3}$ the Lie
algebra $\hol(M,a)$ has dimension $\le11$ at every point. Using this
result would save few arguments in the proof of Lemma \ruf{LA}.
Instead, we preferred to present a self-contained proof of the
proposition.

\KAP{Elf}{The cases $\7s\cong\so(2,3)$ and $\7s\cong\so(1,3)$}

We continue the proof of Theorem II.  So far we have proved that for
the CR-algebra $(\7g,\7q)$ under consideration the simple factor $\s$
of $\7g^{\s}$, see \Ruf{SS}, can only be isomorphic to one of the
simple Lie algebras listed in Lemma \ruf{LA}. In this section we show
that from these the possibilities $\7s\cong \so(2,3)$ and $\7s\cong
\so(1,3)$ cannot occur. Here and in the following upper case roman
numerals refer to the conditions around \ruf{FD}.

We like to mention that for the tube $\5M$ over the future light cone
$\hol(\5M,o)\cong\so(2,3)$ holds and that there exists a copy of
$\so(1,3)$ in $\hol(\5M,a)$ that is also locally transitive. Since
these Lie algebras have dimensions 10 and 6 and since, on the other
hand, there are transitive subalgebras of $\hol(\5M,o)$ of dimension 5
these two Lie algebras do not satisfy the minimality condition \26.
In the following we consider both cases separately:

\noindent{\fett $\7s\cong\so(2,3)$.} The only Lie subalgebras $\7h$ in
the normal real form $\so(2,3)$ with $\codim_\7s\7h\le 3$ are the
3-codimensional maximal parabolic subalgebras. We need to take a
closer look at the structure of these subalgebras. There are up to
isomorphisms of $\7s$ only two such parabolic subalgebras: If
$\Pi(\7a)=\{\alpha,\beta\}$ is a basis of the root system $\Phi(\7a)$
($\7a$ split Cartan subalgebra, $\alpha$ long, $\beta$ short) then

\smallskip $\displaystyle \eqalign{\7h_1:&=(\7a\oplus \7s_\alpha\oplus
\7s_{-\alpha}) \;\;\oplus \;\; \7s_{-\beta}\oplus
\7s_{-\alpha-\beta}\oplus \7s_{-\alpha-2\beta}=\7h_1^{\r}\ltimes
\7h_1^{-\n} \cr \7h_2:&=(\7a\oplus \7s_\beta\oplus
\7s_{-\beta})\;\;\oplus\;\; \7s_{-\alpha}\oplus
\7s_{-\alpha-\beta}\oplus \7s_{-\alpha-2\beta}=\7h_2^{\r}\ltimes
\7h_2^{-\n}}$

\smallskip\noindent are representatives of the corresponding isomorphy
classes. The only instance where the filtration $\7s_o\subset
\7s_{\scriptscriptstyle\7F}\subset \7s_{\scriptscriptstyle\7H}\subset
\7s$ could be nontrivial (recall that $\7s_o,\;
\7s_{\scriptscriptstyle\7F}$ are subalgebras,
$\7s_{\scriptscriptstyle\7H}$ is a
$\ad(\7s_{\scriptscriptstyle\7F})$-stable subspace) arises when
$\7s_o\subset \7s_{\scriptscriptstyle\7F}=\7h_j\subset
\7s_{\scriptscriptstyle\7H}\subset \7s$ for $j=1,2.$ The possibility
`$\,\7h_2=\7s_{\scriptscriptstyle \7F}$' cannot occur since the
adjoint representation of $\7h_2^{\r}$ on $\7h^{\n}_2\cong
\qu{\7s}{\7s_{\scriptscriptstyle\7F}}$ is irreducible, contradicting
the existence of a 2-dimensional
$\ad(\7s_{\scriptscriptstyle\7F})$--stable subspace
$\qqu{\7s_{\scriptscriptstyle\7H}}{\7s_{\scriptscriptstyle\7F}}.$ It
remains the possibility ${\7s_{\scriptscriptstyle\7F}}=\7h_1\supset
\7s_o$. From now on $\7h:=\7h_1$ and we analyze the various
possibilities for $\dim {\7s_{\scriptscriptstyle\7F}}-\dim
\7s_o\in\{0,1,2\}.$ The equation ${\7s_{\scriptscriptstyle\7F}}=\7s_o$
contradicts condition \26$_{2}$ since in that case the proper
subalgebra $(\7h^{\n}\oplus \7s')\ltimes \rad(\7g)$ would be
transitive on $(M,o).$ The case when $\7s_o$ is of codimension 1 in
${\7s_{\scriptscriptstyle\7F}}=\7h$ can also be ruled out: Either
$\7h^{-\n}\subset \7s_o$ and then $(\7h^{\opp}\oplus \7s')\ltimes
\rad(\7g)$ would be a transitive proper subalgebra of $\7g$, or the
intersection of $\7s_o$ with $\7h^{-\n}=(\7s_{-\beta}\oplus
\7s_{-\beta-\alpha}) \oplus \7s_{-\alpha-2\beta}$ is a 2-dimensional
subalgebra. In such a case the image $\pi(\7s_o)$ of the projection
$\pi:\7h=\7h^{-n}\rtimes \7h^{\r}\to \7h^{\r}$ coincides with
$\7h^{\r}.$ But this also leads to a contradiction: Neither the
intersection $\7s_o\cap \7h^{-n}$ can coincide with
$\7s_{-\beta}\oplus \7s_{-\beta-\alpha}$ (since it is not a
subalgebra), nor $(\7s_{-\beta}\oplus \7s_{-\beta-\alpha})\cap \7s_o$
can be 1-dimensional (since $\7h^{\r}$ acts irreducibly on
$(\7s_{-\beta}\oplus \7s_{-\beta-\alpha})$). Finally we are left with
the case $\dim {\7s_{\scriptscriptstyle\7F}}-\dim \7s_o=2, $ i.e.,
$\7s$ is transitive on $(M,o).$ Thus $\7g=\7s$ by assumption \26.  But
then $\dim\7g_{o}=5$ and there always exist proper subalgebras
$\7g'\subset \7g$ with $\7g'+\7g_{0}=\7g$, a contradiction to
condition \26$_{2}$.

\medskip \noindent{\fett $\7s\cong\so(1,3)\cong\7{sl}(2,\CC) $.} We
work out this case by investigating various possibilities for $\dim
\7s_o.$ Assume first that \nline $\Bullet\;\dim \7s_o\ge 3.$ We claim
that then there exists a solvable subalgebra $\7r\subset \7s$ such
that $\7s_o+\7r=\7s:$ The case $\dim \7s_o\ge 4$ is easily settled as
all 4-dimensional subalgebras in $\7s$ are maximal, i.e., they are
Borel subalgebras of $\7{sl}(2,\CC)$ and consequently have nilpotent
complementary subalgebras. There do not exist (real) subalgebras of
$\7{sl}(2,\CC)$ of dimension 5. If $\dim \7s_o=3$ then $\7s$ is either
semisimple or solvable. In the semisimple case we work with an
explicit matrix realization $\7s\subset \CC^{2\times 2}$: Either
$\7s_o \cong \su(2)$ or $\cong\su(1,1),$ i.e.,
$$
\7s_o\cong \Big\{\!\!\pmatrixe{it & \epsilon\overline z \cr z &
-it}\!\!:t\in \RR,z\in \CC
\Big\}\Steil{for}\epsilon=1\steil{or}\epsilon=-1\;.
$$ In both cases the upper triangular Borel subalgebra $\7b^+\subset
\7{sl}(2,\CC)\subset \CC^{2\times 2}$ forms a linear complement of
$\7s_o.$ This cannot happen, since then condition \26$_{2}$ would be
violated. \nline If $\7s_o$ is solvable then $\7s_o$ is contained as
certain 1-codimensional (real) subalgebra in a (complex) Borel
subalgebra $\7b=\7t\ltimes \7b^{\nil}$ of $\7s\cong\7{sl}(2,\CC)$. We
claim that $\7s_o\supset \7b^{\nil}:$ Otherwise (i.e., if
$\dim^{}_{\RR }\7s_o\cap \7b^{\nil}=1$) we would have
$\pi(\7s_o)=\7t$, where $\pi:\7b\to \7t$ is the projection
homomorphism, and consequently $\7s_o$ would contain a certain complex
Cartan subalgebra $\7t'$ of $\7b$ which acts $\RR$-irreducibly on
$\7b^{\nil}$. \nline However, from our claim it follows that the
opposite Borel subalgebra is a complementary subspace of $\7s_o$ in
$\7s$. The case `$\dim \7s_o=3$' now is completely ruled out. Next, we
investigate the case \nline {$\Bullet\;\dim \7s_o=2.$} It follows that
$\7s_o$ is solvable and either complex or totally real. Since it is
immediate that every complex subalgebra in $\7{sl}(2,\CC)$ has a
complementary subalgebra (simple check), it remains only to deal with
the totally real case, i.e., with $\7s_o$ being a real form of a Borel
subalgebra $\7b\subset \7{sl}(2,\CC).$ Let $\tau:\7b\to \7b$ be the
conjugation with $\7b^\tau=\7s_o.$ It is well-known that there exists
a $\tau$-stable Cartan subalgebra $\7t\subset \7b$ and consequently we
have the $\tau$-stable decomposition $\7b=\7t\ltimes
[\7b,\7b]=:\7t\ltimes \7n.$ Select $\4h\in\7t,\;\4e\in \7n$ such that
$[\4h,\4e]=2\4e.$ By construction $\tau(\4h)=a\cd \4h$ and
$\tau(\4e)=b\cd \4e$ for suitable $a,b\in\CC$. Since $\tau$ is an
automorphism, $a=1.$ Since $\tau$ is an involution, $|b|=1.$ This
shows that up to a conjugation we may assume that $\7b$ is the
upper-triangular Borel subalgebra of $\7{sl}(2,\CC)$ and the real
forms $\7s_o\subset\7b^+$ have the following realizations:
$$
\7s_o\cong \Big\{\!\!\pmatrixn{t & sc\cr 0& -t}:t,s\in
\RR\Big\}\steil{for some}c\in\CC^{*}\;.
$$ For every $c\in\CC^{*}$ at least one of the Borel subalgebras
$\CC\big({0\kern5pt 1 \atop 1\kern5pt 0}\big)\oplus \CC\big({1\kern2pt
-1 \atop 1\kern2pt -1}\big)$ or $\CC\big({0\kern2pt -i \atop i\kern7pt
0}\big)\oplus \CC\big({1\kern7pt i \atop i\kern2pt -1} \big)$ then is
complementary to $\7s_o$ in $\7{sl}(2,\CC).$ Finally, we need to deal
with the case \nline $\Bullet\;\dim \7s_o= 1,$ i.e., $\7s$ itself is
locally transitive and thus $\7g=\7s$ by the minimality condition
\26$_{2}$. Consider the following matrix realization:\vskip-20pt
$$
\displaylines{\7l=\7{sl}(2,\CC)\times\7{sl}(2,\CC) \;\subset
\;\CC^{2\times 2}\!\times \CC^{2\times 2} \cr
\sigma(\4x,\4y)=(\overline{\4y},\overline{\4x})\,, \qquad
\7g=\{(\4x,\overline{\4x}):\4x\in \7{sl}(2,\CC)\}\subset \7l\;.}
$$
The 3-dimensional complex subalgebra $\7q\subset\7l$ is either simple
or solvable: \nline {\sl $\7q$ simple:} Then $\7q\cong \7{sl}(2,\CC)$
and $\7q$ is either one of the two factors of $\7{sl}(2,\CC)\times
\7{sl}(2,\CC)$ or $\7q$ is conjugate in $\7l$ to $\7g$. The first case
can be ruled out immediately since then $\7q\cap \7g=0=\7g_o$ which is
absurd.  In the second case, all simple subalgebras $\7q\subset \7l$
which are not ideals are conjugate to each other. We may select the
particular subalgebra $\7q=\{(\4x,-\4x^t):\4x\in
\7{sl}(2,\CC)\}\subset \7l$, where $\4x^t$ is the transpose of
$\4x$. The corresponding Lie subgroup $Q=\{(\4q,(\4q^{-1})^t):\4q\in
\SL(2,\CC)\}\subset L:=\SL(2,\CC)\times\SL(2,\CC)$ is closed and the
map $(\4x,\4y)\mapsto \4x\cd \4y^t$ identifies the quotient with the
affine quadric $\SL(2,\CC)\subset\CC^{2\Times2}$ on which
$S:=\SL(2,\CC)$ (considered as real Lie group) acts by the holomorphic
transformations $(\4s,\4z)\mapsto \4s\4z\overline{\4s}{}^t$. Hence, in
this situation every CR-germ associated with a CR-algebra $(\7s,\7q)$
is globalizable. It is well-known (see, for instance \Lit{FEKA}) that
the hypersurface $S$-orbits in $L/Q$ are either Levi-nondegenerate or
locally CR-equivalent to the tube $\5M$ over the light cone. However,
$\dim \7g=6$ and in this case $(\7g,\7q)$ satisfies \21 -- \25,
\26$_{1}$ and \26$_{2}$, but not \26.\nline {\sl $\7q$ solvable:} We
will show that also this case contradicts the fundamental assumption
\ruf{FD}: \nline If $\7l_o=\CC(\4t,\overline{\4t})$ is ad-semisimple
(and without loss of generality $\4t=\big({t \kern1em \atop \kern.7em
-t}\big )$ for some $t\in \CC^{*}$) then $\7l_o$ is a regular torus in
$\7l.$ Consequently, the centralizer $C_{\7l}(\7l_o)=:\7t_1\times
\7t_2$ is a $\sigma$-stable Cartan subalgebra. Either
$C_{\7l}(\7l_o)\subset \7q$ (but then $\7q+\sigma\7q$ is of
codimension 2 in $\7l$) or $C_{\7l}(\7l_o)\cap \7q=\7l_o.$ In the
latter case denote by $\7l^{\pm\alpha_1},$ $\7l^{\pm\alpha_2}$ the
root spaces with respect to $\7t_1\times \7t_2.$ A direct check shows
that $\7q$ is the direct sum of $\7l_o$ with two further root spaces
(also if $\overline{\4t}=\4t$). Since $\7q$ is solvable, there are 4
possibilities of choosing such pairs of root spaces. In all 4 cases
either $\7q+\sigma\7q$ is too small or $\7f=\7l_o,$ i.e., the
corresponding CR-germ is Levi nondegenerate. \nline It remains to
discuss the case when $\7l_o=\CC(\4n,\overline{\4n})$ is
ad-nilpotent. Since $\7q$ is solvable, it is contained in a Borel
subalgebra $\7b$ of $ \7l.$ Consequently, $\7l_o\subset
\7b^{\nil}=C_\7l(\7l_o)\cong \CC\4n \times \CC\overline{\4n}.$ Let
$\pi:\7b\to \7b/\7b^{\nil}$ be the canonical projection. The image
$\pi(\7q)$ cannot be surjective since there is no 3-dimensional
subalgebra $\CC(\4n,\overline{\4n})\subset \7q \subset \7b$ with this
property. It remains only the possibility $\7b^{\nil}\subset \7q$ but
then $\7q\cap\sigma\7q \supset \7b^{\nil},$ which is too big.

Summarizing, we have for the CR-algebra $(\7g,\7q)$ satisfying
\ruf{FD} that $\7g^{\s}$ must be a finite direct sum of copies of
$\7{sl}(2,\RR)$ and $\su(2)$. In the next section we direct our
attention to factors of type $\7{sl}(2,\RR)$.

\KAP{Zwoelf}{The case $\7s\cong\7{sl}(2,\RR)$}

In this section we continue the proof of Theorem II and consider only
CR-algebras $(\7g,\7q)$ subject to \ruf{FD}, for which the simple
factor $\7s$ of $\7g^{\s}$ is isomorphic to $\7{sl}(2,\RR)$, compare
\Ruf{SS}. As already proved, the remaining simple factors in $\7s'$
(if there are any) are isomorphic to $\su(2)$ or $\7{sl}(2,\RR)$.

From $\hol(\5M,o)\cong\so(2,3)$ for the light cone tube $\5M$ it is
clear that $\hol(\5M,o)$ contains copies of
$\so(2,2)\cong\7{sl}(2,\RR)\Times\7{sl}(2,\RR)$. Fixing a subalgebra
$\7r\subset\7{sl}(2,\RR)$ (which is necessarily non-abelian) the
5-dimensional Lie algebra $\7{sl}(2,\RR)\times\7r$ can be embedded
into $\hol(\5M,o)$, and this can be done in such a way that the image
is a transitive subalgebra. Therefore, the simple factor
$\7{sl}(2,\RR)$ of $\7g^{\s}$ cannot avoided in the classification
proof. However, we will show that this factor only occurs in the
instance described above, that is, in connection with $\5M$.

Our first result is that in $\7g^{\s}$ the simple factor
$\7{sl}(2,\RR)$ can occur at most once. Here and in the following we
repeatedly use the basic fact that each proper subalgebra of
$\7{sl}(2,\RR)$ has a solvable complementary subalgebra.  Furthermore,
through this subsection we denote by $\pi:\7g=(\7s\times \7s')\ltimes
\rad(\7g) \to \7s$ the canonical projection. We analyze various
possibilities for the image of the isotropy subalgebra $\7g_{o}$ under
$\pi$. \nline If $\pi(\7g_o)\ne 0$ then there exists a solvable
complement $\7r\subset \7s$ to $\pi(\7g_o)$ and $(\7r\times
\7s')\ltimes\rad(\7g)$ is a proper transitive subalgebra of $\7g$
violating \26$_{2}$.\nline If $\pi(\7g_o)= 0$ then $\7g_o\subset
\7s'\ltimes \rad(\7g)$ and is there of codimension 2. It follows that
there is no factor $\7h\cong\7{sl}(2,\RR)$ in $\7s'$: Otherwise,
counting dimensions we would have $\dim \7h \cap\7s_o\ge 1$ which
implies that there is a proper transitive subalgebra of $\7g$,
violating \26$_{2}$. This proves that $\7s'$ is a product of factors
which all are isomorphic to $\su(2)$. Finally, we show that $\7s'$
consists of at most {\sl one} such simple factor: Indeed, suppose that
$\7s'=\7h\times\7s''$ for some ideal $\7h\cong\su(2)$ of $\7s'$.
Denote by $\pi_{\7h}:=(\7s\times\7h\times\7s'')\ltimes
\rad(\7g)\to\7h$ the canonical projection. Then $\pi_{\7h}(\7g_{o})$
has codimension $\le2$ in $\7h$. Since $\su(2)$ does not have a
subalgebra of dimension $2$, the dimension of $\pi_{\7h}(\7g_{o})$ can
only be $1$ or $3$. But dimension $3$ violates \26$_{2}$ since then
$\7g=\7g_{o}+\ker(\pi_{\7h})$. It remains the case
$\dim\pi_{\7h}(\7g_{o})=1$. Then $\pi_{\7h}(\7g_o)$ is a torus in
$\7h$ and $\pi(\7g_o)=\7h\cap \7g_o$. Since $\7g':=\7s\Times \7h$ is a
transitive subalgebra of $\7g$ we must have $\7g'=\7g$ by
\26$_{2}$. But this is not possible:

\Lemma{LK} $\7s\cong\7{sl}(2,\RR)$ implies
$\7g\;=\;\7s\!\ltimes\!\rad(\7g)$ and $\7g_o\subset\rad(\7g)\,.$

\Proof By the above discussion we only have to rule out the case
$\7g=\7s\times\7s'$ with $\7g_{o}\subset\7s'=\su(2)$. The key point
here is that the isotropy subalgebra $\7g_o$ is toral in $\su(2)$ and
that there exists a unique $\ad(\7g_o)$--stable subspace
$\7p^\star\subset\su(2)$ on which the adjoint representation of
$\7g_o$ is irreducible. Let $\pi_{\su}:\7s\oplus \su(2)\to \su(2)$ be
the projection onto the second factor. Either $\pi_{\su}(\7F)=\su(2)$
or $\pi_{\su}(\7F)=\7g_o.$ In the first case, we actually have
$\7F=\su(2)$ as $[\7g_o,\7F]\subset \7F$ and $[\7s,\7g_o]=0$. But this
cannot be true: Independently of what exactly $\7H=W\oplus \su(2)$
would be, we always would have $[\7H,\7F]\subset \7F$. But this
violates the nondegeneracy condition \25.  \nline It remains the case
$\pi_{\su}(\7F)=\7g_o$, i.e., $\7F=\7b\oplus \7g_o$, where
$\7b=\7F\cap \7s\subset\7s$ is a 2-dimensional real subalgebra. By a
dimension argument $\dim \7H\cap \su(2) \ge 2.$ But this intersection
must be $\ad(\7g_o)$-stable, which implies $\7H=\7b\oplus \su(2)$,
i.e., $\7H$ is a subalgebra of $\7g=\7s\oplus \su(2)$.  Clearly, this
violates condition \23.\qed

\smallskip In the next Lemma we show that the semidirect product
$\7s\ltimes\rad(\7g)$ in \ruf{LK} actually is a direct product
$\7s\Times\rad(\7g)$ and that the radical has dimension 2, i.e.  $\7g$
has dimension 5. Recall the obvious fact that, up to isomorphy, there
exist precisely two Lie algebras of dimension 2, the abelian Lie algebra
$\RR^{2}=\RR\Times\RR$ and a non-abelian one.

\Lemma{FM} $\7s\cong\7{sl}(2,\RR)$ implies $\7g\;=\;\7s\Times\7r$ with
$\7r:=\rad(\7g)$ being non-abelian and of dimension 2.

\Proof Let $\pi:\7g\to \7s$ be the canonical projection. We use the
same symbol also for the complex extension $\7l\to
\7s^{\CC}$. Because of \21 and $\7g_o\subset\7r$ the image
$\pi(\7F)$ has dimension $\le2$.

\noindent $\Bullet\dim \pi(\7F)=0$ implies $\7F=\7r,$ i.e., $\7F$
is an ideal in $\7g.$ But this violates \25.

\noindent $\Bullet\dim \pi(\7F)=1$ implies
$\pi(\7f)=\pi(\sigma\7f)$ and $\dim\pi(\7f)=1$. Consequently
$\pi(\7q),\pi(\sigma\7q)$ are 2-dimensional (Borel) subalgebras
which generate $\7{sl}(2,\CC)$ as a linear space (otherwise $(M,o)$
would be Levi flat). It follows $\pi^{-1}(\pi(\sigma\7q))\cap
\7q=\7f.$ But this implies $[\7f,\sigma\7q]\subset \7f+\sigma\7q,$
a contradiction to \25.

\noindent $\Bullet\dim \pi(\7F)=2.$ Note that $\pi(\7f),\;
\pi(\sigma\7f)$ are 1-dimensional since $\7l_o\subset \7r^{\CC}$.
Since $\pi(\7f)+\pi(\sigma\7f)$ is a 2-dimensional subalgebra,
$\7t:=\pi(\7f),$ $\7t':=\pi(\sigma\7f)=\sigma\7t$ are tori as follows
with the elementary structure theory of $\7{sl}(2,\CC)$. The case
`$\pi(\7f)=\pi(\7q)$' can be excluded since otherwise we would have
$\pi(\7q+\sigma\7q)=\pi(\7f+\sigma\7f)$ which is a subalgebra, a
contradiction to \23. Hence, the only possibility remaining is
$\pi(\7f)\ne \pi(\7q)$. This implies $\7q\cap \rad= \sigma\7q\cap
\rad=\7l_o$. Furthermore, $\7l_o$ is an ideal in $\7q,\sigma\7q,$ and,
since $[\7q,\sigma\7q]=\7l,$ even an ideal of $\7l.$ By the
effectivity assumption \26$_{1}$, this implies $\7g_o=0,$ i.e., $\7g$
has dimension $5$ and thus
$$\7g\cong\7{sl}(2,\RR)\ltimes_\rho\RR^2\Steil{ or }\7g\cong
\7{sl}(2,\RR)\times \7r\steil{with}\dim\7r=2\,,\Leqno{EM}$$ where
$\rho:\7{sl}(2,\RR)\to \End(\RR^2)$ is the canonical inclusion.

\Claim{NO} The first case in \Ruf{EM}, that is
$\7g\cong\7{sl}(2,\RR)\ltimes_\rho\RR^2$, cannot occur.

\ProofC Let $\pi:\7l\to \7{sl}(2,\CC)$ be the canonical projection and
$\sigma:\7l\to \7l$ the complex conjugation defining the real form
$\7g$ of $\7l$. The possibility \nline $\Bullet\;\pi(\7q)=0\,$
($=\pi(\sigma\7q)\,$) can be excluded since then
$\7q=\rad(\7l)=\sigma\7q$, violating \23. Also

\noindent $\Bullet\dim \pi(\7q)=1$ can be ruled out: Then
$\pi(\7q)\ne\pi(\7f)$ since otherwise
$\pi(\7q+\sigma\7q)=\pi(\7f+\sigma\7f)$ would be a subalgebra,
violating \23. Hence, $\7f=\7q\cap \rad(\7l)$. But this contradicts
\25 since then $[\7f,\sigma\7q]\subset\rad(\7l)=\7f\oplus
\sigma\7f$.

\noindent $\Bullet\dim \pi(\7q)=2$ (i.e., $\7q\cap \rad(\7l)=0$) is
the most involved case. Then $\pi(\7q)=\7b$ is a Borel subalgebra of
$\7{sl}(2,\CC)$. To rule out also this case we need some
preparations. First realize $\7l$ as \vskip-20pt
$$
\7{sl}(2,\CC)\ltimes_{\rho} \CC^2=\{(X,w): X\in \7{sl}(\CC^2), \;w\in
\CC^2 \}\steil{with}[(X,w),(Y,u)]=([X,Y]\;, \;Xu-Yw)
$$ and complex conjugation $\sigma$ given by $(X,w)\mapsto(\overline
X,\overline w)$. Since $(M,o)$ is not Levi flat (compare \23), we
necessarily have $\7b+\sigma\7b=\7{sl}(2,\CC),$ and the 1-dimensional
intersection $\7b\cap \sigma\7b$ is a toral subalgebra
$\7t\subset\7{sl}(2,\CC)$ (we use here the well-known fact that the
[$\sigma$--stable] intersection of any two Borel subalgebras contains
a [$\sigma$--stable] Cartan subalgebra). In the next two paragraphs we
recall some elementary facts from the representation theory of
$\7{sl}(2,\CC)$ which we need to complete our proof.

\Joker{FN}{$\sigma$-adapted $\7{sl}_2$-triples} Let $\4h\in \7t\subset
\7{sl}(2,\CC)$ be the element for which $[\4h,E]=2E$ and $[\4h,F]=-2F$
for every $E\in \7b^{\nil},$ $F\in
\sigma(\7b^{\nil})=(\sigma\7b)^{\nil}.$ Since $\sigma$ interchanges
the eigenspaces of $\ad(\4h),$ we have $\sigma(\4h)=-\4h$. There are
crucial technical points here: We claim that there exists $\4e\in
\7b^{\nil}$ such that $[\4e,\sigma(\4e)]=\4h,$ i.e., $(\4e,\4h,\4f)$
with $\4f:=\sigma(\4e)$ being an $\7{sl}_{2}$-triple in
$\7{sl}(2,\CC)$ (i.e.,
$[\4h,\4e]=2\4e,\;[\4h,\4f]=-2\4f,\;[\4e,\4f]=\4h$). \nline{\bf
Construction:} We have to be careful here about signs: Since $\sigma$
defines a {\sl non-compact} real form of $\7{sl}(2,\CC),$ the
nondegenerate Hermitian 2-form $\kappa(\cdot\,,\sigma(\,\cdot\, ))$
has precisely one negative eigenvalue, where $\kappa$ denotes the
Killing form. Since $\kappa(\4h,\sigma(\4h))=\kappa(\4h,-\4h)<0$, it
follows $\kappa(E,\sigma(E))>0$ and consequently $[\lambda
E,\sigma(\lambda E)]=\4h$ for an appropriately chosen $\lambda\in
\CC^*$ (keeping in mind the general formula $[E,F]=\kappa(E,F)H$ where
$H$, the coroot, is a positive multiple of $\4h$). Define then
$\4e:=\lambda E$ and $\4f:=\sigma(\4e)$. Each $\7{sl}_{2}$-triple
$\4e',\4h',\4f'\in \7s^{\CC}\cong \7{sl}(2,\CC)$ with
$\sigma(\4e')=\4f'$ is called $\sigma$--{\sl adapted} in the
following.

\Joker{FO}{\hskip-3pt} Complexifying $\7{sl}(2,\RR)\ltimes_\rho \RR^2$, we
briefly discuss how the 2-dimensional abelian radical $\CC^2,$
considered as an $\7{sl}(2,\CC)$-module, is related to the real
structure. Let a $\sigma$-adapted $\7{sl}_2$-triple $\4e,\4h,\4f$ be
given.  Let $\CC^2=V_+\oplus V_-$ be the decomposition into $(\pm 1)$
$\4h$-eigenspaces. They are interchanged by $\sigma$. We claim that
there exists a $v_+\in V_+$ with
$v_-:=\sigma(\4e)v_+=\sigma(v_+)\ne0$: Indeed, choose $w_+\in V_+\sm
\{0\}$ arbitrarily. Then there is a $c\in \CC^*$ with $\4f
w_+=c\sigma(w_+)$, and a direct check shows $|c|=1.$ Choose a
$b\in\CC$ with $b^{2}=\overline c$ and put $v_+:=b\cd w_+$.  \nline We
use the linear basis $v_+,v_-$ of the radical $\rad(\7l)\cong\CC^2$ in
the following computations.

\noindent We now resume the proof of Claim \ruf{NO}. For short,
write $\7r^{\CC}  :=\rad(\7l)$
for the abelian radical.
Since  $\7q\cap \7r^{\CC}=0$, we can write
$\7q=\CC\cd(\4e,w_1)\oplus \CC\cd(\4h,w_2)$ for suitable $w_1,w_2\in
\CC^2.$ Clearly, the $w_j$'s are not arbitrary: Write
$w_1=\lambda_1v_+ +\mu_1v_-$ and $w_2=\lambda_2v_+ +\mu_2v_-.$ Note
that $w_2\ne -\overline w_2$ (otherwise $\7q\cap\sigma\7q\ne 0$),
i.e., $\lambda_2\ne -\overline \mu_2.$ The fact that $\7q$ is a
subalgebra imposes a further condition on the coefficients
$\lambda_j,\mu_j:$ A simple computation shows that $\mu_1=0$ and
$\mu_2=-\lambda_1,$ i.e., for each $\lambda,\mu \in \CC$ with
$\lambda\ne \overline \mu$ we have the two 2-dimensional complex
subalgebras $\7q=\7q_{\lambda,\mu}$ and $\sigma\7q:$

\noindent \hfil $\7q=\CC\cd(\4e,\lambda v_+)\oplus
\CC\cd(\4h,\mu v_+-\lambda v_-),\quad ~ \sigma\7q=\CC\cd(\4f,\overline
\lambda v_-)\oplus \CC\cd(-\4h,-\overline \lambda v_++\overline\mu
v_-)\,. $

\noindent Recall the definition $\7f=\{u\in
\7q:[u,\sigma\7q]\subset\7q+\sigma\7q\}$ from \ruf{CA}.iv. A
straightforward calculation shows that in all cases $\7f=0$
holds. This contradicts \24 and proves \ruf{NO}. \qd

\noindent We proceed the proof of Lemma \ruf{FM} by restricting  the
radical of $\7g\cong\7{sl}(2,\RR)\times\7r$ further:

\Claim{FP} The Lie algebra $\7g$ cannot be isomorphic to
$\7{sl}(2,\RR)\times \RR^2$.

\ProofC Assume to the contrary that $\7g=\7{sl}(2,\RR)\times\RR^{2}$
and denote by $\pi_{1}:\7l\to\7{sl}(2,\CC)$, $\pi_{2}:\7l\to\CC^{2}$
the canonical projections.  As in the proof of the previous claim, the
cases $\dim \7q\cap \rad(\7l)\ge 1$ can be ruled out immediately. We
therefore only have to exclude the case $\7q\cap\CC^{2}=0$: In that
case, let $\7b$ be the 2-dimensional $\pi_{1}$-image in
$\7{sl}(2,\CC)$. It follows $\sigma(\7b)\ne \7b$ since otherwise
$\pi_\7s(\7q+\sigma\7q)\subset \7b$ would contradict \23. Consequently
there exist $\4e\in \7b,$ $\4h\in \7b\cap \sigma\7b$ and
$\4f=\sigma(\4e)\in \sigma(\7b)=\pi(\sigma\7q)$ with the properties
described in \ruf{FN}. The fact that $\7q,\sigma\7q$ are subalgebras
and that $\pi(\7q)=\CC\4h \oplus \CC \4e$ determine $\7q,\sigma\7q$ as
follows: $\7q=\CC\cd(\4h,w)\oplus\CC\cd (\4e,0)\steil{and}
\sigma\7q=\CC\cd (\4h,-\overline w) \oplus \CC\cd (\4f,0)\steil{for a}
w\in \CC^2\,.$\nline We may further assume that $w,\overline w$ are
linearly independent, otherwise $\pi^{}_{2}(\7q+\sigma\7q)\ne\CC^{2}$
would contradict \23. As a consequence, $(\4h,0)\not\in
\7q+\sigma\7q$. The definition of $\7f$ then shows $
\7f=\CC\cd(\4h,w)\;$.  On the other hand, for this $\7f$ we have
$[\7f,\sigma\7q]\subset \7f\oplus \sigma\7q$ in contradiction to
\25. This proves the claim and Lemma \ruf{FM}.\qed

\medskip So far we know that under the assumption $\7s\cong
\7{sl}(2,\RR)$
necessarily $\7g\cong\7{sl}(2,\RR)\times\7r$, where $\7r$ is the
2-dimensional non-abelian Lie algebra. To determine the full CR-algebra
$(\7g,\7q)$ we still have to find out how the complex subalgebra
$\7q\subset\7l$ sits inside $\7{sl}(2,\CC)\times\7r^{\CC}$.

\Claim{FQ} Up to a CR-algebra automorphism of $(\7g,\7q)$ the
subalgebra $\7q\subset\7l$ is obtained in the following way: Fix a
linear basis $\4x,\4z$ of $\,\7r^{\CC}$ with $[\4x,\4z]=\4z$ and a
$\sigma$--adapted $\7{sl}_{2}$-triple $\4e,\4h,\4f\in\7{sl}(2,\CC)$
(see \ruf{FN} for the definition). Then
$\7q=\CC(\4h\,,2\4x+\mu\4z)\;\oplus\;\CC(\4e\;, \nu\4z)\;$ for
suitable $\mu,\nu\in\CC$ with $\Im(\mu)=\pm2$ and $|\nu|=1$.

\ProofC The case $\7q=\7r^{\CC}$ can clearly be
excluded. If the intersection $\7q\cap \7r^{\CC}$ would be
1-dimensional then $\7q\cap \7r^{\CC}\;\oplus\;\sigma\7q\cap
\7r^{\CC}=\7r^{\CC}$ (as the sum $\7q+\sigma\7q$ must be
direct). Since $\7f$ must also be 1-dimensional, we have
$\7f=\7r^{\CC}\cap \7q$. But then $[\7f,\sigma\7q]\subset
\7r^{\CC}=\7f\oplus \sigma\7f$, violating condition \25. \nline It
remains the case $\7q\cap \7r^{\CC}=0,$ i.e., $\pi_\7s(\7q)$ is a
Borel subalgebra $\CC \4e \oplus \CC \4h,$ and
$\pi(\sigma\7q)=\CC\4h\oplus \CC \4f$ where $\4h,\4e,\4f\in
\7{sl}(2,\CC)$ are chosen as explained in \ruf{FN}. Since
$\7q\subset \7s^{\CC}\times \7r^{\CC}$ is a Lie sub{\sl algebra}, it
necessarily has the following form:
$$ \diagram{\7q=\CC\cd (\4h,\lambda \4x+\mu \4z)\oplus \CC\cd (\4e,
\nu \4z)\Steil{with}\sigma\7q= \CC\cd (-\4h,\overline\lambda
\4x+\overline \mu \4z)\oplus \CC\cd (\4f, \overline \nu \4z) \cr
\hbox{and} \quad \lambda,\mu,\nu\in \CC\Steil{satisfying}
\lambda=2\steil{if}\nu\ne 0\,.}\Leqno{FR}
$$ {\fett$\nu =0$:} Then $\lambda \4x+\mu \4z$ and $\overline \lambda
\4x+\overline \mu \4z$ must be linearly independent in $\7r^{\CC}$,
i.e., $\lambda\overline\mu\ne \overline\lambda\mu$ (otherwise
$\pi_\7r(\7q{+}\sigma\7q)\ne \7r^{\CC}$). Observe that
$(\4h,0)\not\in\7q+\sigma\7q$ and $\lambda\ne0$. A direct verification
shows $\7f=\CC\cd(\4h,\lambda \4x+\mu \4z)$ if $\lambda\in i\RR$ and
$\7f=0$ otherwise. But then $\7f\oplus \sigma\7q$ is a subalgebra in
contradiction to \25.\nline {\fett $\nu\ne0$:} Possibly after
replacing $\4z$ by $|\nu|\4z$ we may assume $|\nu|=1$ in
\Ruf{FR}. Employing the definition of $\7f$ in \ruf{CA}.iv, a
simple calculation shows
$$
\7f=\cases{\CC(\4h\mp2i\overline\nu\,\4e\,,\,2\4x +\Re(\mu)\4z)\quad
& $\Im(\mu) =\pm 2$\cr 0 & otherwise,}\Leqno{FS}
$$
that is, $\7f$ is 1-dimensional only if $\Im\mu =\pm 2$. This proves
the claim.\qd

\Lemma{ZC} For the CR-algebra $(\7g,\7q)$ in \ruf{FQ} the associated
CR-germ $(M,o)$ is CR-equivalent to $(\5M,a)$, where $\5M$ is the tube
over the future light cone.

\Proof We start by giving a particular representing manifold for the
associated germ, compare Example 6.6 in \Lit{FEKA}. Let $\mu,\nu$ be
the constants occurring in \ruf{FQ} and consider the affine quadric
$Z:=\SL(2,\CC)\subset \CC^{2\times 2}$, on which the group $\widehat
L:=\SL(2,\CC)\times \SL(2,\CC)$ acts holomorphically by $z\mapsto
gzh^{-1}$ for all $(g,h)\in\widehat L$. Then $\widehat
G:=\SL(2,\RR)\times \SL(2,\RR)$ is a real form of $\widehat L$. Via
$$
\eqalign{ \4x:=&{1\over 2}\pmatrix{ 0 &1 \cr 1 & 0},\quad \4z:={1\over
2}\pmatrix{1 &\mkern-10mu -1 \cr 1 &\mkern-10mu -1}\quad
\in\;\7r\qquad\cr \4e:=&{\nu\over 2}\pmatrix{1 & i\cr i &
\mkern-12mu -1},\quad \4h:=\pmatrix{0 & \mkern-10mu -i \cr i &
0},\quad \4f:=\sigma(\4e)\;\;\;\in\;\;\7s^{\CC} }\Leqno{VB}
$$ we consider $\7l=\7s^{\CC}\times\7r^{\CC}$ as a complex subalgebra
of $\widehat\7l=\7{sl}(2,\CC)\times\7{sl}(2,\CC)$ and $\7g$ as a
subalgebra of $\widehat{\7g}$. \nline We recall some basic facts,
compare \Lit{FEKA}: The polynomial function $\psi(z):=\det(z+\overline
z)-2$ on $Z$ is invariant under the action of the group $\widehat
G$. Furthermore, the nonsingular part $M$ of the algebraic subset
$S:=\psi^{-1}\big(\{\pm2\}\big)$ is a (non-connected) hypersurface in
$Z$, locally CR-equivalent to $\5M$. On the other hand, the singular
part of $S$ is a totally real submanifold of $Z$. Consider the point
$$
o:={e^{-\pi i/4}\over4}\pmatrix{4+\mu & -\mu \cr i\mu
&i(4-\mu)}\;\in\;Z
$$
which actually is in $S$ because of $\psi(o)=\Im(\mu)$. A direct
computation shows that the isotropy subalgebra $\7l_{o}$ of
$\7l\subset\widehat\7l$ at $o\in Z$ is $\CC(\4h,2\4x+\mu\4z)\oplus \CC
(\4e,\nu\4z)$ and that the isotropy subalgebra of
$\7g\subset\widehat\7g$ at $o$ is trivial. This implies $a\in M$ and
also that the CR-algebra $(\7g,\7q)$ is associated to the germ
$(M,o)$. Since $M$ is locally CR-equivalent to the tube $\5M$ over the
future light cone the proof for \ruf{ZC} is complete.\qed

\KAP{Dreizehn}{The case $\7s\cong\su(2)$}

The status of our proof so far can be summarized as follows: For the
CR-algebra $(\7g,\7q)$ satisfying \ruf{FD} and $\7g^{\s}\ne 0$ the
simple ideal $\7s$ of $\7g^{\s}$ can only by isomorphic to
$\7{sl}(2,\RR)$ or $\su(2)$, Furthermore, the case
$\7s\cong\7{sl}(2,\RR)$ only occurs if the associated CR-germ is
equivalent to $(\5M,o)$ with $\5M$ being the light cone tube, and then
$\7g\cong \7{sl}(2,\RR)\times \7r$ with nonabelian 2-dimensional $\7r$.
Consequently, only the situation needs to be investigated when
$\7g^{\s}$ only contains simple ideals isomorphic to $\su(2)$. We
assume this throughout this section and state as main lemma:

\Lemma{VD} There is no simple factor of $\7g^{\s}$ isomorphic to
$\su(2)$.\Formend

\noindent The {\sl proof} will be subdivided into several claims.
Note that contrary to $\7{sl}(2,\RR),$ the only non-trivial proper
subalgebras of $\su(2)$ are 1-dimensional tori.  As before let $\7s$
be a fixed simple factor of $\7g^{\s}$ and denote by $\pi:\7g\to\7s$
the canonical projection. Then the image $\pi(\7g_o)$ must be a proper
subalgebra since otherwise $\7g'=\pi^{-1}(0)$ would violate
\26$_{2}$. Our first observation is
\Claim{} For $\7g$ and $\7r:=\rad(\7g)$ only the following cases may
occur\vskip-23pt
$$
\displaylines{ \7g\cong\su(2)\ltimes\7r,\qquad \7g\cong\su(2)\times
\su(2) \qquad\hbox{\rm or} \cr \7g\cong\big(\su(2)\times
\su(2)\big)\ltimes\7r \Steil{ with }\7r\ne
0\steil{and}\7g_o=(\RR\kern1pt \4t_1 \times \RR\kern1pt \4t_2)
\ltimes(\7g_o\!\cap\7r)\steil{for}\4t_{j}\in\su(2)\;.}
$$

\ProofC Write $\7g=(\7s_1\oplus \7s_2\oplus \7s'')\ltimes \7r$ where
$\7s_1,\7s_2\cong\su(2)$, $\7s''$ is the complementary ideal in
$\7g^{\s}$ and $\pi_j$ denotes the projection onto $\7s_1$ or $\7s_2$,
respectively.  Only one of the following possibilities could occur:
\nline $\Bullet\pi_1(\7g_o)=0$. Then $\7s_1\oplus \7s_2$ is already a
transitive subalgebra of $\7g$ and $\7s''=\7r=0$. \nline
$\Bullet\pi_1(\7g_o)\ne 0\ne \pi_2(\7g_o)$. Then either
$\7g=\7s_1\oplus\7s_2$ or $\7g_o\cap (\7s_1\oplus \7s_2)=
(\7g_o\cap\7s_1)\oplus (\7g_o\cap \7s_2)= \RR \4t_1 \oplus \RR \4t_2$
for suitable nonzero $\4t_1,\7t_2\in \su(2)$. In the latter case
$\7g_o\cap (\7s''\ltimes \7r)$ is of codimension 1 in $\7s''\ltimes
\7r$. But then we conclude that $\7s''=0$ as there is no
1-codimensional subalgebras in $\su(2)$. \qd

\noindent For the proof of \ruf{VD} we only need to investigate the
above three types of $\7g$. We repeatedly use the fact that each 1- or
2-dimensional representation of $\su(2)\cong\so(3)$ is trivial and
that each toral subalgebra $\7t\subset \su(2)$ acts irreducibly on
$\qu{\su(2)}{\7t}.$

\Claim{} $\7g\cong(\su(2)\times \su(2))\ltimes\7r$ implies
$\7r=0$.

\ProofC Let $\pi_1,\pi_2$ be the projection onto the first and the
second simple factor, respectively. As observed in the preceeding
Claim, $\7g_o=(\7t_1\times \7t_2)\ltimes \7r_o$ with $\7r_o:=\7g_o\cap
\7r$ and 1-dimensional toral subalgebras $\7t_j.$ Recall that we have
the $\7g_o$-stable filtration $\7g\supset \7H\supset \7F\supset
\7g_o.$ At least one of the images $\pi_j(\7F)$ coincides with
$\7s_j,$ say, for $j=2.$ Since $\7t_j$ acts irreducibly on
$\qu{\7s_j}{\7t_j},$ it follows $\7F=\7s_2\ltimes \7r_o$. Since then
$\7H\cap \7s$ is at least 2-dimensional, the irreducibility of the
action of $\7t_1$ implies $\7H=\7s \times \7s_2\ltimes \7r_o.$ But
this would imply $[\7H,\7H]\subset \7H,$ contradicting \23.\qd

\Claim{} $\7g\not\cong\su(2)\times\su(2)$.

\noindent\rm {\sl Proof of the Claim.} Suppose to the contrary that
$\7g=\su(2)\times\su(2)$. Since $\7F\subset\7g$ is a subalgebra of
dimension 3 and there is no solvable subalgebra of this dimension,
necessarily $\7F\cong \su(2)$. Consequently, either $\7F$ is one of
the simple factors of $\7g$ or $\7F$ is the graph of an automorphism
of $\su(2)$. Both possibilities lead to a contradiction: In the first
case $\7F$ is an ideal in $\7g,$ contradicting \25, and in the second
case $\7F$ acts irreducibly on $\qu{\7g}{\7F}$, violating the
existence of the $\ad(\7F)$-stable proper subspace $
\qu{\7H}{\7F}\subset\qu{\7g}{\7F}$.\qd

\noindent The remaining case $\7g=\su(2)\ltimes\7r$ with
$\7r:=\rad(\7g)$ is the most involved. In this situation the
projection of $\7g_o$ under the canonical projection
$\pi:\7g\to\su(2)$ is of dimension $\le1$. As usual we denote the
canonical projection $\7l\to\7{sl}(2,\CC)$ by the same symbol
$\pi$. We investigate the various possibilities for the $\pi$-images
of the Lie subalgebras $\7l_o\subset\7f\subset\7q$ defined in
\Ruf{FA}.\nline $\Bullet\pi(\7l_o)=\pi(\7f)=\pi(\7q)$.  Since
$\pi(\7l_o)$ is a $\sigma$-stable torus, it would follow
$\pi(\7q+\sigma\7q)=\pi(\7l_o).$ Counting dimensions, this cannot
happen. \nline $\Bullet\pi(\7l_o)\subsetneq\pi(\7f)=\pi(\7q)$. Here,
we have to rule out the 2 possibilities $\dim \pi(\7l_o)\le1$. \nline
If $\dim\pi(\7l_o)=0$ then $\pi(\7F)$ is a 1-dimensional toral
subalgebra in $\su(2),$ and consequently
$\pi(\7f)=\pi(\sigma\7f)=\pi(\7q)=\pi(\sigma\7q)$ is
1-dimensional. This contradicts the fact that $\7q+\sigma\7q$ is a
hyperplane in $\7l.$ It remains the case when $\pi(\7l_o)$ is
1-dimensional. Then $\pi(\7F)=\7s,$ and possibly after replacing
$\7s\cong \su(2)$ by the Levi factor contained in $\7F$ we may assume
that $\7s\subset \7F.$ By dimension reasons $\7r_o:=\7g_o\cap
\7r=\7F\cap \7r$. Further, $\7r_o$ is an ideal in $\7F$ and we have
$\7F=\7s\ltimes \7r_o$ as well as $\7H=\7s\ltimes \7r_\7H$ with
$\7r_\7H:=\7H\cap \7r.$ Since $\dim \7r_\7H-\dim \7r_o=2,$ the
ad-representation of $\7s$ on $\qu{\7r_\7H}{\7r_o}$ is trivial, i.e.,
$[\7s,\7r_\7H]\subset \7r_o.$ Since $\7f=\7l_o\oplus \CC \4y$ for some
$\4y\in \7s^{\CC}\cong\7{sl}(2,\CC)$ and $\sigma\7q=\sigma\7f\oplus
\CC\overline{\4x}$ for some $\overline{\4x}\in \7r^{\CC}_\7H,$ we
would have $[\7f,\sigma\7q]=[\7l_o\oplus \CC \4y,\sigma\7f\oplus \CC
\overline{\4x}]\subset \7f+\sigma\7f,$ in contradiction to \25. \nline
$\Bullet\pi(\7l_o)=\pi(\7f)\subsetneq\pi(\7q).$ The possibility
$0=\pi(\7f)$ can be excluded since otherwise $\7F=\7r$ and
$\7f+\sigma\7f$ would be an ideal in $\7l,$ contradicting condition
\25. Next we deal with the case when
$\pi(\7l_o)=\pi(\7f)=\pi(\sigma\7f)$ is a (1-dimensional) toral
subalgebra. By assumption $\7b:=\pi(\7q)$ is then 2-dimensional, i.e.,
a Borel subalgebra, which implies that $\7b+\sigma{\7b}=\7s^{\CC}\cong
\7{sl}(2,\CC),$ or equivalently $\pi(\7H)=\7s.$ Counting dimensions,
$\7r_\7F:=\7F\cap \7r=\7H\cap \7r$. Further, $[\7r_\7F,\7H]\subset
\7r\cap \7H=\7r_\7F$ and since, due to condition \23, $\7H$ generates
$\7g$ as a Lie algebra, we deduce that $\7r_\7F$ is an ideal in $\7g.$
From $\pi(\7l_o)=\pi(\7f)$ follows that $\7f=\7l_o\oplus \CC\4r$ with
$\4r\in \7r^{\CC}_\7F.$ But this implies $[\7f,\sigma\7q]=
[\7l_0\oplus \CC\4r,\sigma \7q]\subset \sigma\7q+\7r^{\CC}_\7F\subset
\sigma\7q\oplus \7f,$ violating \25. It remains one last possibility
for the flag $\7l_o\subset \7f\subset \7q$ in $\7l$:

\Claim{} If $\pi(\7l_o)\subsetneq\pi(\7f)\subsetneq\pi(\7q)$ then
$\7l_o=0$ and $\7g\cong \su(2)\times\7r$ with $\7r:=\rad(\7g)$ of
dimension 2.

\ProofC The properness of the inclusions implies $\7l_o\cap
\7r^{\CC}=\7f\cap\7r^{\CC}=\7q\cap \7r^{\CC}=\sigma\7q\cap \7r^{\CC}$.
Since $\7q$ and $\sigma \7q$ generate $\7l$ as a Lie algebra, it
follows that $\7l_o\cap \7r^{\CC}$ is an ideal in $\7l,$ or
equivalently $\7g_o\cap \7r\vartriangleleft \7g.$ By the effectivity
assumption \26$_{1}$, we have $\7g_o\cap \7r=0.$ Consequently, since
$\dim\pi(\7l_o)\le1$ the same estimate holds for $\dim\7l_o$. Next we
show that the case $1=\dim_{\RR}\pi(\7g_o)\,$ ($\,=\dim
\pi(\7l_o)=\dim \7l_o$) cannot happen: Assume on the contrary that
$\pi(\7g_o)=1$. Then $\pi(\7F)=\7s$ and, possibly after replacing the
Levi factor $\7s\subset \7g$ by a conjugate one, we may assume
$\7s\subset \7F.$ Counting dimensions yields then $\7F=\7s\cong
\su(2).$ Since $[\7F,\7H]\subset \7H$ by \22, we have a representation
of $\7s$ on $\qu{\7H}{\7F}$. But this yields the contradiction as
$\dim \qu{\7H}{\7F}=2$, compare \21, implies that this representation
is trivial, i.e. $[\7s,\7H]=[\7F,\7H]\subset \7F$, violating
\25. \nline We have proved that the only possibility for $\7g$ is
$\su(2)\ltimes\7r$ with a 2-dimensional radical $\7r.$ Since $\su(2)$
can act only trivially on such an $\7r$, the above semidirect product
is in fact direct. This proves the claim.\qd

\Claim{} The case $\7g\cong\su(2)\times \7r$, $\dim \7r=2$, also
cannot occur.

\ProofC We write $\pi=\pi_\7s$ and $\pi_\7r$ for the projections onto
$\7s^{\CC}\cong \7{sl}(2,\CC)$ and $\7r^{\CC},$ respectively. The
proof analyses various possibilities for $\7f\subset \7q.$ \nline If
$\pi_\7s(\7F)=0$ then $\7F=\7r$ would be an ideal in $\7g$
contradicting \25. It remains the case when $\pi_\7s(\7F)$ as well as
$\pi_\7s(\7f)$ is 1-dimensional. In that situation necessarily
$\pi_\7s(\7q)\ne \pi_\7s (\7f),$ and $\pi(\7q)$ is a Borel subalgebra
in $\7{sl}(2,\CC).$ Further,
$\pi_\7s(\7q)+\sigma(\pi_\7s(\7b))=\7{sl}(2,\CC)$ and
$\pi_\7s(\7q)\cap\pi_\7s(\sigma\7b)$ is a Cartan subalgebra of
$\7{sl}(2,\CC).$ Similar to the situation considered in \ruf{FN} we
also can choose here a {\sl standard triple} $\4e,\4f,\4h\in
\7s^{\CC}$ with $\pi_\7s(\7q)=\CC \4h \oplus \CC \4e $ and
$\sigma(\4h)=-\4h$, but now $\sigma(\4e)=-\4f$ for the Cartan
involution $\sigma$. The radical $\7r$ cannot be abelian: Otherwise,
exactly as in the proof of Claim \ruf{FP}, we obtain a
contradiction. In the remaining non-abelian case we proceed as in the
proof of Claim \ruf{FQ} by investigating the various positions of
$\7q$ in $\7l$: The cases $\dim \7r^{\CC}\cap \7q>0$ are easily ruled
out (same argument as in the proof of \ruf{FQ}, following
\Ruf{FR}). Hence, we may assume that $\7q$ and $\sigma\7q$ are given
by the formula \Ruf{FR}, except that now $\sigma\7q =\CC\cd
(-\4h,\overline\lambda \4x+\overline \mu \4z)\oplus \CC\cd (-\4f,
\overline \nu \4z)$, i.e., the sign in front of $\4f$ has
changed. This slight difference is precisely the reason why in case
$\nu\ne 0,$ contrary to \Ruf{FS}, the CR-germ $(M,o)$ associated to
$(\7g,\7q)$ would be Levi nondegenerate, as shown by a simple
computation. This contradicts our fundamental assumption and concludes
the proof of the claim as well as the proof of Proposition
\ruf{VD}.\qed

\KAP{Vierzehn}{Reduction to the case where $\7g$ is solvable and of
dimension 5}

Striking the balance for the proof of Theorem II obtained so far,
we have shown the following:

\noindent{\sl Let $(M,o)$ be an arbitrary
locally homogeneous $2$-nondegenerate CR-germ of dimension $5$ and let
$(\7g,\7q)$ be an associated CR-algebra. If the Lie algebra $\7g$ is
not solvable then $M$ is locally
CR-equivalent at $o\in M$ to the tube $\5M$ over the future light cone.}

For the rest of the proof we therefore assume that every CR-algebra
$(\7g,\7q)$ under consideration satisfies the fundamental assumption
\ruf{FD} and that $\7g$ is solvable. As main result of this section we
show

\Lemma{FT} The solvable Lie algebra $\7g$ has dimension $5$, i.e.,
$\7g_o=0$.\Formend

\noindent The {\sl proof} of this Lemma will be subdivided into several
steps. Recall that by definition the {\sl nilradical}
$\7g^{\nil}$ of $\7g$ is the maximal nilpotent ideal in $\7g$. It is
well-known that $\7g^{\nil}$ contains the commutator subalgebra
$[\7g,\7g]$, and each element $\xi\in\7g^{\nil}$ is ad-nilpotent in
$\7g.$ Similarly, we denote the (complex) nilradical of $\7l$ by
$\7l^{\nil}.$ We retain the notation from \ruf{FA} and \ruf{CA}.

\Claim{FU} ${\7g}_o\subset \7g^{\nil}$.

\ProofC Since $[\7g,\7g]\subset \7g^{\nil},$ the quotient
$\qu{\7g}{\7g^{\nil}}$ is abelian. Let $\pi:\7g\to
\qu{\7g}{\7g^{\nil}}$ be the projection. Assume that $\7g_o\not\subset
\7g^{\nil}$, i.e., $\pi(\7g_o)\ne 0.$ Select a subalgebra $\7r\subset
\qu{\7g}{\7g^{\nil}}$ (possibly 0) which is a complement of
$\pi(\7g_o)$. But then $\7g':=\pi^{-1}(\7r)$ would violate \26$_{2}$.
Consequently we necessarily have $\pi(\7g_o)=0$, i.e., $\7g_o\subset
\7g^{\nil}$ as claimed. \qd
\Parag{XY} Recall that for the CR-algebra $(\7g,\7q)$ under
consideration there exist flags $\7l_o\subset\7f\subset\7q$ and
$\7l_o\subset\sigma\7f\subset\sigma\7q$ of subalgebras in $\7l$, as
defined in \ruf{FA}. For the subsequent considerations we select the
following elements in $\7l$: $\4y\in\7f\setminus \7l_o$ and $\4x\in
\7q\setminus \7f,$ and write $\overline{\4x}:=\sigma(\4x)$ and
$\overline{\4y}:=\sigma(\4y),$ for short. Then

\noindent \hfil$\7f= \CC\4y \oplus \7l_o,\quad \7q=\CC\4x\oplus
\CC\4y\oplus \7l_o, \quad \sigma\7f= \CC\overline{\4y} \oplus
\7l_o\Steil{and} \sigma\7q=\CC\overline{\4x}\oplus
\CC\overline{\4y}\oplus \7l_o\,.$

\noindent The inclusion $\7l_o\subset \7l^{\nil}$ is guaranteed by
\ruf{FU}. Hence, $\7l_o$ acts by ad-nilpotent endomorphisms on
$\7l.$ In particular, $[\7l_o,\7f]\subset\7l_o$ and
$[\7l_o,\7q]\subset \7f.$ The condition \23 means
that $[\4x,\overline{\4x}]\not\in \7q +\sigma\7q$, and
\25 is equivalent to $[\overline{\4x},\4y]\not\in
\7f+\sigma\7q.$ \nline Let $\7z\subset \7l^{\nil}$ be the center of
the nilradical (which is nontrivial if $\7l \ne 0$), to which we refer
as to the {\sl nilcenter} of $\7l$. As for every characteristic
ideal, we have $\sigma(\7l^{\nil})=\7l^{\nil}$ and $\sigma(\7z)=\7z.$

\Claim{FV} {\rm (i)} The nilcenter $\7z$ of $\7l$ is not contained in
$\7q+\sigma\7q.$ \ite=2 \1 $\7l_o=0$ if $\dim \7z\ge 2$.

\Proof {\fett ad (i):} Assume that (i) is not true and let
$\4x,\4y,\overline{\4x},\overline{\4y}\in \7l$ be as in \ruf{XY}.  Let
$\zeta:=a\4x+\overline a\overline{\4x}+b\4y+\overline
b\overline{\4y}+\gamma_o,$ $\gamma_o\in\7g_o,$ $a,b\in \CC$, be an
arbitrary element in $\7z^\sigma.$ Then
$[\overline{\4x},\zeta]\in\7z\subset \7q+\sigma\7q$ since $\7z$ is an
ideal in $\7l$.  On the other hand
$$
[\overline{\4x},a\4x+\overline a\overline{\4x}+b\4y+\overline
b\overline{\4y}+\gamma_o
]=a[\overline{\4x},\4x] + [\overline{\4x},b\4y+\overline
b\overline{\4y}+\gamma_o
]\equiv a[\overline{\4x},\4x]\equiv0\;\mod\; \7q{+}\sigma\7q\,,
$$
which implies $a=0.$ This shows that $\7z\subset \FFF.$ Given then
$\zeta=b\4y+\overline b\overline{\4y}+\gamma_o\in \7z^\sigma$, the
inclusion $[\overline{\4x},\zeta]\subset\7z$ holds since $\7z$ is an
ideal. On the other hand
$$
[\overline{\4x},b\4y+\overline
b\overline{\4y}+\gamma_o\;
]=b[\overline{\4x},\4y]+[\overline{\4x},\overline
b\overline{\4y}+\gamma_o]\equiv b[\overline{\4x},\4y]\equiv0\;\mod \;
\7f{+}\sigma\7q\;.
$$ Hence, $b=0$ as the consequence of the above equation and \25. But
this cannot be true since then $\7z\subset \7l_o$ would be a
nontrivial ideal of $\7l$, violating \26$_{1}$. \nline {\fett ad
(ii):} Assume $\dim \7z\ge 2,$ i.e. $\7z\cap \QQQ \ne 0$. Recall that
$\7l_o\subset \7l^{\nil}$ by \ruf{FU} and consequently
$[\4y,\7l_o]\subset\7l_o,$ $[\overline{\4y},\7l_o]\subset\7l_o.$ Let
an arbitrary $\zeta :=a\4x+\overline a\overline{\4x}+b\4y+\overline
b\overline{\4y}+\gamma_o \in \7z^\sigma\cap \QQQ$ be given. Computing
its bracket with $\7l_o$ we obtain:
$$
[\zeta,\7l_o]=a[\4x,\7l_o]+\overline a[\overline{\4x},\7l_o]\equiv0
\;\mod\;\7l_o \;.
$$ Either $a=0$ for all such $\zeta,$ and then $\7z\cap \QQQ=\7z\cap
\FFF$, or $a\ne 0,$ and then
$[\4x,\7l_o]\subset\7l_o\supset[\overline{\4x},\7l_o].$ In the latter
case it follows that $\7l_o$ is an ideal in $\7l$ (due to
\ruf{FC}.{\ninecmss III}, $\4x,\4y,\overline{\4x},\overline{\4y}$ and
$\7l_o$ generate $\7l$ as a Lie algebra), hence, $\7l_o=0$ as claimed.
The other possibility would be $\7z\cap (\7q+\sigma\7q)=\7z\cap
(\7f+\sigma\7f)$ and we show that this cannot happen: Given an
arbitrary $\zeta=b\4y+\overline b\overline{\4y}+\gamma_o\in
\7z^\sigma\cap (\7f+\sigma\7f)$, note that
$$
[\overline{\4x},\zeta]\in \7z\cap \QQQ= \7z\cap\FFF.
$$
More explicitly, $[\overline{\4x}\,,\,b\4y+\overline
b\overline{\4y}+\gamma_o ]=b[\overline{\4x},\4y]+[\overline{\4x}\,,\,
\overline b\overline{\4y}+\gamma_o]\in
b[\overline{\4x},\4y]+\sigma\7q$. But then the above equation together
with \25 would imply $b=0,$ i.e., $\7z\cap \QQQ=\7z\cap\7l_o\ne 0.$
This is absurd, as
$$
[\7q+\sigma\7q,\7z\cap\7l_o]\subset \7z\cap
(\7q{+}\sigma\7q)=\7z\cap\7l_o \ne 0,
$$
i.e., since $\7q+\sigma\7q$ generates $\7l$ as a Lie algebra, $\7z\cap
\7l_o$ would be a nontrivial ideal in $\7l$. \qd

\noindent
In remains the case when the nilcenter is 1-dimensional.

\Claim{FW}Suppose $\dim\7z=1.$ Then
$\7l=\7z\oplus \QQQ\,$ and $\,\7g_o=\7l_o=0.$

\Proof Since $\7z\not\subset \7q+\sigma\7q$ by \ruf{FV}, the sum $\7z
+ \QQQ $ is direct. Recall that the recursively defined subspaces
$C_0(\7l^{\nil}):=0,$ $C_k:=C_k(\7l^{ \nil}):=\{u\in \7l^{
\nil}:[u,\7l^{ \nil}]\subset C_{k-1}(\7l^{ \nil})\}$ for every $k>0$
form the ascending central series of $\7l^{\nil}.$ Clearly, $\7z=C_1$
and $\sigma(C_k)=C_k$ for all $k.$ Either $\7z=C_1=C_2=\7l^{\nil}$ and
consequently $\7l_o=0$ (due to \26$_{1}$ $\7l_o$ must be a proper
subalgebra of the 1-dimensional algebra $\7l^{ \nil}$) or $C_1\ne
C_2.$ In the latter case $C_2\cap \QQQ\ne 0.$ Let $\eta=a\4x+\overline
a \overline{\4x}+b\4y+\overline b\overline{\4y}+\gamma_o\in
C_2^\sigma\cap \QQQ$ be arbitrary. Since $[\eta,\7l_o]\subset \7z\cap
\QQQ=0,$ we have
$$
[\eta,\7l_o]\equiv a[\4x,\7l_o]+\overline
a[\overline{\4x},\7l_o]\equiv0\;\mod\;\7l_{o}\,.
$$
If $a\ne 0$ then $[\4x,\7l_o],[\overline{\4x},\7l_o] \subset \7l_o$
and $\7l_o$ is an ideal in $\7l,$ i.e., $\7l_o=0$ by \26$_{1}$. If
$a=0$ for every choice of $\eta\in C_2^\sigma\cap \QQQ$ as above, then
$C_2\cap \QQQ=C_2\cap \FFF.$ This possibility can be ruled out as
follows: For a nonzero $\eta=b\4y+\overline
b\overline{\4y}+\gamma_o\in C_2^\sigma\cap\FFF$ we have $[\eta,\4x]\in
\QQQ \cap C_2=\FFF\cap C_2,$ i.e.,
$$
[b\4y+\overline b\overline{\4y}+\gamma_o,\4x]=\overline
b[\overline{\4y},\4x] \;+
([b\4y+\gamma_o,\4x])\equiv0\;\mod\;\FFF\,,
$$
which is only possible if $\overline b=0.$ This would imply $\QQQ\cap
C_2=\7l_o\cap C_2.$ On the other hand, the identity $\QQQ\cap\
C_2=\7l_o\cap C_2$ implies that $\7l_o\cap C_2$ is an ideal of $\7l.$
The effectivity of $(\7g,\7q)$ forces then $\7l_o\cap C_2=0,$
contradicting $\QQQ\cap C_2\ne 0.$ This completes the proof of \ruf{FW}
and, together with \ruf{FU}, \ruf{FV} also the proof
of Lemma \ruf{FT}. \qed

\KAP{Fuenfzehn}{The existence of a 3-dimensional abelian ideal in
$\7g$ suffices}

The proof of Theorem II has brought us to the point where we may and
do henceforth assume that the Lie algebra $\7g$ in the CR-algebra
$(\7g,\7q)$ satisfying \ruf{FD} is solvable. The
subalgebra $\7q\subset\7l=\7g\oplus i\7g$ then necessarily is of
complex dimension $2$.

The main result \ruf{GS} of this section states that the CR-germ
associated to $(\7g,\7q)$ is represented by the tube $F+i\RR^3$ over
an affinely homogeneous surface $F\subset \RR^3$ if $\7g$ is
isomorphic to a semidirect product $\7h\ltimes \7r$ with $\7h$ being a
2-dimensional Lie subalgebra and $\7r\cong \RR^3$ being an abelian
ideal.  Once this Main Lemma is proved, our proof of the
classification theorem will be complete as soon as we can show that
every 5-dimensional solvable Lie algebra $\7g$ occurring in $(\7g,\7q)$
indeed is isomorphic to a semidirect product as above. This will be
achieved in the final section \ruf{Sechzehn}.  In this section we only
prove the partial result that if $\7q$ is abelian, then also the
commutator $[\7g,\7g]$ is abelian and 3-dimensional. Moreover, there
exists an abelian subalgebra $\7h\subset \7g$ such that
$\7g=\7h\ltimes [\7g,\7g]$.

Since there is no general structure theory for solvable Lie algebras,
we develop ad hoc methods and describe the structure constants in
$\7l=\7g^{\CC}$ with respect to a particularly chosen basis. Every
CR-algebra $(\7g,\7q)$ under consideration gives rise to the
1-dimensional subalgebras $\7f$ and $\sigma\7f$ of the 2-dimensional
subalgebras $\7q$ and $\sigma\7q,$ respectively. We construct a basis
of $\7l$ which reflects the conditions \21-\25 and investigate the
various possibilities for the values of the corresponding structure
constants. Select a non-zero $\4z\in \7z^{-\sigma}\setminus \QQQ$
(this is possible due to Claim \ruf{FW}) and an $\4x\in \7q\setminus
\7f$ such that for $\overline{\4x}:=\sigma\4x$ the congruence
$[\4x,\overline{\4x}]\equiv \4z\; \mod\;\7q +\sigma\7q$ holds (this is
possible due to \23). By \25 it is further possible to select $\4y\in
\7f\setminus \7l_o$ such that for $\overline{\4y}:=\sigma\4y$ the
structure equations of $\7l$ are of the following form (in particular,
the coefficient in front of $\overline{\4x}$ in the second equation is
$b_2=1$):
$$
\diagram{[\4x,\overline\4x]&\;=\;&\4z+\;& a_1\4x &\;- \;&\overline
a_1\overline \4x&\;+\;& a_2\4y &\;-\;&\overline a_2\overline \4y \cr
[\4x,\overline{\4y}]&\;=\;& & b_1\4x&+& \phantom{\;1\,}\overline\4x
&+&b_3\4y&+&\hfill b_4 \overline{\4y} \cr [\4y,\overline{\4y}]&\;=\;&
& & & & &\hfill c\,\4y & - &\hfill\overline c\,\overline{\4y} \cr
[\4y,\4x]&\;=\;& & d_1\4x & & &+ & d_2\4y & & \cr
[\4z,\eta]&\in&\7z(\7l^{\nil}) \rlap{ \qquad for every $\eta\in \7l\;.$}
}\Leqno{FI}
$$
The brackets $ [\4y,\overline{\4x}]$ and
$[\overline{\4y},\overline{\4x}]$ are completely determined by the
above 5 equations due to the fact that $\sigma:\7l\to \7l$ is an
antilinear Lie algebra automorphism. Of course, not for all values of
the constants the above identities give rise to a Lie {\sl
algebra}. In fact the above structure constants $a_1,...\,,d_2,$ are
subject to further constraints, imposed by (1) the Jacobi identity,
(2) our assumption that $\7l$ is solvable and, (3) our assumption
$\4z\in \7z(\7l^{\nil}).$ \nline Conversely, let a 5-dimensional
solvable complex Lie algebra $\7l=\CC\4z\oplus \CC \4x\oplus\CC
\overline{\4x}\oplus \CC\4y\oplus \CC\overline{\4y}$ be given with
structure equations as in \Ruf{FI} (together with
$[\overline{\4x},\4y]= \4x +\overline b_1\overline\4x + \overline
b_4\4y +\overline b_3\overline{\4y} $ and
$[\overline{\4y},\overline{\4x}]=\overline d_1\overline\4x + \overline
d_2\overline{\4y}\,$) for certain $a_1,..\,,d_2\in\CC$ and
$\4z\in \7z^{-\sigma},$ where $\7z$ is the nilcenter of $\7l$. Define
$$
\7g:=\RR \,i\4z \oplus \RR(\4x+\overline{\4x})\oplus
\RR(i\4x-i\overline{\4x})\oplus \RR(\4y+\overline{\4y})\oplus
\RR(i\4y-i\overline{\4y})\Steil{ and }\7q:=\CC\4x\oplus \CC\4y \;\;.
$$
Then the CR-algebra $(\7g,\7q)$ satisfies the fundamental assumption
\ruf{FD}.

\medskip As already mentioned, besides the geometrically motivated
conditions \ruf{FC}, which already are incorporated in \Ruf{FI},
further conditions affect the particular values of the structure
constants, for example those given by the Jacobi identity. We use
$\Jac{\xi_1}{\xi_2}{\xi_3}:=[[\xi_1,\xi_2],\xi_3]+
[[\xi_2,\xi_3],\xi_1]+[[\xi_3,\xi_1],\xi_2]$ as shorthand. The
identity $\Jac{\4x}{\4y}{\overline{\4y}}=0$ implies
$$
|c|=1, \qquad d_1=\overline c -\overline b_1 \qquad \hbox{and}\quad
\overline c\re2b_1\in\RR\,.\Leqno{FY}
$$

\Remark{VJ} From $c\ne 0$, i.e., $[\4y,\overline{\4y}]\ne0$,
follows that the solvable subalgebra $\7f+\sigma\7f,$ and in turn
$\7l,$ cannot be nilpotent. \Formend

\medskip\noindent Keeping in mind the identities \Ruf{FY} it is
possible to readjust the basis $\4x,...\,,\4z$ of $\7l$ as follows:
Write $c=e^{2}$ for the coefficient $c$ in the third equation of
\Ruf{FI} and replace $\4x$ by $e\4x$ as well as $\4y$ by $c\4y.$ After
this replacement the structure equations \Ruf{FI} keep their form,
only $c$ changes to $c=1$ and $b_{1}$ becomes real.

\noindent {\sl Notational agreement:} For the rest of this section we
use $\alpha,\beta,\gamma$ and $\delta$ to denote real numbers, while
$a_j, b_j,c_j$ and $d_j$ stand for complex numbers. In particular, we
write $\beta_1:=b_1$ to underline that the structure
constant $b_1$ is real.

\noindent The Jacobi identity $\Jac{\4x}{\4y}{\overline{\4y}{}}=0$
implies further relations between the coefficients in the structure
equations \Ruf{FI}:\vskip-25pt
$$
\diagram{c=1,\qquad\qquad \beta_1:=b_1\in \RR \qquad\cr
d_1=1-\beta_1\;,\qquad \overline b_3-d_2=b_4(\beta_1-1)\;,\qquad
-\beta_1(b_3+d_2)=b_4-\overline b_4 \;.}\Leqno{FZ}
$$ \Joker{GQ}{Perfect basis.} Summarizing, there exists a basis
$\4z,\4x,\overline{\4x},\4y,\overline{\4y}$
($\overline{\4x}:=\sigma(\4x), \;\overline{\4y}:=\sigma(\4y)$) of
$\7l=\7g^{\CC}$ with $\4z\in \7z^{-\sigma}\setminus
(\7q{\oplus}\sigma\7q)$ and $\7q=\CC\4x\oplus \CC\4y$, $\7f=\CC \4y$
such that the corresponding structure constants in \Ruf{FI} satisfy
\Ruf{FZ}. We call each such basis {\sl perfect}.

\def\FFF{(\7f\mkern-1mu\oplus\mkern-1mu\sigma\7f)}
\def\QQQ{(\7q\mkern-1mu\oplus\mkern-1mu\sigma\7q)}

There are still more constraints for the structure
constants given by the Jacobi identity for further triples of elements
and also by the fact that $\7z$ is in the center of
$\7l^{\nil}$. Later on, we characterize these additional conditions
more explicitly. For now, we elaborate the particular structure of the
nilcenter:

\Lemma{GR} The nilcenter $\7z\subset \7g$ has dimension 1 or 3.

\Proof We closely analyze the conditions in \ruf{FC} in order to get
the dimension estimates. It is clear that $1\le \dim \7z\le 4$, see
Remark \ruf{VJ}. Assume $\dim \7z=4$. Then $\dim\7z\cap \QQQ=3$
follows by Lemma \ruf{FV}. But since $\7f+\sigma\7f$ is not abelian,
i.e., $\dim \7z\cap \FFF=1,$ there exist elements in $\7z$ of the form
$\4x+\eta_1,\overline{\4x}+\eta_2$ with $\eta_j\in \7f+\sigma\7f.$
This leads to a contradiction: Indeed, on the one side
$[\4x+\eta_1,\overline{\4x}+\eta_2]=0$, and on the other side \25
implies $[\4x+\eta_1,\overline{\4x}+\eta_2]
\equiv[\4x,\overline{\4x}]\not\equiv0\;\mod\;{\7q +\sigma\7q}$. Hence,
we have proved that $1\le \dim \7z\le 3.$

\noindent
The main difficulty is to rule out the possibility $\dim\7z=2$.
Select a perfect basis in $\7l$ as described in
\ruf{GQ} keeping in mind \Ruf{FI} and \Ruf{FZ}. We have to deal with 2
subcases.

\Claim{VF} If $[\4x,\4y]\ne 0$, i.e., $\7q$ is not abelian, then
$\dim \7z=1$.

\ProofC We first show that $\7z\cap \QQQ=\7z\cap \FFF$. Select
an arbitrary
$\4z':=\lambda\4x+\overline\lambda\overline{\4x}+\mu\4y+\overline\mu
\overline{\4y} \in \7z^\sigma\cap \QQQ.$ We have to investigate the 2
possibilities `$\beta_1\ne 1$' and `$\beta_1=1$': In the first
case we get

\noindent \hfil
$[\4z',(1-\beta_1)\4x+ d_2\4y]\equiv\overline\lambda(1-\beta_1)
[\overline{\4x},\4x]\equiv0\;\mod\;\7q +\sigma\7q$

\noindent since $[\4y,\4x]\in \7l^{\nil}$. Therefore$\lambda=0$ by the
condition \23, i.e., $\4z'\in \7f +\sigma\7f.$ If $\beta_1=1,$ i.e,
$d_1=0$, then by our assumption $d_2\ne 0,$ and in turn $\4y\in
\7l^{\nil}.$ Hence,
$[\4z',\4y]\equiv\overline\lambda[\overline{\4x},\4y]\equiv0\;\mod\;
\7f +\sigma\7f,$ which, together with \25 also forces
$\lambda=0$. This proves $\7z\cap \QQQ=\7z\cap \FFF$. \nline We claim
that this identity can only h old if both sides vanish, i.e., $\dim
\7z=1$ by Lemma \ruf{FV}.i: Assuming to the contrary that $\7z\cap
(\7q\oplus\sigma\7q)\ne 0$, then on the one hand there exists
$\4z'=\mu\4y+\overline\mu\overline{\4y} \in \7z\cap \FFF$ with $\mu\ne
0.$ On the other hand,
$[\4x,\4z']=[\4x,\mu\4y+\overline\mu\overline{\4y}]\in
\7z\cap\QQQ=\7z\cap \FFF$ which in view of \25 is only possible if
$\mu=0$. This shows $\dim\7z=1$.

\noindent It remains to rule out the second subcase:

\Claim{VG} If $[\4y,\4x]=0$ , i.e., $\7q$ is abelian, then the
commutator $[\7l,\7l]$ is a 3-dimensional abelian ideal. Further, if
$\dim \7z\ge 2$ then $[\7l,\7l]=\7l^{\nil}$ and consequently
$\7z=[\7l,\7l]$ is 3-dimensional.

\ProofL This is the most involved case. The assumption $[\4y,\4x]=0$,
i.e, $d_1=d_2=0$ implies $\beta_1=1,$ $b_3=0$ and $b_4:=\beta_4\in
\RR$, see the table below:
$$
\diagram{[\4x,\overline\4x]&\;=\;&\4z \;+\;\;\;& a_1\4x &\;- \;
&\overline a_1\overline \4x&\;+\;& a_2\4y &\;-\;&\overline
a_2\overline \4y\cr [\4x,\overline{\4y}]&\;=\;& &\hfill \4x&+&
\hfill\overline\4x && &+ &\beta_4\overline{\4y} \cr
[\4y,\overline{\4y}]&\;=\;& & & & & &\hfill\4y & -
&\hfill\overline{\4y} \cr [\4y,\4x]&\;=\;&0\hfill \cr [\4z,\eta]&\;=
\;&\rlap{$ c_\eta \4z+ \4q_\eta$\quad $\eta\in \7l$, $\;\4q_\eta\in
\7z\cap \QQQ\;.$ }\hfill }\Leqno{VA}
$$
We need to analyze the relations between the structure constants in
more detail. Let\nline
\centerline{$\4q_{\4x}=z_1\4x+z_2\overline{\4x}+z_3\4y+z_4\overline{\4y}
\Steil{and}\4q_{\4y}=w_1\4x+w_2\overline{\4x}+w_3\4y+w_4\overline{\4y}$.}
\nline The fact that the $\4q_\eta$ commute with all elements of
$[\7l,\7l]\subset \7l^{\nil}$ (and in particular with
$\4y-\overline{\4y}$) implies $z_1=z_2$, $w_1=w_2$ and
$z_4=z_1\beta_4-z_3$, $ w_4=w_1\beta_4-w_3$.  The Jacobi identity
$\Jac{\4x}{\overline{\4x}}{\4y}=0$ yields
$$
w_j=0,\; c_{\4y}=1,
\Steil{i.e.,}[\4z,\4y]=[\4z,\overline{\4y}]=\4z,\qquad
a_1=i\alpha_1,\qquad a_2=\alpha_2+{i\over 2}\alpha_1\beta_4\;,
\Leqno{}
$$
for some $\alpha_1,\alpha_2\in \RR$, and
$\Jac{\overline{\4x}}{\4y}{\4z}=0$ implies
$$
\Im z_1=0, \qquad \Re z_3=\qu{z_1\beta_4}2 \Steil{ and }\Re
c_{\4x}=-\qu{\beta_4}2\;.\Leqno{}
$$
Summarizing,
$$
[\7l,\7l]=\CC
\4z\;\oplus\;\CC(\4x+\overline{\4x}+\hbox{\Klein$\displaystyle
{\beta_4\over 2}$} \4y+\hbox{\Klein$\displaystyle {\beta_4\over
2}$}\overline{\4y})\;\oplus\; \CC(\4y-\overline{\4y})\;, \Leqno{VK}
$$
and, using the above relations between the structure constants, a
simple check shows that this ideal is abelian. Moreover, the subset
$$
\7h:=\RR \; i\big(\4x-\overline{\4x}+\hbox{\Klein$\displaystyle
{\beta_4\over 2}$} \4y-\hbox{\Klein$\!\displaystyle {\beta_4\over
2}$}\overline{\4y}\big) \;\oplus \;\RR (\4y+\overline{\4y}) \Leqno{VL}
$$
is an abelian subalgebra of $\7g$
(and $\7g=\7h\oplus [\7g,\7g]$; in fact $\7g=\7h\ltimes [\7g,\7g]$).
\nline
It should be noted that the nilcenter $\7z$ may be 1-dimensional,
and then $\7l^{\nil}$ properly contains
$[\7l,\7l]$.

\noindent We next show that the situation `$\dim \7z=2$' does not
occur: Assume to the contrary that $\dim \7z= 2$. Then $\7z\cap \QQQ$
is nonzero. Let $\4z':=\lambda\4x + \overline \lambda\overline{\4x} +
\mu\4y +\overline \mu\overline{\4y}\in \7z^\sigma\cap \QQQ$ be
arbitrary. Either, for all such $\4z'$ the coefficient $\lambda$ is 0,
i.e., $\7z\cap \QQQ=\7z\cap \FFF$, and then this case can be ruled out
by a similar argument as in the proof of the above claim. Or, there
exists $\4z'$ with $\lambda\ne 0.$ In such a situation the identity
$[\4z',\4y-\overline{\4y}]=0$ implies $\lambda\in \RR^*$ and without
loss of generality we may suppose $\4z'=\4x + \overline{\4x} + z_3\4y
+\overline z_3\overline{\4y}$ for some $z_3\in \CC$. The condition
$[\4z,\4z']=0$ gives $\Re z_3=\qu{\beta_4}2$. But then
$\4z'=(\4x+\overline{\4x}+\qu{\beta_4}2\cd
(\4y+\overline{\4y}))+i\,\Im z_3 \cd (\4y-\overline{\4y})\in
[\7l,\7l]$, compare \Ruf{VK}. We claim, that in the situation under
consideration, $\7l^{\nil}=[\7l,\7l]$. To prove this, we simply
compute the ad-action of $\4y+\overline{\4y}$ and $\4x-\overline{\4x}$
on $[\7l,\7l]$. Since $[\4x-\overline{\4x},\4z' ]\in \CC \4z\oplus \CC
\4z'$, the relation $\beta_4^2+4(\Im z_3)^2=4\alpha_1 \Im
z_3+4\alpha_2$ must also be fulfilled. Once again, a simple
computation yields
$$
\ad(\4y+\overline{\4y})\rest {[\7l,\7l]}=-2\cd\id, \qquad
[\4x-\overline{\4x},\4z']=2\4z+2i(\alpha_1-\Im\,z_3)\4z'\;.
$$
The above identities show that for every $\4v:=u_1\cd
(\4x-\overline{\4x})+u_2(\4y+\overline{\4y})+u_3\eta,$ $\,u_j\in \CC$,
$\eta\in [\7l,\7l]$, the condition $[\4v,\4z']=0$ implies
$u_1=u_2=0$, i.e., the centralizer $C_{\7l}(\4z')$ coincides with
$[\7l,\7l]$. This proves $[\7l,\7l]=\7l^{\nil}=C_{\7l}(\4z')$. But
this is absurd, since then the nilcenter $\7z$ would coincide with
the 3-dimensional abelian ideal $[\7l,\7l]$, contrary to our
assumption `$\dim \7z=2$'.

\noindent Finally, we need to investigate the case $\dim \7z=3.$ We
claim that $\7z=[\7l,\7l]=\7l^{\nil}$: To see this, it enough to show
that $\7l^{\nil}$ is 3-dimensional, as, due to \Ruf{VK}, $[\7l,\7l]$
is 3-dimensional, too. As already mentioned (see the sentence
following \Ruf{FY}) $\7l^{\nil}$ can be at most 4-dimensional. But the
4-dimensional case can be excluded, otherwise $\7l^{\nil}=\7z\oplus
\CC \4n$ for $\4n\in \7l^{\nil}\setminus [\7l,\7l]$ which would imply
that $\7l^{\nil}$ is abelian. Hence the nilradical is
3-dimensional. This proves \ruf{VG} and Lemma \ruf{GR}.\qed

\medskip The next statement is one of the key points in our
classification of 5-dimensional 2-nondegenerate homogeneous
CR-germs. Before stating it we first fix some notation. Given a vector
space $V$, write $\aff(V)$ for the Lie algebra consisting of affine
maps of $V$. This Lie algebra (as well as the corresponding Lie group
$\Aff(V)$) has the natural semidirect product structure:
$\aff(V)=V\rtimes \gl(V)$ (with $\gl(V)=\{X\in \aff(V):X(0)=0\}$). Let
$\pi:\aff(V)\to \gl(V)$ be the projection homomorphism. We use similar
notation on the Lie group level and write, for instance,
$\pi:\Aff(V)=V\rtimes \GL(V)\to \GL(V)$ for the corresponding group
homomorphism. Sometimes we simply write $\psi^{\lin}:=\pi(\psi)$ for
the linear part of an element in $\aff(V)$ or $\Aff(V).$

\ifarx\eject\fi

\MLemma{GS} Let $(\7g,\7q)$ be a CR-algebra satisfying the fundamental
assumption \ruf{FD} and let $\7g$ be solvable and of dimension 5.
Suppose that there exists a 3-dimensional abelian ideal $\7v\subset
\7g$ and a 2-dimensional subalgebra $\7h\subset \7g$ with $\7h\cap
\7v=\{0\}$. Then the associated CR-germ $(M,o)$ is locally
CR-equivalent to a tube $F\times i\RR^3 \subset \CC^3$, where
$F\subset \RR^3$ is an affinely homogeneous surface.\Formend

\noindent The {\sl proof} is divided in several steps which give more
precise (but also more technical) information concerning the structure
of the Lie groups corresponding to $\7g$, $\7l$ and $\7q$ and a
realization of the CR-germ $(M,o)$:

\Claim{} The adjoint representation $\ad:\7h\to \gl(\7v)$ is
faithful. Consequently, identifying $\7h$ with the subalgebra
$\ad(\7h)\subset\gl(\7v)$, the Lie algebras $\7g$ and $\7l$ can be
realized as Lie subalgebras of affine transformations:
$\;\;\7g=\7v\rtimes \7h=\7v\rtimes \ad(\7h)\subset \aff(\7v)\cong
\aff(\RR^3)\;$ {and} $\;\7l=\7v^{\CC}\!\rtimes \7h^{\CC}\subset
\aff(\7v^{\CC})\cong \aff(\CC^3).$

\ProofC Let $\7n\subset \7h$ be the kernel of the adjoint
representation $\ad:\7h\to \gl(\7v).$ The case $\dim \7n=1$ can be
excluded, otherwise $\7n\oplus\7v=\7g^{\nil}=\7z$ would be
4-dimensional, contradicting Lemma \ruf{GR}. The case $\dim
\7n=2$ can be also excluded: Otherwise $\7g=\7v\times \7h$ would be
abelian or contain a 4-dimensional abelian nilradical which in both
cases would contradict \23 -- \25.\qd

Write $V\cong \RR^3$ for a vector group with Lie algebra $\7v$ and
$E:=V^{\CC}$ for its complexification. Let $H_{\GL}\subset \GL(V)$ and
$H_{\GL}^{\CC}\subset \GL(E)$ be the Lie subgroups, corresponding to
the Lie algebras $\ad(\7h)$ and $\ad(\7h^{\CC})$, respectively. Since
$\GL(V)\cong \GL(3,\RR)$ contains no compact torus of dimension
$\ge2$, each subgroup, in particular $H_{\GL}$, is closed. This is in
general not true for the complex subgroup $H_{\GL}^{\CC}$. Let
$H^{\CC}$ be the simply connected Lie group with Lie algebra
$\7h^{\CC}$, $\pr:H^{\CC}\to H_{\GL}^{\CC}\subset\GL(\7v^{\CC})$ the
homomorphism induced by $\ad:\7h^{\CC}\to
\ad(\7h^{\CC})\subset\gl(\7v^{\CC})$ and $L:=V^{\CC}\rtimes
H^{\CC}$. For simplicity, for each $h\in H^{\CC}$ we also write
$h^{\lin}\subset \GL(E)$ instead of $\pi(h)$. Let $G=V\rtimes H\subset
L$ be the real form. Since every 2-dimensional Lie algebra is
solvable, we deduce that also $\7l=\7v^{\CC}\!\rtimes \7h^{\CC}$ (as
well as $\7g$, $L$ and $G$) is solvable.

\Claim{} Let $Q\subset L$ be the subgroup corresponding to the Lie
subalgebra $\7q\subset \7l$. Then $Q$ is closed and $Q\cap
V^{\CC}=\{e\}$. Hence $L=V^{\CC}\!\rtimes Q$ is a semidirect product.

\ProofC Let $\pi:\7l\to \7h^{\CC}$ be the projection homomorphism. Our
first observation is that $\pi(\7q)=\7h^{\CC}$: The case
`$\pi(\7q)=0$' can clearly be excluded as in such a situation
$\7q\subset \7v^{\CC}$, and in turn $\7q+\sigma\7q\subset \7v^{\CC}$,
which is absurd. \nline The possibility $\dim \pi(\7q)=1$ can be
ruled out as follows: We may assume that $\pi(\7q\oplus
\sigma\7q)=\7h^{\CC}$ (otherwise $\7q\oplus \sigma\7q$ would be a
subalgebra). Either $\pi(\7f)=0,$ i.e., $\7q\cap \7v^{\CC}=\7f$, and
in turn $(\7q\oplus\sigma\7q)\cap \7v^{\CC}=\7f\oplus \sigma\7f$: This
leads to a contradiction since then $[\7F,\7H]\subset \7H\cap
\7v=\7F$, violating \25. Or, $\pi(\7f)=\pi(\7q)$. But then
$\7q=\7f\oplus \CC \4x$ with a nonzero $\4x\in \7q\cap \7v^{\CC}$, and
in turn $[\7q,\sigma\7q]\subset \7q\oplus\sigma\7q$, contradicting
\23. Summarizing, $\pi(\7q)=\7h^{\CC}$. \nline On the group level,
since $L$ is simply connected and solvable, every connected subgroup
is closed. The restriction of $\pi$ to $Q$ induces a surjective
homomorphism $Q\to H^{\CC}$. Since both groups are 2-dimensional,
this homomorphism is a covering. Our assumption that $H^{\CC}$ is
simply connected finally implies that
$\pi\lower1pt\hbox{\mathsurround=0pt$|$}_{Q}:Q\to H^{\CC}$ is an
isomorphism. In particular $Q\cap V^{\CC}=Q\cap \ker\pi=\{e\}$.\qd

\Claim{VM}With respect to the identification $Z:=L/Q=V^{\CC}\cong
\CC^3$ the real form $G$ acts on $L/Q$ by affine transformations and
$V\subset G$ by translations.

\ProofC The existence of the decomposition $L=V^{\CC}\!\rtimes Q$
implies that there are well-defined functions $v\colon L\to V^{\CC}$,
$q\colon L\to Q$ such that $\ell=v(\ell)\cd q(\ell)$ for every
$\ell\in L.$ Let $g=w\cd h\in V\rtimes H=G$ (with $w\in V, h\in H$) be
arbitrary. Then, for any $z\in V^{\CC}$ we have
$$
g\cd zQ=w\cd h \cd zQ =w \cd v(h)\cd q(h)\cd zQ = w \cd v(h)\cd
(q(h)\cd z \cd q(h)^{-1})Q
$$
and $ q(h)\cd z \cd q(h)^{-1}= v(h)^{-1}h \cd z \cd h^{-1}v(h)=
h^{\lin}(z)$. Hence, with respect to the identification $V^{\CC}=L/Q$
(induced by the inclusion $V^{\CC}\into L$ such that $0$ corresponds
to the point $eQ\in L/Q$), the action of $G$ can written as follows:
$$
g\cdot z= h^{\lin}(z) + v(h) +w \qquad g=w\cd h\in L,\;\; z\in V^{\CC}\,.
\Leqno{VN}
$$
In particular, the subgroup $V\subset G$ acts by translations
$z\mapsto z+w$. \qd

\noindent
Consequently $M:= G\cdot 0=VH\cdot 0= V +
F\subset V\oplus iV$ with $F:=M\cap iV$.
It should be noted, however, that in general $F :=(G\cdot
0)\cap iV \ne H\cdot 0$. Nevertheless, as we shortly will see, $F$ is
affinely homogeneous under a slightly different subgroup of
$\Aff(iV)$. Clearly, multiplying a tube manifold $F+ iV\subset
V\oplus iV$, $F\subset V$, by the imaginary unit $i$ we get the
CR-equivalent realization $V+ iF=V+ F'$ with $F'=iF\subset
iV$. The latter form of a tube manifold is more suitable in our
particular setup, and we keep this notation until the end of the
proof of the Main Lemma.

\Claim{} Retaining the previous notation, there exists a subgroup
$B\subset iV\!\rtimes \GL(iV)=\Aff(iV)$ such that $F:= (G\cdot 0)\cap
iV =B\cdot 0$.

\ProofC Let $\pr^{i}:V\oplus iV\to iV$ be the linear projection. A
glance at \Ruf{VN} shows that $F=\pr^{i}\{v(h): h\in H\}.$ In order to
determine $v(h)$ more explicitly we need to analyze the position of
$Q$ in $V^{\CC}\!\rtimes H^{\CC}$ in greater detail: Since $\7h$ is
2-dimensional, there exists a basis $\4s,\4n\in \7h$ such that
$[\4s,\4n]=\epsilon \4n$ with $\epsilon\in \{0,1\}$. Recall that the
projection map $\pi:\7q\to \7h^{\CC}$ is an isomorphism. Consequently
there exist $\4w_{\4s},\4w_{\4n}\in \7v^{\CC}=V^{\CC}$ such that the
elements $\4w_{\4s}+\4s$ and $\4w_{\4n}+\4n$ in
$\7l=V^{\CC}\oplus\7h^{\CC}$ generate $\7q$. Since then $\sigma\7q=\CC
(\overline {\4w}_{\4s}+\4s)\oplus \CC (\overline{\4w}_{\4n}+\4n)$, we
must have $\4w_{\4s}\ne \overline {\4w}_{\4s}$ and $\4w_{\4n}\ne
\overline {\4w}_{\4n}$ (otherwise $\7q\cap\sigma\7q \ne 0$). Let
$\exp:\7l\to L$ and $\Exp:\ad(\7h^{\CC})\to \GL(V^{\CC})$ be the
exponential maps (i.e., $\Exp({\ad \4v})=\pi(\exp \4v)$ with $\pi$ as
in the paragraph preceeding \ruf{GS}). Furthermore, let $\Psi$ be the
entire function defined by $$\Psi(z)={e^{z}-1\over
z}=\sum^{\infty}_{k=0}{z^{k}\over(k+1)!}z^{k}\,.$$
Then for $\4S^t:=\Psi(t\ad(\4s))$ and $\4N^u:=\Psi(u\ad(\4n))$
a simple computation shows:
$$
\eqalign{ H&=\{\exp(t\4s)\cd \exp(u\4n):t,u\in \RR\} \cr Q&=\{\exp
t(\4w_{\4s}+\4s)\cdot\exp u(\4w_{\4n}+\4n): t,u\in \CC \}\cr &=\{
\4S^t(t\4w_{\4s})\cd\exp(t\4s)\cdot
\4N^u(u\4w_{\4n})\cd\exp(u\4n):t,u\in \CC \} \cr
&=\{\big(\4S^t(t\4w_{\4s})\cd \Exp(t\ad(\4s))\,(\4N^u(u\4w_{\4n}))\big)
\cdot \exp (t\4s)\cd\exp(u\4n):t,u\in \CC \} \;\;\subset\;\; V^{\CC}
\cdot H^{\CC}=L\,. } \Leqno{VO}
$$
The explicit form of $v(h)$ (compare the proof of Claim \ruf{VM}) can
be read off the last line in \Ruf{VO}:
$$
v(h)=v(\exp (t\4s)\exp(u\4n))
=\big(\4S^t(t\4w_{\4s})\big)^{-1}\cdot\big(\Exp(t\ad(\4s))\,
(\4N^u(u\4w_{\4n}))\big)^{-1}\,.
$$
Since $\ad(\4s), \4N^u$ and $\4S^t$ are real operators, it follows for
$h=\exp(t\4s)\cdot\exp( u\4n)$ as before:
$$
\eqalign{ \pr^{i}(v(h))\;&=\;\big(\4S^t(t\4w_{\4s})\big)^{-1}\!\cdot
\big(\Exp(t\ad(\4s))\, (\4N^u(u\4w_{\4n}))\big)^{-1}\cr&=\;\; \exp
t(-\4w^{i}_{\4s}+\4s)\cd \exp u(-\4w^{i}_{\4n}+\4n)\cdot 0 \;\subset
iV\;.\cr }\Leqno{VP}
$$
(Using additive notation, $\pr^{i}(v(h))=-\4S^t(t\4w^{i}_{\4s})-
\Exp(t\ad(\4s))\,(\4N^u(u\4w^{i}_{\4n}))\;\subset iV\,$). Define

\hfil $\7b:=\RR(-\4w^{i}_{\4s}+\4s)\oplus
\RR(-\4w^{i}_{\4n}+\4n)\;\;\subset \; \7l=V^{\CC}\!\rtimes
\7h^{\CC}$\hfil

\noindent and check that this is a Lie algebra. Then $B:=\exp
\RR(-\4w^{i}_{\4s}+\4s)\cd \exp \RR(-\4w^{i}_{\4n}+\4n)$ is the
subgroup of $L$ with Lie algebra $\7b$, and \Ruf{VP} shows that
$F=\pr^{i}\{v(h):h\in H\}=B\cdot 0$. This finishes the proof of the
claim and of the Main Lemma.\qed

\KAP{Sechzehn}{The final steps}

Our final step toward the complete classification of all
5-dimensional 2-nondegenerate and homogeneous CR-germs is to deduce
that each 5-dimensional solvable Lie algebra $\7g$ which occurs in a
CR-algebra $(\7g,\7q)$ subject to \ruf{FD}, also satisfies the
assumptions of the preceeding Main Lemma \ruf{GS}:

\Lemma{VQ} Let $(\7g,\7q)$ be a CR-algebra satisfying \ruf{FD} and
suppose that $\7g$ is solvable and of dimension 5.  Then there exists a
semidirect product decomposition $\7g=\7v\rtimes \7h$ with a
3-dimensional abelian ideal $\7v\subset \7g$ and a 2-dimensional
subalgebra $\7h$.  \Formend

\Proof We have already observed in Lemma \ruf{GR} that the nilcenter
$\7z$ has dimension 1 or 3. If $\dim \7z=3$ then, due to \Ruf{VF}, we
can apply \ruf{VG}: Consequently, we can choose $\7v=[\7g,\7g]$
(compare \Ruf{VK}) and $\7h$ as defined in \Ruf{VL}.

\noindent The situation $\dim \7z=1$ requires some more elaborated
work. We classify all CR-algebras $(\7g,\7q)$ under consideration in
terms of the corresponding structure equations with respect to some
perfect basis $\4x,...\,,\4z$ of $\7l$ (see \ruf{GQ}): Given
$(\7g, \7q)$, let the corresponding structure equations be as in
\Ruf{FI}, taking into account \Ruf{FZ}. To handle the various sets of
relations between the structure constants, we divide the class of
CR-algebras under consideration into the subclasses A, B and C, see
below.

\noindent {\fett Case A: $\beta_1\ne \pm 1$}. In this situation it is
possible to assume $a_1=0$ (simply replace $\4x$ by
$\4x+\lambda\4y$ with $\lambda=u+iv$ defined by $u:=\qu{\Re
a_1}{(1-\beta_1)}$ and $v:=\qu{\Im a_1}{(1+\beta_1)}$). The structure
equations then read
$$
\diagram{[\4x,\overline\4x]&\;=\;&\4z+\;& & & &\;+\;& a_2\4y &\;-\;
&\overline a_2\overline \4y \cr [\4x,\overline{\4y}]&\;=\;& &\kern1em
\beta_1\4x&+& \;\overline\4x &+&b_3\4y&+&b_4\overline{\4y} \cr
[\4y,\overline{\4y}]&\;=\;& & & & & & \;\;\4y & - & \;\;\overline{\4y}
\cr [\4y,\4x]&\;=\;& & \kern-10pt(1-\beta_1)\4x & & &+ &d_2\4y & &
\rlap{\qquad\qquad $d_2=\overline b_3+(1-\beta_1)b_4$} \cr
[\4z,\eta]&\;= \;&c_\eta\4z \rlap{ \qquad for every $\eta\in \7l\;$
,}}\qquad\qquad\qquad\qquad \Leqno{GT}
$$
and we now work out further constraints imposed on the constants: An
explicit evaluation of the Jacobi identity
$\Jac{\4x}{\overline{\4x}}{\4y}=0$ yields $c_\4y=c_{\overline{\4y}}
=2\beta_1-1.$ Furthermore, we obtain the equations $\overline
b_3(1+\beta_1)=b_4(\beta_1^2-\beta^{}_2)$ and $\overline
b_3(1+\beta_1)=(b_4-\overline b_4)(\beta_1-1)$ which imply
$b_4(1-\beta_1)=\overline b_4.$ \nline In order to investigate most
conveniently further relations between the structure constants, we
deal separately with the following 3 subcases:

\medskip\centerline
{ {\bf AI:} $b_4\in i\RR^{*}$ and $\beta_1=2,$\hfil
{\bf AII:} $b_4\in \RR^{*}$ and $\beta_1=0,$\hfil
{\bf AIII:} $b_4=0.$}

\noindent {\bf Ad AI:} Put $\beta_{4}:=-ib_{4}$. In this particular
situation the identity $\Jac{\4x}{\overline{\4x}}{\4y}=0$ implies
$b_3=-{2\over 3} i\beta_4,$ $a_2=-{2\over 9}\beta_4^2$ and
$\Jac{\4x}[\4y]{\4z}=0$ implies $c_\4x=-i\beta_4.$ There are no more
conditions imposed by the Jacobi identity and the structure
equations of $\7l$ are now
$$
\diagram{[\4x,\overline\4x]&\;=\;&\4z+\;& & & &\;-\;& 2 \gamma^2\4y
&\;+\;&2 \gamma^2\overline{\4y} \cr [\4x,\overline{\4y}]&\;=\;&
&\kern1em 2\4x&+& \; \overline\4x &-&2i \gamma\4y&+&3i
\gamma\overline{\4y} \cr [\4y, \overline{\4y}]&\;=\;& & & & & &
\;\;\4y & - & \;\;\overline{\4y} \cr [\4y,\4x]&\;=\;& & \kern-3pt-\4x
& & &- &i \gamma \4y & & \cr [\4z,\4x]&\;= \;&\rlap{$-3i \gamma\,\4z$
\qquad\qquad $[\4z,\4y]=-3\4z\;,$}\hfill }\qquad\qquad\qquad\qquad
\Leqno{GU}
$$ where $ \gamma:=\beta_4/3\in \RR^{*}$. Keeping in mind Lemma
\ruf{GS}, the structure of $\7l$, $\7g$ and the position of $\7q$ is
determined by \Ruf{GU}. This can be seen more clearly by decomposing
$\7g$ and $\7l$ into the eigenspaces of $\ad(\4s),$ where
$\4s:=-(\4y+\overline{\4y}/2:$ Define the following elements
$\4n,\4v_1,\4v_2,\4v_3$ from $\7g:$
$$
\diagram{\4n&\;:= \;& i\4x-i\overline{\4x}-\;\gamma\4y-\;
\gamma\overline{\4y} \cr
\4v_1&:=& \qquad\qquad \hbox{\klein$\displaystyle{i\over 2}$}
\4y-\hbox{\klein$\displaystyle{i\over 2}$}\overline{\4y} \cr
\4v_2&:=& \;\4x+\;\overline{\4x}-2i\gamma\4y+2i \gamma\overline{\4y}\cr
\4v_3&:=& 2i\4z\;.\vphantom{{}^{\textstyle|}}\hfill}
$$
It is clear that $\4s,\4n,\4v_1,\4v_2,\4v_3$ form a basis of $\7g.$
The bracket relations are:
$$
[\4s,\4v_k]=k\4v_k,\;[\4s,\4n]=\4n,\quad [\4n,\4v_1]=\4v_2,\;
[\4n,\4v_2]=\4v_3,\;[\4n,\4v_3]=0.
$$ Further, $\7v:=\RR\4v_1\oplus \RR\4v_2\oplus \RR\4v_3\cong \RR^3$
is an abelian ideal in $\7g$ with $[\7g,\7g]=\7g^{\nil}=\RR\4n\oplus
\7v$. Hence, $\7g$ has the structure of the semidirect product
$\7v\rtimes\7h$ with $\7h=\RR\4s\oplus\RR\4n$ and the Main Lemma
applies.

\noindent {\bf Remark.} A direct verification shows that for every
$\gamma\in\RR^*$ $(\7g,\7q)$ in \Ruf{GU} is associated to Example
\ruf{EV}.

\noindent We show that in the next case the Lie algebra cannot be
solvable, hence, this case can be discarded.

\noindent {\bf Ad AII:} Write $\beta_4:=b_4$.  The Jacobi identity
implies $a_2=0=b_3$ and $d_2=\beta_4$, and \Ruf{VA} reads \vskip-18pt
$$
\diagram{[\4x,\overline\4x]&\;=\;&\4z\;\;\cr
[\4x,\overline{\4y}]&\;=\;& \;\;& \kern1em & \kern1em \;&
\;\overline\4x & \;& \;&+&\beta_4\overline{\4y} \cr
[\4y,\overline{\4y}]&\;=\;& & & & & & \kern.7em \;\4y & - & \kern.7em
\; \overline{\4y} \cr [\4y,\4x]&\;=\;& & \;\4x & & &+ &\beta_4 \4y & &
\cr [\4z,\4x]&\;= \; &\rlap{$\beta_4\4z$ \qquad\qquad
$[\4z,\4y]=-\4z\;.$}\hfill }\qquad\qquad\qquad\qquad \Leqno{}
$$
But then the linear span of the vectors
$$
\displaylines{ \4e^+:=\4y-\overline{\4y}, \qquad\qquad
\4h:=-\frac1{\beta_4}\big (\4x+\overline{\4x}+(\beta_4-1)\4y \;+
\;(\beta_4+1) \overline{\4y}\big),\kern1em \cr \4e^-:=
\frac1{4\beta_4^2}\big(2\4z+(2-2\beta_4)\4x+(2+2\beta_4)
\overline{\4x}-(1-\beta_4)^2\4y+(1+\beta_4)^2\overline{\4y}\big) }
$$
is a copy of $\7{sl}(2,\CC)$ in $\7l$, that is, $\7l$ is not solvable.

\noindent{\bf Ad AIII:}. The condition $b_4=0$ together with the Jacobi
identity $\Jac{\4x}{\overline{\4x}}{\4y}=0$ implies $a_2=b_3=0$ and
$c_\4y=2\beta_1-1$. Since $(1-\beta_1)\ne 0$, the identity
$[\4z,[\4y,\4x]]=0$ implies $[\4z,\4x]=0$, see the table below.
$$
\diagram{[\4x,\overline\4x]&\;=\;&\4z\;& \cr
[\4x,\overline{\4y}]&\;=\;& &\kern1em \beta_1\4x&+& \;\overline\4x &
\cr
[\4y,\overline{\4y}]&\;=\;& & & & & & \;\;\4y & - & \;\;\overline{\4y}
\cr
[\4y,\4x]&\;=\;& & \kern-6pt(1-\beta_1)\4x & & & & & &
\cr
[\4z,\4x]&\;= \;&0\hfill
\rlap{\qquad,\qquad $[\4z,\4y] \;=\; (2\beta_1-1)\4z\;.$} \hfill
}\qquad\qquad\qquad\qquad
\Leqno{GV}
$$
Select the following basis of $\7g$:
$$
\diagram{\4n\,&\;:= \;& \hbox{\klein$\displaystyle{i\over
2}$}(\4x-\overline{\4x}) \qquad,\qquad\quad \4s\,:=
\hbox{\Klein$\displaystyle{1\over 2-2\beta_1}$}(\4y+\overline{\4y})\cr
\4v_1&:=& \qquad\qquad i\4y\;- i\overline{\4y} \hfill\cr \4v_2&:=&
\;\;\4x \; + \;\overline{\4x} \hfill\cr \4v_3&:=&
i\4z\;.\vphantom{\big|}\hfill}
$$
One checks immediately that $\7v:=\RR \4v_1\oplus\RR \4v_2\oplus\RR
\4v_3$ is an abelian ideal and $\7h:= \RR \4s\oplus\RR \4n$ is a
subalgebra with $[\4s,\4n]=\4n$. Further,
$$
[\4s,\4v_j]= \hbox{\klein$\displaystyle{2-j+(j-1)\beta_1\over
\beta_1-1}$}\4v_j\steil{for}j=1,2,3\Steil{and}[\4n,\4v_1]=\4v_2,\;
[\4n,\4v_2]=\4v_3,\;[\4n,\4v_3]=0.
$$
Hence, $\7g=\7v\rtimes \7h$ as claimed, and the Main Lemma applies.
Also in this case for all $\beta_1\ne \pm 1$ the CR-algebra
$(\7g,\7q)$ is associated to Example \ruf{EV}.

\medskip It remains to discuss the cases $\beta_1=\pm1.$

\noindent {\fett Case B: $\beta_1=1$.} Plugging $\beta_1$ into
\Ruf{FZ} gives $d_2=\overline b_3$ and $d_1=0.$ A direct check shows
that $\Jac{\4x}{\overline{\4x}}{\4y}=0$ implies $0=b_3=d_2,$ i.e.,
$[\4x,\4y]=0$. Lemma \ruf{VG} then gives that $\7g$
is isomorphic to the semidirect product $\7v\rtimes\7h$ with
$\7v=[\7g,\7g]$ and the abelian subalgebra $\7h$ as defined in
\Ruf{VL}.

\medskip\noindent {\fett Case C: $\beta_1=-1$.} We proceed as in the
preceeding cases. Starting from the structure equations \Ruf{FI} with
respect to some perfect basis
$\4x,\overline{\4x},\4y,\overline{\4y},\4z\,$, we first evaluate
\Ruf{FZ} for this particular value of $\beta_1$. We get $d_1=2$ and
$d_2=\overline b_3+2b_4.$ Next, the Jacobi identity
$\Jac{\4x}{\overline{\4x}}{\4y}=0$ implies $b_4:=\beta_4\in \RR$ and
$a_1=\beta_4.$ Further, $\Jac{\4x}{\4y}{\overline{\4y}}=0$ implies
$b_3=-\beta_4 +i\gamma,$ $\gamma\in \RR.$ Next,
$\Jac{\4x}{\overline{\4x}}{\4y}=0$ determines the value of $a_2:$
$a_2=(3\beta_4^2+\gamma^2)/4-i\beta_4\gamma/2.$ Finally
$\Jac{\4x}{\4y}{\4z}=0$ implies
$[\4z,\4x]=(\quu{3\beta_4}2-\quu{3i\gamma}2)\4z,$ see diagram below:
$$
\diagram{[\4x,\overline\4x]&\;=\;&\4z+\;& \beta_4\4x & -
&\beta_4\overline{\4x} &\;+\;& a_2\4y &\;-\;&\overline a_2\overline
\4y \rlap{\qquad
$a_2=\quu{(3\beta_4^2+\gamma^2)}4-\quu{i\gamma\beta_4}2 $} \cr
[\4x,\overline{\4y}]&\;=\;& &\kern1em -\4x&+& \;\overline\4x &
&+b_3\4y&+&\beta_4\overline{\4y} \rlap{\qquad $b_3=-\beta_4+i\gamma$}
\cr [\4y,\overline{\4y}]&\;=\;& & & & & & \;\;\4y & - &
\;\;\overline{\4y} \cr [\4y,\4x]&\;=\; & & \kern-4pt 2\4x & & &-
&b_4\4y & & \cr [\4z,\4x]&\;= \;\rlap{$
(\quu{3\beta_4}2-\quu{3i\gamma}2) \4z$ \qquad $[\4z,\4y]=-3\4z\;\;$.}
}\qquad\qquad\qquad\qquad\qquad\qquad\Leqno{GX}
$$
Select the following basis of $\7g$:
$$
\diagram{\4n\;&\;:= \;& 2i\4x-2i\overline{\4x} \;-ib_3\4y
\;+\;i\overline b_3\overline{\4y} \qquad\qquad \4s\;:=
-\hbox{\klein$\displaystyle{1\over 4}$}(\4y+\overline{\4y})\cr
\4v_1&:=& \qquad\qquad\kern1.3em \hbox{\Klein$\displaystyle{i\over
2}$}\4y\;\;- \hbox{\Klein$\displaystyle{i\over 2}$}\overline{\4y}
\hfill\cr \4v_2&:=& \;2\4x + 2\overline{\4x} -\;b_3\4y \;- \;\;\;
\overline b_3 \overline{\4y} \hfill\cr \4v_3&:=&
4i\4z\;.\vphantom{\big|}\hfill}
$$
A direct computation (using \Ruf{GX}) shows that $\7v:=\RR\4v_1\oplus
\RR\4v_2\oplus \RR\4v_3$ is an abelian ideal and
$$
[\4s,\4n]=\4n, \qquad \diagram{[\4s,\4v_1]=
-\hbox{\klein$\displaystyle{1\over 2}$}\4v_1 , \quad[\4s,\4v_2]=
\hbox{\klein$\displaystyle{1\over 2}$}\4v_2 , \quad[\4s,\4v_3]=
\hbox{\klein$\displaystyle{3\over 2}$}\4v_3 , \cr \kern.6em
[\4n,\4v_1]=\;\4v_2 , \kern2.4em [\4n,\4v_2]= \;\4v_3 , \kern1.6em
[\4n,\4v_3]= \;0\;.\qquad }
$$ Again, this shows that $\7g$ is isomorphic to the semidirect
product $\7v\rtimes \7h$ with $\7h=\RR \4n\oplus \RR \4s$. Actually,
an explicit realization of the corresponding CR-manifold $M$ along the
lines of proof of the Main Lemma shows that $(\7g,\7q)$ is associated
to the tube over the future light cone. \qed

\medskip We close by stating that Lemma \ruf{GS} together with Lemma
\ruf{VQ} finishes the proof of the classification theorem.

\Partskip
{\gross\noindent References\write\lst{\def\csname Refa\endcsname{\folio}}}
\bigskip

{\Klein
\parindent 15pt\advance\parskip-1pt

\def\Springer{Ber\-lin-Hei\-del\-berg-New York: Sprin\-ger~}

\font\klCC=txmia scaled 900

\Ref{ANFR}Andreotti, A., Fredricks, G.A.: Embeddability of real analytic Cauchy-Riemann manifolds. Ann. Scuola Norm. Sup. Pisa Cl. Sci. {\bf 6} (1979), 285--304.
\Ref{ANHI}Andreotti, A., Hill, C.D.: Complex characteristic coordinates and tangential Cauchy-Riemann equations. Ann. Scuola Norm. Sup. Pisa Sci. Fis. {\bf 26} (1972), 299-324.
\Ref{BELO}Beloshapka, V.K.: Symmetries of Real Hypersurfaces in Complex 3-Space. Math. Notes {\bf 78} (2005), 156--163.
\Ref{BHUR}Baouendi, M.S., Huang, X., Rothschild, L.P.: Regularity of CR mappings between algebraic hypersurfaces. Invent. Math. {\bf 125} (1996), 13-36.
\Ref{BERO}Baouendi, M.S., Ebenfelt, P., Rothschild, L.P.: {\sl Real Submanifolds in Complex Spaces and Their Mappings}. Princeton Math. Series {\bf 47}, Princeton Univ. Press, 1998.
\Ref{BAHR}Baouendi, M.S., Huang, X., Rothschild, L.P.: Regularity of CR mappings between algebraic hypersurfaces. Invent. Math. {\bf 125} (1996), 13-36.
\Ref{BARZ} Baouendi, M.S., Rothschild, L.P., Zaitsev, D.: Equivalences of real submanifolds in complex space. J. Differential Geom. 59 (2001), 301--351.
\Ref{BOGG}Boggess, A.: {\sl CR manifolds and the tangential Cauchy-Riemann complex.} Studies in Advanced Mathematics. CRC Press. Boca Raton~ Ann Arbor~ Boston~ London 1991.
\Ref{BOUR}Bourbaki, N.: {\sl \'El\'ements de math\'ematique: groupes et alg\`ebres de Lie.} Chapitre I-IX. Herman, Paris 1960.
\Ref{BUSH}Burns, D., Shnider, S.: Spherical hypersurfaces in complex manifolds. Invent. Math. {\bf 33} (1976), 223--246.
\Ref{CART}Cartan, \'E.: Sur la g\'eom\'etrie pseudo-conforme des hypersurfaces de l'espace de deux variables complexes. Annali di Matematica Pura ed Applicata {\bf 11/1} (1933) 17--90.
\Ref{CHMO}Chern, S.S., Moser, J.K.: Real hypersurfaces in complex manifolds. Acta. Math. {\bf 133} (1974), 219-271.
\Ref{DAYA}Dadok, J.,Yang, P.: Automorphisms of tube domains and spherical hypersurfaces. Amer. J. Math. {\bf 107} (1985), 999-1013.
\Ref{DOKR}Doubrov, B., Komrakov, B., Rabinovich, M.: Homogeneous surfaces in the three-dimenaional affine geometry. In: {\sl Geometry and Topology of Submanifolds, VIII,} World Scientific, Singapore 1996, pp. 168-178
\Ref{EAEZ}Eastwood, M., Ezhov, V.: On Affine Normal Forms and a Classification of Homogeneous Surfaces in Affine Three-Space. Geometria Dedicata {\bf 77} (1999), 11-69.
\Ref{EBFT}Ebenfelt, P.: Normal Forms and Biholomorphic Equivalence of Real Hypersurfaces in \hbox{\klCC\char131}$^3$. Indiana J. Math. {\bf 47} (1998), 311-366.
\Ref{EBEN}Ebenfelt, P.: Uniformly Levi degenerate CR manifolds: the 5-dimensional case. Duke Math. J. {\bf 110} (2001), 37--80.
\Ref{FELS}Fels, G.: Locally homogeneous finitely nondegenerate CR-manifolds. To appear. {\tt arXiv:math.CV/0606032}
\Ref{FEKA}Fels, G., Kaup, W.: CR-manifolds of dimension 5: A Lie algebra approach. J. Reine Angew. Math, to appear. {\ninett http://arxiv.org/pdf/math.DS/0508011}
\Ref{GAME}Gaussier, H., Merker, J.: A new example of a uniformly Levi degenerate hypersurface in \hbox{\klCC\char131}$^3$. Ark. Mat. {\bf 41} (2003), 85--94.
\Ref{HOCH}Hochschild, G.: {\sl The structure of Lie groups.} Holden-Day, Inc., San Francisco-London-Amsterdam 1965
\Ref{KNAP}Knapp, A.W.: {\it Lie groups beyond an introduction}. Second edition. Progress in Mathematics {\bf 140}. Birkh\"auser Boston, Inc., Boston, MA, 2002.
\Ref{HART}Hartshorne, R.: {\sl Algebraic Geometry}. Graduate Texts in Mathematics {\bf 52}. \Springer 1977
\Ref{HORM}H\"ormander, L.: {\sl An Introduction to Complex Analysis in Several Variables}. Third edition. North-Holland Publishing Co. Amsterdam - New York, 1990.
\Ref{ISMI}Isaev, A.V., Mishchenko, M.A.: Classification of spherical tube hypersurfaces that have one minus in the Levi signature form. Math. USSR-Izv. {\bf 33} (1989), 441-472.
\Ref{KAZT}Kaup, W., Zaitsev, D.: On local CR-transformations of Levi degenerate group orbits in compact Hermitian symmetric spaces. J. Eur. Math. Soc. {\bf 8} (2006), 465-490.
\Ref{LOBO}Loboda, A.V.: Homogeneous real hypersurfaces in \hbox{\klCC\char131}$^3$ with $2$-dimensional isotropy groups. Tr. Mat. Inst. Steklova {\bf 235} (2001), 114-142 and Proc. Steklov Inst. Math {\bf 235} (2001), 107-135.
\Ref{LOBP}Loboda, A.V.: Homogeneous nondegenerate surfaces in \hbox{\klCC\char131}$^{3}$ with two-dimensional isotropy groups. (Russian) translation in  Funct. Anal. Appl. {\bf 36}  (2002), 151--153.
\Ref{LOBQ}Loboda, A.V.: On the determination of a homogeneous strictly pseudoconvex hypersurface from the coefficients of its normal equation. (Russian)  translation in  Math. Notes  {\bf 73}  (2003), 419--423.
\Ref{MENA}Medori, C., Nacinovich, M.: Algebras of infinitesimal CR automorphisms.  J. Algebra  {\bf287}  (2005), 234--274.
\Ref{PALA}{Palais, R.~S.:} {\sl A global formulation of the Lie theory of transformation groups.} Mem. AMS 1957.
\Ref{SAKA}Sakai, T.: {\sl Riemannian geometry.} Translations of Mathematical Monographs, 149. American Mathematical Society, Providence, RI, 1996.
\Ref{SEVL}Sergeev, A.G., Vladimirov, V.S.: Complex analysis in the future tube. In: {\sl Several Complex Variables II,} Encyclopedia Math. Sci {\bf 8}, \Springer 1994.
\Ref{TANA}Tanaka, N.: On the pseudo-conformal geometry of hypersurfaces of the space of $n$ complex variables. J. Math. Soc. Japan {\bf 14} (1962), 397-429.
\bigskip
}

\medskip\centerline{Mathematisches Institut, Universit\"at
T\"ubingen, Auf der Morgenstelle 10, 72076 T\"ubingen, Germany}
\smallskip
\centerline{e-mail: {\ninett gfels@uni-tuebingen.de, kaup@uni-tuebingen.de}}

\closeout\aux\closeout\lst\bye